\documentclass[11pt,a4paper]{article}
\usepackage[francais,english]{babel}
\usepackage{amsmath}
\usepackage{amssymb}
\usepackage{latexsym}
\usepackage{amsfonts}
\newcommand{\ghg}{{\widehat{{\mathfrak g}}}}
\newcommand{\np}{{\widehat{{\mathfrak n}}_{+}}}
\newcommand{\bm}{{\widehat{{\mathfrak b}}_{-}}}
\newcommand{\nm}{{\widehat{{\mathfrak n}}_{-}}}
\newcommand{\NP}{{\widehat N_{+}}}
\newcommand{\BM}{{\widehat B_{-}}}
\newcommand{\BP}{{\widehat B_{+}}}
\newcommand{\GHG}{{\widehat G}}
\newcommand{\NN}{{\mathbb N}}
\newcommand{\EE}{{\mathbb E}}
\newcommand{\CC}{{\mathbb C}}
\newcommand{\ZZ}{{\mathbb Z}}
\newcommand{\bs}{\bigskip}
\newcommand{\q}[1]{\left[ #1\right] }

\renewcommand{\deg}[1]{{\partial}^{o}#1}
\newcommand{\degp}[1]{{\partial}_{p}^{o}#1}
\newcommand{\dega}[1]{{\partial}_{1}^{o}#1}
\newcommand{\degb}[1]{{\partial}_{2}^{o}#1}
\newcommand{\degi}[1]{{\partial}_{i}^{o}#1}
\newcommand{\dego}[1]{{\partial'}^{o}#1}
\newcommand{\ootimes}{{\bar\otimes}}
\newcommand{\DDelta}{{\bar\Delta}}

\renewcommand{\a}{\alpha}

\renewcommand{\b}{\beta}

\newcommand{\C}[2]{\left[ \begin{array}{c}#1\\#2\end{array}\right] }
\newcommand{\cnp}[2]{\left( \begin{array}{c}#1\\#2\end{array}\right) }

\newcommand{\tab}{\hspace*{\fill}}
\newcommand{\ok}{\tab $\Box$\bs}

\newtheorem{lem}{Lemme}[section]
\newtheorem{prop}{Proposition}[section]
\newtheorem{thm}{Th{\'e}or{\`e}me}[section]
\newtheorem{cor}{Corollaire}[section]

\newenvironment{dem}{{\noindent\bf D{\'e}monstration. }}{\ok}
\newenvironment{pr}{{\noindent\bf Preuve. }}{}

\title{Th{\'e}orie de Toda sur r{\'e}seau et espaces homog{\`e}nes
quantiques}
\author{C.Grunspan}
\date{}
\numberwithin{equation}{section}

\begin{document}
\maketitle

\newpage
\selectlanguage{french}
\begin{abstract}
La th{\'e}orie de Toda continue pour $\widehat{sl_2}$
est la consid{\'e}ration de deux op{\'e}rateurs de Vertex
$V_{\pm}(z)=
{\sideset{_o^o}{}\exp\bigl(\pm\beta\varphi (z)\bigr)}_{o}^{o}$
sur des modules de Fock,
o{\`u} $\varphi (z)$ est un champ libre.
Ces op{\'e}rateurs v{\'e}rifient pour $|z|=1$ 
et $q=\exp(i\pi{\beta}^{2})$
des relations de $q-$commutation remarquables.
En outre, leurs modes z{\'e}ro,
$V_{\pm}=\oint V_{\pm}(z)dz$ satisfont les relations de Serre
quantiques ([1]).
Ceci conduit {\`a} l'{\'e}tude 
du syst{\`e}me sur r{\'e}seau correspondant dont la limite classique
a {\'e}t{\'e} trait{\'e}e dans [2].
Succintement, on consid{\`e}re des op{\'e}rateurs $x_i$ et $y_i$ 
soumis {\`a} des lois de
$q$-commutation.
Les sommes $\Sigma^{+}=\sum x_i$ et $\Sigma^{-}=\sum y_i$
sont les analogues discrets de $V_{+}$ et $V_{-}$.
Sur l'alg{\`e}bre engendr{\'e}e par les variables $x_i$ et $y_i$, 
on a une action de
$U_q(\bm)$, o{\`u} $\bm$ est la sous-alg{\`e}bre
de Borel n{\'e}gative de $\widehat{sl_2}$.
On obtient ainsi 
une structure de $U_q(\bm)$-alg{\`e}bre-module.
Un des probl{\`e}mes qui se pose, est de caract{\'e}riser 
cette alg{\`e}bre-module en termes
de $U_q(\bm)$, ou, du moins de caract{\'e}riser 
au moyen du groupe quantique $U_q(\bm)$ uniquement,
une sous $U_q(\bm)$-alg{\`e}bre-module de 
$\CC[x_i^{\pm 1},y_i^{\pm 1},1\leq i\leq n]_q$
qui soit int{\`e}gre et qui poss{\`e}de le m{\^e}me corps 
de fractions que
$\CC[x_i^{\pm 1},y_i^{\pm 1},1\leq i\leq n]_q$.
A isomorphisme pr{\`e}s, nous obtiendrons 
cette sous-alg{\`e}bre module comme un espace homog{\`e}ne
quantique sur $\BM$.
Les r{\'e}sultats d{\'e}montr{\'e}s g{\'e}n{\'e}ralisent ceux de [2].
\end{abstract}

\selectlanguage{english}
\begin{abstract}
The continuous Toda theory for $\widehat{sl_2}$ is the consideration
of two vertex operators $V_{\pm}(z)=
{\sideset{_o^o}{}\exp\bigl(\pm\beta\varphi (z)\bigr)}_{o}^{o}$
on Fock modules where $\varphi (z)$ is a free field. These operators
satisfy for $|z|=1$ and $q=\exp(i\pi{\beta}^{2})$ some remarkable
$q-$commutations relations. Furthermore, their zero modes,
$V_{\pm}=\oint V_{\pm}(z)dz$ satisfy the quantum Serre relations
([1]). This leads to the study of the corresponding system on lattice
whose classical limit has been treated in [2].
Briefly, we consider some operators $x_i$ and $y_i$ which obey to some
laws of $q-$ commutations. The sums $\Sigma^{+}=\sum x_i$ and
$\Sigma^{-}=\sum y_i$ are the discrete analogous of $V_{+}$ and
$V_{-}$. On the algebra spanned by these variables, we have an action
of $U_q(\bm)$, where $\bm$ is the negative sub-Borel algebra of
$\widehat{sl_2}$. Thus we get a structure of
$U_q(\bm)$-module-algebra. One of the problems we have, is to
understand this module-algebra in terms of $U_q(\bm)$, or, at least,
to caracterize by means of the quantum group $U_q(\bm)$ only,
a sub $U_q(\bm)$-module-algebra of 
$\CC[x_i^{\pm 1},y_i^{\pm 1},1\leq i\leq n]_q$ which possesses
the same fration field as $\CC[x_i^{\pm 1},y_i^{\pm 1},1\leq i\leq
n]_q$. We shall obtain this sub-module-algebra as a quantum homogeneous
space over $\BM$. The results proved below generalize [2].
\end{abstract}

\bigskip
\noindent
{\bf Classification AMS~:} 35Q53, 81R50.

\bigskip
\centerline{{\bf Mots cl{\'e}s}}

\smallskip
\noindent
Th{\'e}orie de Toda sur r{\'e}seau, morphismes de Poisson, feuilletage
symplectique, cellule de Schubert quantique, espace homog{\`e}ne
quantique, produit tensoriel tordu.

\bigskip
\centerline{{\bf Key words}}

\smallskip
\noindent
Toda theory on lattice, Poisson morphisms, symplectic foliation,
quantum Schubert cell, quantum homogeneous space, twisted tensorial
product.

\newpage
\selectlanguage{french}
\setcounter{section}{-1}
\section{Introduction}
\noindent
Une version quantique du syst{\`e}me sur r{\'e}seau de sinus-Gordon
a {\'e}t{\'e} introduite pour la premi{\`e}re fois
par Izergin et Korepin
([3] et [4]).
Ce syst{\`e}me a depuis {\'e}t{\'e} {\'e}tudi{\'e} et d{\'e}velopp{\'e} par de
nombreux auteurs.
Voir par exemple 
[5], [6], [7] et [8].
Au niveau classique, B. Enriquez et B. Feigin
ont montr{\'e} une bijection entre l'espace des phases du syst{\`e}me sur
r{\'e}seau et le groupe de Lie-Poisson $\widehat{\textnormal{SL}_2}$.
Nous nous proposons dans cet article de g{\'e}n{\'e}raliser 
ces r{\'e}sultats, et d'{\'e}tablir une bijection entre l'espace des
phases quantique du syst{\`e}me sur r{\'e}seau et un espace homog{\`e}ne
quantique
pour le groupe quantique $U_q(\widehat{\textnormal{sl}_2})$
que nous caract{\'e}riserons enti{\`e}rement dans le cas o{\`u} $q$ n'est
pas une racine de l'unit{\'e}.
\section{Le syst{\`e}me sur r{\'e}seau}
\label{modlib}

\subsection{D{\'e}finition de $A_q^{(n)}$.}
\label{modlib2}

\noindent
Soient $n\in{\NN}^*$ un entier fix{\'e}, et $q$ une ind{\'e}termin{\'e}e. 
On s'int{\'e}resse {\`a} l'alg{\`e}bre $A_q^{(n)}$ engendr{\'e}e sur  
$\CC[q,q^{-1}]$ 
par les $x_i^{\pm 1}$ et $y_i^{\pm 1}, 1\leq i\leq n$
soumis aux relations:
$x_ix_j=qx_j x_i,\, y_iy_j=qy_jy_i$,
$x_iy_j=q^{-1}y_jx_i,\, y_ix_j=q^{-1}x_jy_i$ pour $i<j$,
et $x_iy_i=q^{-1}y_ix_i$.

En utilisant par exemple des extensions de Ore ([9]),
on montre que $A_q^{(n)}$ est un anneau int{\`e}gre (non commutatif)
qui poss{\`e}de un corps de fractions,
et que c'est un $\CC[q,q^{-1}]$ module libre
dont une base est donn{\'e}e par la famille 
$\prod\limits_{1}^{n}x_i^{a_i}y_i^{b_i}$, avec 
$\forall i,\, (a_i,b_i)\in{\ZZ}^2$.

Pour $i\in\lbrace 1,\ldots,n\rbrace$, soit $\CC[x_i,x_{i}^{-1}]$
(resp. $\CC[y_i,y_{i}^{-1}]$) l'alg{\`e}bre engendr{\'e}e sur
$\CC[q,q^{-1}]$ par $x_i$ (resp. $y_i$) et
son inverse, que l'on rend gradu{\'e}e par $\deg{x_i}=1$ (resp. $\deg{y_i}=-1$).
Alors, en notant~: 
  $$
  \CC[x_1,x_{1}^{-1}]\ootimes
  \CC[y_1,y_{1}^{-1}]\ootimes\ldots\ootimes
  \CC[x_n,x_{n}^{-1}]\ootimes\CC[y_n,y_{n}^{-1}]
  $$ 
le produit tensoriel $q^{-1}$
tordu (ordonn{\'e}) de ces alg{\`e}bres gradu{\'e}es, on montre que
l'application~:
\begin{equation}\label{debiso}
  \begin{array}{rcl}
    f:\quad\CC[x_1,x_{1}^{-1}]\ootimes
           \ldots\ootimes\CC[y_n,y_{n}^{-1}]
           &\longrightarrow&A_q^{(n)}\\
    u_1\ootimes\ldots\ootimes v_n&
    \longmapsto&u_1\ldots v_n
  \end{array}
\end{equation}
est un isomorphisme d'alg{\`e}bres (voir appendice).

L'alg{\`e}bre $A_q^{(n)}$ est gradu{\'e}e en posant
\begin{equation}\label{gradu}
\deg{x_i}=-\deg{y_i}=1.
\end{equation}
L' isomorphisme $f$ ci-dessus est gradu{\'e}.

On posera pour $a$ et $b$ homog{\`e}nes,
\begin{equation}
[a,b]_q=ab-q^{(\deg{a})(\deg{b})}ba.
\end{equation}

On note {\'e}galement  $\Sigma^{+}=\sum x_i$ et $\Sigma^{-}=\sum y_i$.

\subsection{L'alg{\`e}bre $U_q(\bm)\subset
  U_q(\widehat{{\textnormal{sl}}_2})$.}
\label{gradu2}
\noindent
Soit $U_q(\bm)$ l'alg{\`e}bre de Hopf engendr{\'e}e sur
$\CC[q,q^{-1}]$
par $e_{\pm},k,k^{-1}$ et les relations~:
$kk^{-1}=k^{-1}k=1,\, ke_{\pm}k^{-1}=q^{\pm 1} e_{\pm},$
ainsi que les relations de Serre quantique entre $e_{\pm}$ et $e_{\mp}$~:
$$
e_{\pm}^3e_{\mp}-(q+1+q^{-1})(e_{\pm}^2e_{\mp}e_{\pm}
-e_{\pm}e_{\mp}e_{\pm}^2)-e_{\mp}e_{\pm}^3=0
$$
La comultiplication sur $U_q(\bm)$ est donn{\'e}e par~:
$$
\Delta e_{\pm}=k^{\pm 1}\otimes e_{\pm}+e_{\pm}\otimes 1,
\textnormal{ et }
\Delta k^{\pm 1}=k^{\pm 1}\otimes k^{\pm 1}.
$$
Nous n'utiliserons pas l'antipode dans cet article.

On dispose d'une graduation sur $U_q(\bm)$
d{\'e}finie par~:
$$
\deg{e_{\pm}}=\pm 1\quad\textnormal{ et }\quad\deg{k^{\pm 1}}=0.
$$
\subsection{ L'alg{\`e}bre $U_q(\nm)\subset U_q(\bm)$.}
\noindent
Par d{\'e}finition, $U_q(\nm)$ est la sous-alg{\`e}bre
de $U_q(\bm)$ engendr{\'e}e par $e_{\pm}$.

En notant par $U_q(\nm)\ootimes U_q(\nm)$
le produit tensoriel $q^{-1}$ tordu de $U_q(\nm)$
par lui m{\^e}me (voir appendice), on dispose d'une comultiplication tordue 
sur $U_q(\nm)$
qui est un morphisme d'alg{\`e}bres~:
\begin{equation}
\begin{array}{rcl}
\DDelta:\quad U_q(\nm)&\longrightarrow&
U_q(\nm)\ootimes U_q(\nm)\\
e_{\pm}&\longmapsto&1\ootimes e_{\pm}+e_{\pm}\ootimes 1
\end{array}
\end{equation}
On peut montrer que l'existence de $\DDelta$ {\'e}quivaut aux relations
de Serre quantique.
\subsection{ Action de $U_q(\bm)$ sur $A_q^{(n)}$.}\label{action}
\noindent
Pour $i\in\lbrace 1,\ldots,n\rbrace$, on dispose de deux morphismes
d'alg{\`e}bres gradu{\'e}es~:
$$
\begin{array}{rcl}
\varphi_i:\quad U_q(\nm)&\longrightarrow&\CC[x_i,x_{i}^{-1}]\\
e_{+}&\longmapsto&x_i\\
e_{-}&\longmapsto&0
\end{array}
$$
et
$$
\begin{array}{rcl}
\psi_i:\quad U_q(\nm)&\longrightarrow&\CC[y_i,y_{i}^{-1}]\\
e_{+}&\longmapsto&0\\
e_{-}&\longmapsto&y_i
\end{array}
$$
Donc, en composant avec $\DDelta^{(2n-1)}$ et en utilisant
l'isomorphisme $f$ de (\ref{debiso}), on obtient un morphisme
d'alg{\`e}bres gradu{\'e}es~:
\begin{equation}\label{ecran}
  \begin{array}{rcl}
    \varphi:\quad U_q(\nm)&\longrightarrow&A_q^{(n)}\\
    e_{\pm}&\longmapsto&\displaystyle\Sigma^{\pm}
  \end{array}
\end{equation}
Par suite, les op{\'e}rateurs $\Sigma^{\pm}$ satisfont les relations
de Serre quantique.
Donc, les relations de Serre quantique {\'e}tant homog{\^e}nes 
en $e_{+}$ et $e_{-}$, on peut d{\'e}finir pour
$(\lambda,\mu)\in{\CC}^{2}$, deux
actions $\sigma_{\lambda}$ et $\tau_{\mu}$ de $U_q(\nm)$ 
sur $A_q^{(n)}$ par~:
 $$
 \forall P\in A_q^{(n)}\,\textnormal{homog{\`e}ne},
 \quad \sigma_{\lambda}(e_{\pm}).P=
 \lambda\Sigma^{\pm}P,
 $$
et~:
 $$
 \forall P\in A_q^{(n)}\,\textnormal{homog{\`e}ne},
 \quad \tau_{\mu}(e_{\pm}).P=
 \mu q^{\pm \deg{P}}P\Sigma^{\pm}.
 $$
On remarque que~:
  $$
  \forall (x,y)\in{U_q(\nm)}^{2}\,\textnormal{homog{\`e}nes} ,\quad 
  \sigma_{\lambda}(x)\circ\tau_{\mu}(y)=q^{-\deg{x}\deg{y}}
  \tau_{\mu}(y)\circ\sigma_{\lambda}(x).
  $$
Par suite, on peut d{\'e}finir une action 
$\sigma_{\lambda}\ootimes\tau_{\mu}$ de $U_q(\nm)\ootimes
U_q(\nm)$ sur $A_q^{(n)}$, par~:
  $$
  \forall (x,y)\in{U_q(\nm)}^{2},\quad
  (\sigma_{\lambda}\ootimes\tau_{\mu})(x\ootimes y)=
  \tau_{\mu}(y)\circ\sigma_{\lambda}(x).
  $$
En composant avec $\DDelta$ qui est un morphisme d'alg{\`e}bres,
on obtient une action $f_{\lambda,\mu}$ de $U_q(\nm)$ sur $A_q^{(n)}$.
En prenant $\lambda=-\mu=1$, cette action est donn{\'e}e par~:
  $$
  \forall P\in A_q^{(n)}\,\textnormal{homog{\`e}ne},\quad 
  f_{1,-1}(e_{\pm}).P=[\Sigma^{\pm},P]_{q}.
  $$
On peut {\'e}tendre $f_{\lambda,\mu}$ en une action de 
$U_q(\bm)$ sur $A_q^{(n)}$ par~:
  $$
  \forall P\in A_q^{(n)}\,\textnormal{homog{\`e}ne},\quad 
  k.P=q^{\deg{P}}P\quad
  $$
Il est facile de voir que l'on obtient ainsi une structure de
$U_q(\bm)$-alg{\`e}bre module sur $A_q^{(n)}$.
Le calcul montre que pour $i\in\lbrace 1,\ldots,n\rbrace,
\, f_{1,-1}(e_{\pm}).\, x_{i}$ tout comme $f_{1,-1}(e_{\pm}).\, y_{i}$
est une somme de termes tous divisibles par $q-1$.
Par suite, en posant pour $P\in A_q^{(n)}$ homog{\`e}ne~:
  \begin{align}
  \label{ack}
  k.P&=q^{\deg{P}}P\\
  \label{acepm}
  e_{\pm}.P&=\displaystyle{1\over q-1}[\Sigma^{\pm},P]_{q},
  \end{align}
on obtient encore une action de $U_q(\bm)$ sur $A_q^{(n)}$.
Cette action est gradu{\'e}e pour $\deg{}$,
et munit $A_q^{(n)}$ d'une structure de $U_q(\bm)$-alg{\`e}bre
module. Dans toute la suite, nous allons nous int{\'e}resser {\`a}
d{\'e}crire cette structure.

\subsection{ D{\'e}finitions {\'e}l{\'e}mentaires des $q$ nombres et
factorielles quantiques.}\label{defdeb1}
\noindent

On d{\'e}finit successivement~:
\begin{itemize}
  \item pour $n\in\NN,\, [n]=\displaystyle{q^n -1\over q-1}
        =\sum_{k=0}^{n-1} q^k\in{\mathbb Z}[q,q^{-1}]$\,;
  \item pour $n\in\NN,\, [n]!=\displaystyle\prod\limits_{k=1}^{n}\, [k]$\,;
  \item pour $(n,p)\in{\ZZ}^2$,
        $$\C{n}{p}=\left\{
        \begin{array}{ll}
          0,&\text{si }p<0\text{ ou si }p>\text{Max}(0,n)\,;\\
          1,&\text{si }p=0\,;\\
          \displaystyle{[n]!\over [p]! [n-p]!},&\text{si }0<p\leq n\,;
        \end{array}
        \right.$$
  \item pour des entiers $N,a_1,\ldots,a_N$, 
        $$
        F_q(a_1,\ldots,a_N)=\displaystyle\prod\limits_{i=1}^{N-1}
        \C{a_i+a_{i+1}-1}{a_{i+1}}.
        $$
\end{itemize}

\subsection{ D{\'e}finition des $u_i$ et de $U_q^{(n)}$.}
\label{defuqn}

\noindent
Pour $i\in\NN$, soit $u_i$ l'{\'e}l{\'e}ment de 
$A_q^{(n)}$ d{\'e}fini par~:
\begin{equation}\label{defui}
  u_i= \displaystyle\sum\limits_{\alpha_1,\ldots,\,\alpha_{2n-1}
  \atop{\alpha_1+\ldots+\alpha_{2n-1}=i-1}}
  F_q(\alpha_{2n-1},\ldots,\alpha_1)
  (x_1y_1)^{-\alpha_1}
  \ldots (x_n y_n)^{-\alpha_{2n-1}}y_n^{-1}
\end{equation}
Cette d{\'e}finition a bien un sens car 
s'il existe $i\in\lbrace 1,\ldots,N\rbrace$
tel que $\a_i<0$, alors $F_q(\a_1,\ldots,\a_N)=0$.

On note $U_q^{(n)}$ la sous-alg{\`e}bre de $A_q^{(n)}$ engendr{\'e}e 
par les $u_i$ pour\break $i\in\lbrace 1,\ldots,2n\rbrace$.

Nous verrons que les {\'e}l{\'e}ments
$\displaystyle\prod\limits_{i=1}^{2n}u_i^{\a_i},\,
\a_i\in\NN$ forment une base de $U_q^{(n)}$.

\subsection{ D{\'e}finition de ${\cal A}_q\simeq \CC[\NP]_q$.}\label{defaq}

\noindent
Soient $k=\CC[q^{{1\over 4}},q^{-{1\over 4}}]$,
\begin{equation}\label{mat2}
  H = \begin{pmatrix}
        q^{-{1\over 4}}&0&0&0\\
        0&q^{{1\over 4}}&0&0\\
        0&0&q^{{1\over 4}}&0\\
        0&0&0&q^{-{1\over 4}}
      \end{pmatrix} \in M_4(k)\simeq M_2(k)^{{\otimes}^2},
\end{equation}
\begin{equation}\label{mat1}
  R({\lambda,\mu})=
  \begin{pmatrix}
    1&0&0&0\\
    0&{{\lambda -\mu}\over q^{-{1\over 2}}\lambda -q^{{1\over 2}}\mu}&
    {{(q^{-{1\over 2}}-q^{{1\over 2}})\mu}\over 
    q^{-{1\over 2}}\lambda -q^{{1\over 2}}\mu}&0\\
    0&{{(q^{-{1\over 2}}-q^{{1\over 2}})\lambda}\over 
    q^{-{1\over 2}}\lambda -q^{{1\over 2}}\mu}&
    {{\lambda -\mu}\over q^{-{1\over 2}}\lambda -q^{{1\over 2}}\mu}&0\\
    0&0&0&1
  \end{pmatrix} \in M_4(k)\simeq M_2(k)^{{\otimes}^2},
\end{equation}
et ${\cal A}$ l'alg{\`e}bre libre sur 
$\CC[q,q^{-1}]$
engendr{\'e}e par les $a_{i,j}^{(r)}$ pour $i,j\in\lbrace 1,2\rbrace$
et $r\in\NN$. 

On pose~:
\begin{align*}
a_{i,j}(\lambda)&=\sum\limits_{r=0}^{+\infty}a_{i,j}^{(r)}
\lambda^r,\\
{\cal L}(\lambda)&=[a_{i,j}(\lambda)],\\
{\cal L}^1(\lambda)&={\cal L}(\lambda)\otimes Id,\\
\textnormal{ et }\quad
{\cal L}^2(\mu)&=Id\otimes {\cal L}(\mu).
\end{align*}
Par d{\'e}finition, ${\cal A}_q$ est l'alg{\`e}bre quotiente de
${\cal A}$ par l'id{\'e}al engendr{\'e} par les relations~:
\begin{equation}\label{RR1}
a_{2,1}^{(0)}=0,
a_{1,1}^{(0)}=
a_{2,2}^{(0)}=1,
\end{equation}
ainsi que 
\begin{equation}\label{RR2}
R({\lambda,\mu}){\cal L}^1(\lambda)H{\cal L}^2(\mu)=
{\cal L}^2(\mu)H{\cal L}^1(\lambda)R({\lambda,\mu}),
\end{equation}
et 
\begin{equation}\label{RR3}
a_{1,1}(q\lambda)\bigl[
a_{2,2}(\lambda)-
a_{2,1}(\lambda)a_{1,1}(\lambda)^{-1}a_{1,2}(\lambda)\bigr]=1.
\end{equation}
Le calcul montre que ces relations sont {\`a} coefficients
dans $\CC[q,q^{-1}]$.
Donc, ceci a bien un sens. La derni{\`e}re relation est celle du
d{\'e}terminant quantique.

On montre que l'on peut d{\'e}finir une graduation sur ${\cal A}_q$ par 
$\deg{a_{i,j}^{(r)}}=i-j$. Il suffit pour cela de consid{\'e}rer
l'op{\'e}rateur $\lambda\displaystyle{\partial\over{\partial\lambda}}$.

On peut {\'e}galement d{\'e}finir une comultiplication tordue sur
${\cal A}_q$ par~:
$$
\begin{array}{rcl}
\DDelta:\quad {\cal A}_q&\longrightarrow&{\cal A}_q\ootimes{\cal A}_q\\
{\cal L}(\lambda)&\longmapsto&{\cal L}(\lambda)\ootimes{\cal L}(\lambda)
\end{array}
$$
Pour une d{\'e}monstration de ce fait, voir l'appendice.

L'alg{\`e}bre ${\cal A}_q$ est isomorphe {\`a}
$U_q(\nm)$ par l'isomorphisme gradu{\'e}~:
\begin{equation}\label{isonp}
\begin{array}{rcl}
U_q(\nm)&\longrightarrow&{\cal A}_q\\
e_{+}&\longmapsto&a_{2,1}^{(1)}\\
e_{-}&\longmapsto&a_{1,2}^{(0)}.
\end{array}
\end{equation}
Cet isomorphisme respecte les comultiplications tordues.
\subsection{ D{\'e}finition de ${\cal A}_q^{(n)}$.}\label{defaqn}
\noindent
On d{\'e}finit tout d'abord l'alg{\`e}bre ${\bar{{\cal A}}}_q^{(n)}$
comme {\'e}tant le quotient de ${\cal A}_q$ par l'id{\'e}al engendr{\'e} par
les {\'e}l{\'e}ments~: 
$a_{1,1}^{(r-1)},a_{1,2}^{(r-1)},a_{2,1}^{(r)},a_{2,2}^{(r)}$
pour $r> n$.

L'alg{\`e}bre ${\cal A}_q^{(n)}$ est l'alg{\`e}bre engendr{\'e}e par
g{\'e}n{\'e}rateurs~:
$$
a_{i,j}^{(k)},\, 0\leq k\leq n,\, {a_{1,1}^{(n-1)}}',\,
{a_{1,2}^{(n-1)}}',\, {a_{2,1}^{(n)}}',\,
{a_{2,2}^{(n)}}',
$$
et relations~:
\begin{align*}
a_{2,1}^{(0)}=a_{1,1}^{(n)}=a_{1,2}^{(n)}&=0,\\
a_{1,1}^{(0)}=a_{2,2}^{(0)}&=1,\\
\textnormal{(\ref{RR2}) et (\ref{RR3})}\,
\textnormal{ avec }{\cal L}(\lambda)=[a_{i,j}(\lambda)]\,&
\textnormal{ et }\,
a_{i,j}(\lambda)=\displaystyle\sum\limits_{k=0}^{n}
a_{i,j}^{(k)}\lambda^k,\\
a_{1,1}^{(n-1)}{a_{1,1}^{(n-1)}}'=
{a_{1,1}^{(n-1)}}' a_{1,1}^{(n-1)}&=1,\\
a_{1,2}^{(n-1)}{a_{1,2}^{(n-1)}}'=
{a_{1,2}^{(n-1)}}'a_{1,2}^{(n-1)}&=1,\\
a_{2,1}^{(n)}{a_{2,1}^{(n)}}'=
{a_{2,1}^{(n)}}'a_{2,1}^{(n)}&=1,\\
a_{2,2}^{(n)}{a_{2,2}^{(n)}}'=
{a_{2,2}^{(n)}}'a_{2,2}^{(n)}&=1.
\end{align*}
L'alg{\`e}bre ${\cal A}_q^{(n)}$ est trop grosse pour que l'on puisse
esp{\'e}rer obtenir une cellule de Schubert quantique.
Il faudrait imposer en plus les relations entre les $a_{i,j}^{(k)}$,
et ${a_{1,1}^{(n-1)}}',\, {a_{1,2}^{(n-1)}}',\,
{a_{2,1}^{(n)}}',\, {a_{2,2}^{(n)}}'$.
Bien que ce ne soit pas notre sujet d'{\'e}tude, nous reviendrons sur ce
point {\`a} la fin de notre article en {\'e}non\c cant une conjecture.

Quoiqu'il en soit, l'alg{\`e}bre ${\cal A}_q^{(n)}$ d{\'e}finie ci-dessus
contient la sous-alg{\`e}bre ${\cal U}_q^{(n)}$ d{\'e}finie ci-dessous, 
et c'est cette derni{\`e}re alg{\`e}bre qui va nous int{\'e}resser.

Notons que l'on a une fl{\`e}che naturelle~:
$$
{\bar{{\cal A}}}_q^{(n)}\longrightarrow {\cal A}_q^{(n)},
$$
et que si $C$ est une alg{\`e}bre et 
$$
f:\quad {\bar{{\cal A}}}_q^{(n)}\longrightarrow C
$$
est un morphisme tel que
$f\bigl(  a_{1,1}^{(n-1)}\bigr),\ldots,
f\bigl( a_{2,2}^{(n)}\bigr)$ soient inversibles,
alors, il existe~:
$$
g:\quad {\cal A}_q^{(n)}\longrightarrow C,
$$
un morphisme rendant commutatif le diagramme~:
$$
\begin{array}{rcl}
{\bar{{\cal A}}}_q^{(n)}&\longrightarrow&{\cal A}_q^{(n)}\\
f\searrow&&\swarrow g\\
&C&
\end{array}
$$
\subsection{ D{\'e}finition de ${\cal U}_q^{(n)}$.}\label{defuqng}
\noindent
Par d{\'e}finition, ${\cal U}_q^{(n)}$ est la sous-alg{\`e}bre de
${\cal A}_q^{(n)}$ engendr{\'e}e par les $2n$ premiers coefficients 
de $a_{2,2}(\lambda)^{-1}a_{2,1}(\lambda)$ d{\'e}velopp{\'e} en puissances
de $\lambda^{-1}$.
Ceci a bien un sens car 
$a_{2,2}(\lambda)$ est un polyn{\^o}me en $\lambda$
de degr{\'e} $n$ dont le coefficient dominant
$a_{2,2}^{(n)}$ est inversible.

Nous verrons que ${\cal U}_q^{(n)}$ est isomorphe {\`a} $U_q^{(n)}$
et s'identifie {\`a} un espace
homog{\`e}ne quantique dont on exprimera la limite classique.

\section{Le th{\'e}or{\`e}me principal}\label{letheo}

\noindent
Nous allons montrer le th{\'e}or{\`e}me suivant~:

\begin{thm}$\mbox{}$\vspace{-3mm}\\
\begin{itemize}
  \item[{\bf T1.}]
    Les $u_i,\, i\geq 1$ peuvent s'obtenir comme d{\'e}veloppement de la fraction
    continue quantique suivante~:
    \begin{multline}\label{fracq1}
      \Biggl ( 1+\biggl ( 1+\Bigl ( 1+\ldots\bigl (
      1+(\lambda x_1 y_1)^{-1}\bigr )^{-1}(\lambda y_1 x_2)^{-1}
      \Bigr )^{-1}\ldots\\
      \ldots (\lambda y_{n-1} x_n)^{-1}\biggr )^{-1}
      (\lambda x_n y_n)^{-1}\Biggr )^{-1}y_n^{-1}\,
    \end{multline}
    \centerline{$=\displaystyle\sum\limits_{k=0}^{\infty}
    (-1)^k u_{k+1}\lambda^{-k}$}
  \item[{\bf T2.}]
    La sous-alg{\`e}bre $U_q^{(n)}$ engendr{\'e} par les $u_i$ pour
    $i\in\lbrace 1,\ldots,2n\rbrace$ est une sous
    $U_q(\bm)$-alg{\`e}bre-module de $A_q^{(n)}$ qui poss{\`e}de
    le m{\^e}me corps de fractions que $A_q^{(n)}$.
    Elle est isomorphe {\`a} l'alg{\`e}bre d{\'e}finie par g{\'e}n{\'e}rateurs~:
    $t_i,\, 1\leq i\leq 2n$, et relations~:
    $$
      \forall i<j,\quad 
        q [t_i,t_j]=(q-1)\sum_{k=i}^{j-1}t_kt_{i+j-k}
    $$
  \item[{\bf T3.}]
    La sous-alg{\`e}bre ${\cal U}_q^{(n)}$
    de ${\cal A}_q^{(n)}$
    s'identifie {\`a} un espace homog{\`e}ne quantique et est isomorphe
    {\`a} $U_q^{(n)}$.
  \item[{\bf T4.}]
    Sp{\'e}cialisons $q$ en un nombre complexe non nul, non racine de
    l'unit{\'e}, et adjoignons {\`a} $U_q(\bm)$ une racine carr{\'e}e de
    $k$, de sorte que $A_q^{(n)}$ reste une
    $U_q(\bm)$-alg{\`e}bre-module, en posant, en plus de (\ref{ack})
    et (\ref{acepm})~:
      $$
       k^{{1\over 2}}.P=q^{{1\over 2}\deg{P}}P,
      $$
    pour $P$ homog{\`e}ne. Alors, 
    il existe $I_q^{(n)}$ un id{\'e}al de Hopf {\`a} droite gradu{\'e} 
    de $U_q(\bm)$,
    et un couplage de Hopf gradu{\'e} non d{\'e}g{\'e}n{\'e}r{\'e}:
     $$
      {U_q(\bm)}/{I_q^{(n)}}\times U_q^{(n)}\longrightarrow\CC
     $$
    qui induisent un monomorphisme d'alg{\`e}bres~:
     $$
      U_q^{(n)}\hookrightarrow {\Bigl( {{U_q(\bm)}/I_q^{(n)}}\Bigr)}^{*}
     $$
    tel que ${\Bigl( {{U_q(\bm)}/I_q^{(n)}}\Bigr)}^{*}$ soit la compl{\'e}tion
    formelle de $U_q^{(n)}$.
    Autrement dit, $U_q^{(n)}$ est un module coinduit~:
     $$
      U_q^{(n)}\simeq{\bigl( \CC\otimes_{I_q^{(n)}} U_q(\bm)\bigr)}^*
     $$
    o{\`u} $\CC$ est muni d'une structure de $I_q^{(n)}$-module
    {\`a} droite triviale.
  \item[{\bf T5.}]
    En termes de nouvelles r{\'e}alisations de Drinfeld ([10]), $I_q^{(n)}$ est
    l'id{\'e}al {\`a} droite engendr{\'e} par~:
    $h[0]-1,h[-k],k>0;e[-l],l\geq 1;f[-m],m\geq 2n$.
\end{itemize}
\end{thm}

\section{Plan de la d{\'e}monstration}

\noindent
La d{\'e}monstration se fera en sept {\'e}tapes, {\'E}1,...,{\'E}7.

\begin{itemize}
  \item[{\bf {\'E}1.}]
    Etude du cas classique.

    Celle-ci nous conduira naturellement {\`a} consid{\'e}rer les limites
    classiques des $u_i$ d{\'e}finis ci-dessus.
  \item[{\bf {\'E}2.}]
    D{\'e}monstration de T1.
  \item[{\bf {\'E}3.}]
    Calcul de l'action de $U_q(\bm)$ sur les $u_i,\, 1\leq i\leq 2n$.

    En vu de montrer que $U_q^{(n)}$ est un $U_q(\bm)$-module,
    nous calculons directement l'action de $e_{\pm}$ sur les 
    $u_i,\, 1\leq i\leq 2n$.
  \item[{\bf {\'E}4.}]
    Construction d'un morphisme d'alg{\`e}bres de ${\cal A}_q^{(n)}$
    dans $A_q^{(n)}$.

    Nous construisons d'abord un morphisme de  ${\cal A}_q$
    dans $A_q^{(n)}$.
    Ce morphisme n'est rien d'autre que l'application $\varphi$
    de (\ref{ecran}) lue de mani{\`e}re groupique grace
    {\`a} l'isomorphisme (\ref{isonp}).
  \item[{\bf {\'E}5.}]
    Relations de commutations entre les $u_i,\, 1\leq i\leq 2n$.

    Celles-ci sont obtenues grace aux relations de commutation
    dans ${\cal U}_q^{(n)}$ via le morphisme pr{\'e}c{\'e}dent et le 
    plongement de ${\cal U}_q^{(n)}$ dans ${\cal A}_q^{(n)}$.
  \item[{\bf {\'E}6.}]
    D{\'e}monstration de T2 et de T3.

    On montre qu'une base de $U_q^{(n)}$ est donn{\'e}e par la famille
    $\displaystyle\prod\limits_{i=1}^{2n}u_i^{\a_i}$ pour $\a_i\in\NN$.
  \item[{\bf {\'E}7.}]
    Calcul de $I_q^{(n)}$ et fin de la d{\'e}monstration.

    L'alg{\`e}bre $U_q(\bm)$ {\'e}tant isomorphe {\`a}  $\CC[\BP]_q$,
    on cherche $I_q^{(n)}$ dans $\CC[\BP]_q$.
    Nous y parvenons apr{\`e}s avoir mis en {\'e}vidence le fait que
    l'action de $U_q(\bm)$ dans $A_q^{(n)}$ est bi-gradu{\'e}e,
    et examin{\'e} l'action des $a_{i,j}^{(k)}$ dans $U_q^{(n)}$.
\end{itemize}

\setcounter{section}{0}
\renewcommand{\thesection}{{\'E}{\arabic{section}}}

\section{Le cas classique (voir {\'e}galement [2])}

\noindent
Nous rappelons le proc{\'e}d{\'e} utilis{\'e} dans [2]
pour obtenir une sous-alg{\`e}bre de 
$\CC[x_i^{\pm 1},y_i^{\pm 1}]$ ayant
le m{\^e}me corps des fractions que $\CC[x_i^{\pm 1},y_i^{\pm 1}]$,
stable par l'action de $\bm$ et facilement exprimable 
en termes de $U\bm$.

\subsection{}\label{schubclass}

\noindent
Soit $\lambda$ une ind{\'e}termin{\'e}e, et 
$\widehat{G}=SL_2\Bigl( \CC\bigl((\lambda^{-1})\bigr)\Bigr)$.
On pose~:
$$
e=\begin{pmatrix}
0&1\\0&0\end{pmatrix},\quad f=\begin{pmatrix}
0&0\\1&0\end{pmatrix},\quad h=\begin{pmatrix}
1&0\\0&-1\end{pmatrix}
$$
et
$$\begin{array}{lcl}
  \ghg &=& {\textnormal sl}_2\otimes_{\CC}\CC\bigl((\lambda^{-1})\bigr)\,;\\
  \bm  &=& \bigl( e\otimes\lambda^{-1}\CC[[\lambda^{-1}]]\bigr)
           \oplus\bigl( f\otimes\CC[[\lambda^{-1}]]\bigr)
           \oplus\bigl( h\otimes\CC[[\lambda^{-1}]]\bigr)\,;\\
  \np  &=& \bigl( e\otimes\CC[\lambda]\bigr)
           \oplus\bigl( f\otimes\lambda\CC[\lambda]\bigr)
           \oplus\bigl( h\otimes\lambda\CC[\lambda]\bigr)\,;\\
  \BM  &=& \exp(\bm)\textnormal{ et }\NP=\exp(\np).
\end{array}$$
On a~: 
\begin{equation}\label{decompog}
\ghg=\bm\oplus\np.
\end{equation}
Soit $P_0$ le bichamp de Poisson d{\'e}fini sur
$\GHG$ par $P_0=r^{g}-r^{d}$ o{\`u} $r$ est la r-matrice
trigonom{\'e}trique standard~:
$$
r=\displaystyle{1\over 2}
\displaystyle{\lambda+\mu\over\lambda-\mu}h\otimes h
+\displaystyle{2\over\lambda-\mu}
(\lambda e\otimes f+\mu f\otimes e).
$$
Soit $P_1$ le bichamp d{\'e}fini par $P_1=h^g\otimes h^d-h^d\otimes
h^g$.
Le crochet de Lie $[h^g,h^d]$ est nul. Par suite, le crochet de Schouten 
$[P_1,P_1]=0$. De m{\^e}me, le fait que $r$ soit 
${\textnormal ad}_h$-invariant entra{\^\i}ne que le crochet de Schouten
$[P_0,P_1]=0$. D'o{\`u} une structure de Poisson sur $\GHG$
donn{\'e}e par $P=-\displaystyle{1\over 2}P_0
+\displaystyle{1\over 4}P_1$.

Le sous-groupe $\BM$ est une sous-vari{\'e}t{\'e} de Poisson pour la
structure d{\'e}finie par $P_0$ comme pour celle d{\'e}finie par $P$.
Le sous-groupe $\NP$ n'est pas une sous-vari{\'e}t{\'e} de Poisson
pour la structure d{\'e}finie par $P_0$. Par contre, il l'est pour celle
d{\'e}finie par $P$.
Dans la suite, lorsque nous parlerons de
structure de Poisson sur $\GHG$ ( ou sur l'une de ses sous-vari{\'e}t{\'e}s),
il s'agira toujours de la structure d{\'e}finie par $P$.

En g{\'e}n{\'e}ralisant un r{\'e}sultat de DeConcini, Kac 
et Procesi [11] qui d{\'e}coule lui m{\^e}me d'un th{\'e}or{\`e}me de
Semenov-Tian-Shansky [12], 
on montre que les feuilles symplectiques de $\NP$ sont les 
$\BM w\BM\cap\NP$ o{\`u} $w\in W, W$ d{\'e}signant le groupe de Weyl 
de $\GHG$.

\smallskip

Notons $w_n=\begin{pmatrix}
\lambda^{-n}&0\\
0&\lambda^{n}
\end{pmatrix},
S=\BM\cap w_n^{-1}\BM w_n,\, {\cal U}^{(n)}=
{}_{{}_{\displaystyle S}\!\diagdown}\!\! \BM$
l'ensemble de ses classes {\`a} gauche dans $\BM$, et
$$
\begin{array}{rcccl}
\pi:\quad \BM w_n\BM&\longrightarrow&
{}_{{}_{\displaystyle \BM}\!\diagdown}\!\! \BM w_n\BM&
\overset{\sim}{\longrightarrow}&
{}_{{}_{\displaystyle S}\!\diagdown}\!\! \BM=\, {\cal U}^{(n)}\\
aw_nb&\longmapsto&\textnormal{cl}(a w_n b)&\longmapsto&
\textnormal{cl}(b)
\end{array}
$$
Posons {\'e}galement ${\cal A}^{(n)}=\BM w_n\BM\cap\NP$.
On v{\'e}rifie que  ${\cal A}^{(n)}$ est l'ensemble des matrices
$$
\begin{pmatrix}
a&b\\
c&d
\end{pmatrix}
\in SL_2\bigl( \CC[\lambda]\bigr)
$$
telles que
\begin{equation}\label{jor1}
  a(0)= d(0) = 1,\, c(0) = 0,
\end{equation}
et
\begin{equation}\label{jor2}
  \begin{array}{l@{\;}c@{\;}l@{\;}c@{\;}l}
    \deg{a} &=& \deg{b} &=& n-1\,;\\
    \deg{c} &=& \deg{d} &=& n.
  \end{array}
\end{equation}
On montre facilement que la restriction $\pi'$ de $\pi$ {\`a}
${\cal A}^{(n)}$ est injective, et que 
$$
S=\Bigl\lbrace\, \begin{pmatrix}
a&,\lambda^{-1}b\\
\lambda^{-2n}c&d
\end{pmatrix}
,\, a,b,c,d\in\CC[[\lambda^{-1}]]\Bigr\rbrace .
$$
D'o{\`u} l'on voit que ${\cal U}^{(n)}$
s'identifie {\`a}
${\CC}_{2n-1}[\lambda^{-1}]={\CC[[\lambda^{-1}]]}/{\lambda^{-2n}
\CC[[\lambda^{-1}]]}$ en associant {\`a}
$\rho\in{\CC}_{2n-1}[\lambda^{-1}]$
la classe de la matrice
$\begin{pmatrix}
1&0\\
\rho&1
\end{pmatrix}$ dans ${\cal U}^{(n)}$.

Dans cette identification, l'application $\pi'$ n'est autre
que~: 
\begin{equation}\label{ident1}
M=\begin{pmatrix}
a&b\\
c&d
\end{pmatrix}
\longrightarrow\textnormal{cl}(\displaystyle{c\over d})\quad
\pmod{\lambda^{-2n}}
\end{equation}
Le sous-groupe $S$ {\'e}tant une sous-vari{\'e}t{\'e} de Poisson de
$\BM$, on peut munir l'espace homog{\`e}ne
${\cal U}^{(n)}={}_{{}_{\displaystyle S}\!\diagdown}\!\! \BM$
d'une structure de Poisson telle que $\pi$ et $\pi'$
soient de Poisson.

L'action d'habillage de $\bm$ sur $\NP$ est
d{\'e}finie par~:
$$
\forall x\in\bm\,\forall n\in\NP,\,
X_x(n)=(nxn^{-1})_+n,
$$
en notant 
par $y_+$ la projection de $y\in\ghg$ sur $\np$
parall{\`e}lement {\`a} $\bm$, suivant la 
d{\'e}composition (\ref{decompog}).

Par $\pi',$ la restriction de cette action d'habillage de $\bm$ sur
la feuille symplectique ${\cal A}^{(n)}=\BM w_n\BM\cap\NP$ de $\NP$
se retrouve en l'action
par translation {\`a} droite de
$\bm$ sur 
${}_{{}_{\displaystyle \BM}\!\diagdown}\!\! \BM w_n\BM
\simeq {}_{{}_{\displaystyle S}\!\diagdown}\!\! \BM={\cal U}^{(n)}$.

Soit ${\mathfrak s}$ l'alg{\`e}bre de Lie de $S$, et 
$I_{\textnormal{cl}}^{(n)}$ l'id{\'e}al
{\`a} droite de $U\bm$
engendr{\'e} par ${\mathfrak s}$.
Tous les {\'e}l{\'e}ments de $\bm$ {\'e}tant primitifs dans 
$U\bm$ muni de sa structure d'alg{\`e}bre de Hopf naturelle,
$I_{\textnormal{cl}}^{(n)}$ est un id{\'e}al de Hopf. 
De plus, d'apr{\`e}s la forme de $S$,
on a 
\begin{equation}\label{iclass}
I_{\textnormal{cl}}^{(n)}=
<h-1,h\otimes\lambda^{-k},k> 0;e\otimes\lambda^{-l},l\geq 1;
f\otimes\lambda^{-m},m\geq 2n>.
\end{equation}
Si $X$ est un champ de vecteur engendr{\'e} par translation {\`a} droite 
d'un {\'e}l{\'e}ment de ${\mathfrak s}$,
alors clairement~:
$$
\forall f\in\CC[\, {\cal U}^{(n)}],\quad
\textnormal{ev}_{S}(X.f)=0,
$$
en posant $\textnormal{ev}_{S}(g)=g(S)$ pour 
$g\in\CC[{}_{{}_{\displaystyle S}\!\diagdown}\!\! \BM]$.
De plus, on a un couplage non-d{\'e}g{\'e}n{\'e}r{\'e}~:
\begin{equation}\label{coupla1}
\begin{array}{rcl}
U\bm {}_{/\!{}_{\displaystyle I_{\textnormal{cl}}^{(n)}}}\times
\CC[\, {\cal U}^{(n)}]&\longrightarrow&\CC\\
({\bar x},P)&\longmapsto&\textnormal{ev}_{S}(x.P)
\end{array}
\end{equation}
avec un plongement d'alg{\`e}bres~:
$$
\CC[\, {\cal U}^{(n)}]
\hookrightarrow {\Bigl( U\bm {}_{/\!{}_{\displaystyle 
I_{\textnormal{cl}}^{(n)}}}\Bigr)}^{*}
$$
tel que 
${\Bigl( U\bm {}_{/\!{}_{\displaystyle 
I_{\textnormal{cl}}^{(n)}}}\Bigr)}^{*}$
soit la compl{\'e}tion formelle
de $\CC[\, {\cal U}^{(n)}]$.
\subsection{}
\noindent
On munit ${\CC}^{2n}$ d'une structure de Poisson 
par la limite pour $q\rightarrow 1$ du $q$-crochet $[\, ,\, ]_q$ 
et nous ne nous
int{\'e}resserons qu'{\`a} cette structure.

Explicitement, si $x_i$ et $y_i$ d{\'e}signent les fonctions
coordonn{\'e}es ($1\leq i\leq n$), on a~: 
$
\{ x_i,x_j\} =x_ix_j,\,
\{ y_i,y_j\} =y_iy_j,\,
\{ x_i,y_j\} =-x_iy_j,\,
\{ y_i,x_j\}=-y_ix_j \textnormal{ pour }i<j,
$
et $\{ x_i,y_i\} =-x_iy_i.$

La structure de $U_q(\bm)$-alg{\`e}bre-module sur $A_n$
se traduit
en une action par champ de vecteurs de $\bm$
sur $\CC^{2n}$.
Pour tout $P\in\CC[x_i,y_i,1\le i\le n]$ homog{\`e}ne 
($\deg{x_i}=-\deg{y_i}=1$),
on a~:
$$
k.P=\displaystyle{1\over 2}(\deg{P})P\textnormal{ et }e_{\pm}.P=
\{ \displaystyle\Sigma^{\pm},P\}
$$
\subsection{}\label{phicl}
\noindent
Soit maintenant $\varphi_{\textnormal{cl}}$ l'application suivante~:
$$
\begin{array}{rcl}
\varphi_{\textnormal{cl}}:\quad \CC^{2n}&\longrightarrow&\NP\subset\GHG\\
(x_1,\ldots,y_n)&\longmapsto&\displaystyle\prod_{i=1}^{n}
\begin{pmatrix}
1&0\\
\lambda x_i&1
\end{pmatrix}
\begin{pmatrix}
1&y_i\\
0&1
\end{pmatrix}
\end{array}
$$
On montre que $\varphi_{\textnormal{cl}}$ est un morphisme de Poisson.
Or, ${{\CC}^{*}}^{2n}$ est une feuille symplectique
de $\CC^{2n}$. Donc, d'apr{\`e}s ce qui pr{\'e}c{\`e}de, 
$$
\exists w\in W,\, \phi\bigl({{\CC}^{*}}^{2n}\bigr)\subset
\BM w\BM\cap\NP.
$$
En regardant pr{\'e}cis{\'e}ment l'image de ${{\CC}^{*}}^{2n}$
par l'application $\varphi_{\textnormal{cl}}$, on voit que
n{\'e}cessairement 
$w=w_n$.

Par ailleurs, en notant
$\varphi_{\textnormal{cl}}^*$ le comorphisme associ{\'e} {\`a} 
$\varphi_{\textnormal{cl}}$,
on a avec des notations {\'e}videntes, 
$\varphi_{\textnormal{cl}}^*(a_{2,1}^{(1)})=\Sigma^{+}$ et 
$\varphi_{\textnormal{cl}}^*(a_{1,2}^{(0)})=\Sigma^{-}$.
Apr{\`e}s calculs, il en r{\'e}sulte que l'action de $\bm$
sur ${{\CC}^{*}}^{2n}$ par champs de vecteurs
se retrouve
en l'action d'habillage de $\bm$ sur $\NP$.

On est dans la situation suivante o{\`u} toutes les fl{\`e}ches sont des
morphismes de Poisson, et leurs comorphismes des morphismes de $U\bm$-
modules~:
\begin{equation}\label{diagcl}
\begin{array}{rcl}
{{\CC}^{*}}^{2n}&\overset{\varphi_{\textnormal{cl}}}{\longrightarrow}&
\quad {\cal A}^{(n)}\\
\pi'\circ\varphi_{\textnormal{cl}}\searrow&&\swarrow\pi'\\
&{\cal U}^{(n)}&
\end{array}
\end{equation}
Grace {\`a} l'identification entre ${\cal U}^{(n)}$
et $\CC_{2n-1}[\lambda^{-1}]$, notons $\a_i\in\CC[\, {\cal U}^{(n)}]$
la fonction coordonn{\'e}e qui donne le coefficient devant
$\lambda^{-i+1}$. 
Alors, 
$$
\CC[\a_1,\ldots,\a_{2n}]=
\CC[\, {\cal U}^{(n)}],
$$
et, d'apr{\`e}s (\ref{ident1}),
l'application
$(\pi'\circ\varphi_{\textnormal{cl}})^*$ n'est autre que~:
\begin{equation}\label{alphaiclass}
\begin{array}{rcl}
\CC[\, {\cal U}^{(n)}]&\longrightarrow&\CC[x_i^{\pm 1},y_i^{\pm 1}]\\
\a_i&\longmapsto&u_{\textnormal{i,cl}}=
\textnormal{coeff devant }\lambda^{-i+1}
\textnormal{ de }\varphi_{\textnormal{cl}}^*\Biggl(
\displaystyle{a_{2,1}(\lambda)\over a_{2,2}(\lambda)}
\Biggr)
\end{array}
\end{equation}
Notons $U^{(n)}$ la sous-alg{\`e}bre de 
$\CC[x_i^{\pm 1},y_i^{\pm 1}]$ engendr{\'e}e par les 
$u_{\textnormal{i,cl}}$. C'est un $U\bm$-sous-module de
$\CC[x_i^{\pm 1},y_i^{\pm 1}]$.

Le calcul montre que
\begin{equation}
\displaystyle{a_{2,1}\over a_{2,2}}
\bigl( \varphi_{\textnormal{cl}}(x_1,\ldots,y_n)\bigr)=
\displaystyle{1\over1+\displaystyle{(\lambda x_n y_n)^{-1}\over
1+\displaystyle{(\lambda y_{n-1}x_n)^{-1}\over{\displaystyle{
\ddots
(\lambda y_1x_2)^{-1}\over
1+(\lambda x_1y_1)^{-1}}}}}}y_n^{-1}.
\end{equation}
\begin{prop}\label{cl2}
Soient $t_i,1\leq i\leq n$ des variables commutatives.
Alors, dans $\CC[[t_i,\, 1\leq i\leq n]]$, on a~:
$$
\cfrac{1}{1-\cfrac{t_1}{1-\cfrac{t_2}{\ddots\cfrac{t_{n-1}}{1-t_n}}}}=
\displaystyle\sum_{\alpha_1,\ldots,\alpha_n}
F(\alpha_1,\ldots,\alpha_n) t_1^{\alpha_1}\ldots t_n^{\alpha_n}
$$
avec par d{\'e}finition,
$F(\a_1,\ldots,\a_n)=\displaystyle\prod_{i=1}^{n-1}
\cnp{\alpha_i+\alpha_{i+1}-1}{\alpha_{i+1}}$.
\end{prop}
\begin{dem}
Par r{\'e}currence sur $n$ {\`a} partir de la formule~:
$$
{\Bigl(\displaystyle{1\over 1+X}\Bigr)}^n=
\displaystyle\sum\limits_{p=0}^{\infty}(-1)^p
\cnp{n+p-1}{p}X^p
$$
vraie dans $\CC[[X]]$ pour tout $n\in\NN$.
\end{dem}

Il en r{\'e}sulte que pour $i\in\lbrace 1,\ldots,2n\rbrace,$
\begin{equation}\label{expruicl}
u_{i,\textnormal{cl}}=
\displaystyle\sum\limits_{\alpha_1,\ldots,\,\alpha_{2n-1}
\atop{\alpha_1+\ldots+\alpha_{2n-1}=i-1}}
F(\alpha_{2n-1},\ldots,\alpha_1)
(x_1y_1)^{-\alpha_1}
\ldots (x_n y_n)^{-\alpha_{2n-1}}y_n^{-1}.
\end{equation}
De plus, en d{\'e}finissant par la m{\^e}me formule 
$u_{i,\textnormal{cl}}$ pour $i>2n$, on a~:
\begin{equation}\label{somuicl}
\displaystyle{a_{2,1}\over a_{2,2}}
\bigl(
\varphi_{\textnormal{cl}}(x_1,\ldots,y_n)\bigr)=\displaystyle\sum\limits_{k=0}^{\infty}
(-1)^k u_{k+1,\textnormal{cl}}\lambda^{-k}.
\end{equation}
On constate que pour $k\in\lbrace 1,\ldots,
\EE (\displaystyle{n-1\over 2})\rbrace$,
$$
u_{2k}\in\CC[x_{n-k+1}^{-1},y_{n-k+1}^{-1},\ldots,y_n^{-1}]
\textnormal{ et }
u_{2k+1}\in\CC[y_{n-k}^{-1},x_{n-k+1}^{-1},\ldots,y_n^{-1}].
$$
En particulier, on voit que $(\pi'\circ\varphi_{\textnormal{cl}})^*$
est injective et birationnelle.
Donc $U^{(n)}$ poss{\`e}de le m{\^e}me corps des fractions que 
$\CC[x_i^{\pm 1},y_i^{\pm 1}]$ et est un $U\bm$-module
isomorphe {\`a} $\CC[\, {\cal U}^{(n)}]$.

De plus, $(\pi'\circ\varphi_{\textnormal{cl}})^*$ injective entra{\^\i}ne que
$\pi'\circ\varphi_{\textnormal{cl}}$ est dense. Or $\pi'$ est injective.
Donc, $\varphi_{\textnormal{cl}}$ et $\pi'$ sont denses.
Par suite, en termes d'alg{\`e}bres de fonctions, le diagramme (\ref{diagcl})
peut se r{\'e}ecrire en le diagramme commutatif suivant
o{\`u} toutes les fl{\`e}ches sont injectives~:
\begin{equation}\label{diagclfonc}
\begin{array}{rcl}
\CC[\, {\cal U}^{(n)}]&\overset{{\pi'}^{*}}{\longrightarrow}&
\, \CC[\, {\cal A}^{(n)}]\\
({\pi'}\circ\varphi_{\textnormal{cl}})^{*}\searrow&&\swarrow
\varphi_{\textnormal{cl}}^{*}\\
&\CC[x_i^{\pm 1},y_i^{\pm 1}]&
\end{array}
\end{equation}
L'image de $\CC[\, {\cal U}^{(n)}]$ est $U^{(n)}$.
On peut caract{\'e}riser explicitement les alg{\`e}bres de Poisson
$\CC[\, {\cal A}^{(n)}]$
et $\CC[\, {\cal U}^{(n)}]$~:

\noindent
- Pour $\CC[\, {\cal A}^{(n)}]$, on d{\'e}finit tout d'abord
l'alg{\`e}bre (commutative) sur $\CC$ donn{\'e}e par g{\'e}n{\'e}rateurs~:
$$
a_{1,1}^{(0)},\ldots,a_{1,1}^{(n-1)},
a_{1,2}^{(0)},\ldots,a_{1,2}^{(n-1)},
a_{2,1}^{(0)},\ldots,a_{2,1}^{(n)},
a_{2,2}^{(0)},\ldots,a_{2,2}^{(n)},
$$
et relations~:
\begin{align*}
\forall i,i',j,j',k,k',\, a_{i,j}^{(k)}
a_{i',j'}^{(k')}&=a_{i',j'}^{(k')}a_{i,j}^{(k)}\\
a_{1,1}^{(0)}-1=a_{2,2}^{(0)}-1=a_{2,1}^{(0)}&=0\\
a_{1,1}(\lambda)a_{2,2}(\lambda)-
a_{2,1}(\lambda)a_{1,2}(\lambda)&=1,
\end{align*}
avec
$$
\forall i,j,\, a_{i,j}(\lambda)=\sum_k a_{i,j}^{(k)}\lambda^{k}.
$$
L'alg{\`e}bre $\CC[\, {\cal A}^{(n)}]$ s'obtient alors
en localisant suivant la partie multiplicative
engendr{\'e}e par 
$a_{1,1}^{(n-1)},a_{1,2}^{(n-1)},
a_{2,1}^{(n)},a_{2,2}^{(n)}$.
Ceci est du au fait que dans ${\cal A}^{(n)}$, les degr{\'e}s
sont fix{\'e}s \'gaux {\`a} $n-1$ ou $n$.

Par ailleurs, apr{\`e}s calculs, on v{\'e}rifie que 
les crochets de Poisson entre les $a_{i,j}^{(k)}$
sont donn{\'e}s par les formules suivantes~:
\begin{align}
\notag
\forall i,j,\quad \{ a_{i,j}(\lambda),a_{i,j}(\mu)\} &=0\\
\notag
(\lambda -\mu)
\{ a_{1,1}(\lambda),a_{1,2}(\mu) \} &=
\lambda a_{1,1}(\lambda)a_{1,2}(\mu)
-\lambda a_{1,1}(\mu)a_{1,2}(\lambda)\\
\notag
(\lambda -\mu)
\{ a_{1,1}(\lambda),a_{2,1}(\mu) \} &=
-\lambda a_{1,1}(\lambda)a_{2,1}(\mu)+\mu
a_{1,1}(\mu)a_{2,1}(\lambda)\\
\notag
(\lambda -\mu)
\{ a_{1,1}(\lambda),a_{2,2}(\mu) \} &=
-\lambda a_{1,2}(\lambda)a_{2,1}(\mu)+\mu 
a_{1,2}(\mu)a_{2,1}(\lambda)\\
\notag
(\lambda -\mu)
\{ a_{1,2}(\lambda),a_{2,1}(\mu) \} &=
-\mu a_{1,1}(\lambda)a_{2,2}(\mu)+\mu 
a_{1,1}(\mu)a_{2,2}(\lambda)\\
\notag
(\lambda -\mu)
\{ a_{1,2}(\lambda),a_{2,2}(\mu) \} &=
-\mu a_{1,2}(\lambda)a_{2,2}(\mu)+\mu 
a_{1,2}(\mu)a_{2,2}(\lambda)\\
\label{vadon}
(\lambda -\mu)
\{ a_{2,1}(\lambda),a_{2,2}(\mu) \} &=
\mu a_{2,1}(\lambda)a_{2,2}(\mu)-\lambda
a_{2,1}(\mu)a_{2,2}(\lambda).
\end{align}

\noindent
- L'alg{\`e}bre $\CC[\, {\cal U}^{(n)}]$ des fonctions sur l'espace
homog{\`e}ne ${\cal U}^{(n)}$ s'identifie par l'application injective
$(\pi')^{*}$ {\`a} la sous-alg{\`e}bre de $\CC[\, {\cal A}^{(n)}]$
engendr{\'e}e par les $2n$ premiers coefficients du d{\'e}veloppement
en $\lambda^{-1}$ de 
$v(\lambda^{-1}):=\displaystyle{a_{2,1}(\lambda)\over
a_{2,2}(\lambda)}$.
Ces coefficients sont les images des $\a_i$ d{\'e}finis
plus haut. Ils sont alg{\'e}briquement ind{\'e}pendants,
et on a $\CC[\, {\cal U}^{(n)}]=\CC[\a_{1},\ldots,\a_{n}]$
(en identifiant les $\a_i$ avec leurs images).

En exploitant l'{\'e}galit{\'e} (\ref{vadon}) et la propri{\'e}t{\'e} de Leibnitz
du crochet de Poisson, on obtient la relation~:
\begin{equation}\label{avoir+}
(\lambda -\mu)
\bigl\lbrace
v(\lambda^{-1}),v(\mu^{-1})
\bigr\rbrace
=\bigl(
\mu v(\lambda^{-1})-\lambda v(\mu^{-1})
\bigr)
\bigl(
v(\lambda^{-1})-v(\mu^{-1})
\bigr),
\end{equation}
ce qui donne apr{\`e}s d{\'e}veloppement~:
\begin{equation}\label{ahoha}
\forall i<j,\, \{ \a_{i},\a_{j}\}
=\displaystyle\sum\limits_{k=i}^{j-1}
\a_{k}\a_{i+j-k}.
\end{equation}
Ces formules d{\'e}crivent l'alg{\`e}bre de Poisson 
$\CC[{\cal U}^{(n)}]$.
\bigskip

Revenons au diagramme (\ref{diagclfonc}).
La fonction $(\pi'\circ\varphi_{\textnormal{cl}})^*$
envoie $\textnormal{ev}_{S}\in \CC[\, {\cal U}^{(n)}]$ 
sur la fonction $\textnormal{ev}\in\CC[x_i^{\pm 1},y_i^{\pm 1}]$
qui {\`a} tout polyn{\^o}me $P\in\CC[x_i^{\pm 1},y_i^{\pm 1}]$
associe son terme constant.
Donc, d'apr{\`e}s (\ref{coupla1}),
on obtient le couplage non-d{\'e}g{\'e}n{\'e}r{\'e}e suivant~:
\begin{equation}\label{coupla2}
\begin{array}{rcl}
U\bm {}_{/\!{}_{\displaystyle 
I_{\textnormal{cl}}^{(n)}}}\times
U^{(n)}&\longrightarrow&\CC\\
({\bar x},P)&\longmapsto&\textnormal{ev}(x.P)
\end{array}
\end{equation}
avec un plongement d'alg{\`e}bres~:
$$
U^{(n)}
\hookrightarrow {\Bigl( U\bm {}_{/\!{}_{\displaystyle 
I_{\textnormal{cl}}^{(n)}}}\Bigr)}^{*}
$$
tel que 
${\Bigl( U\bm {}_{/\!{}_{\displaystyle 
I_{\textnormal{cl}}^{(n)}}}\Bigr)}^{*}$
soit la compl{\'e}tion formelle
de $U^{(n)}$. Autrement dit, $U^{(n)}$ est le module coinduit~:
\begin{equation}\label{coindclass}
U^{(n)}\simeq {\bigl( \CC\otimes_{I_{\textnormal{cl}}^{(n)}}
U\bm\bigr)}^*
\end{equation}

Ceci cl{\^o}ture le cas classique. Passons au cas quantique.

\section{D{\'e}monstration de T2}

\noindent
Elle d{\'e}coule de la proposition suivante qui est une
g{\'e}n{\'e}ralisation de la proposition \ref{cl2}~:

\begin{prop}
\label{fracq2}
Soit $q$ une ind{\'e}termin{\'e}e, et
$A_n=\CC[q,q^{-1}]$.
Alors, dans ${A\lbrace\lbrace t_1,\ldots,t_n\rbrace\rbrace}/I_n$,
o{\`u} $I_n$ est l'id{\'e}al engendr{\'e} par les {\'e}l{\'e}ments
$t_i t_{i+1}-qt_{i+1}t_i$ pour $i\in\lbrace 1,\ldots,n\rbrace$
et les $t_it_j-t_jt_i$ pour 
$\arrowvert i-j\arrowvert\geq 2$, on a~:

   \begin{multline*}
\biggl (1-\Bigl (1-\bigl (1-\ldots (1-t_n)^{-1}t_{n-1}\bigr )^{-1}
\ldots t_2\Bigr )^{-1}t_1 \biggr )^{-1}\\
     =\displaystyle\sum_{\alpha_1,\ldots,\alpha_n}
F_q(\alpha_1,\ldots,\alpha_n)t_n^{\alpha_n}\ldots t_1^{\alpha_1}.
   \end{multline*}
\end{prop}

 \begin{dem}
Pour la d{\'e}finition de $A_n={A\{\{ t_1,\ldots,t_n \}\}}/I_n$,
voir l'appendice. 

\noindent Dans le cas classique o{\`u} $q=1$, on retrouve la proposition
\ref{cl2}.
\smallskip

\noindent D{\'e}montrons la propri{\'e}t{\'e} par r{\'e}currence sur $n$.
\smallskip

\noindent $\bullet$ Si $n=1$, le r{\'e}sultat est clair.
\smallskip

\noindent $\bullet$
Supposons la propri{\'e}t{\'e} vraie au rang $n$.
Notons pour tout $k,\, v_k$ la valuation sur $A_k$ correspondant
{\`a} la graduation telle que 
$\forall i,\, \deg{t_i}=1$, et
$F_{k}$ la fraction en question. 
Clairement, 
l'application~:
  $$
  \iota_{n}:\quad\begin{array}{rcl}
  A_{n-1}&\longrightarrow&A_{n+1}\\
  t_i&\longmapsto&t_{i+1}
  \end{array}
  $$
est un plongement d'alg{\`e}bres valu{\'e}es.
Par cons{\'e}quent, l'hypoth{\`e}se de\break
r{\'e}currence montre que~: 
  $$
  v_{n+1}\bigl(\iota_{n}(F_{n})\bigr)=v_{n}(F_{n})\geq 0.
  $$
D'o{\`u}~: 
  $$
  v_{n+1}\bigl(\iota_{n}(F_{n})t_1\bigr)\geq 1.
  $$
D'apr{\`e}s l'appendice, il s'ensuit que
$1-\iota_{n}(F_{n})t_1$ est inversible dans $A_{n+1}$, 
ce qui prouve l'existence de
$F_{n+1}$.

\medskip
\noindent
Posons~:
  $$
   v_1=t_1,\ldots, v_{n-1}=t_{n-1},\, v_n=(1-t_{n+1})^{-1}t_n.
  $$
Alors, 
  $$
   \forall i,j,\quad \arrowvert i-j\arrowvert \geq 2\Longrightarrow
   v_iv_j=v_jv_i.
  $$
De plus, 

       \begin{align*}
     v_{n-1}v_n&=t_{n-1}(1-t_{n+1})^{-1}v_n\\
               &=(1-t_{n+1})^{-1}t_{n-1}t_n\\
               &=q(1-t_{n+1})^{-1}t_nt_{n-1}\\
               &=qv_nv_{n-1}.
       \end{align*}

Donc, en appliquant l'hypoth{\`e}se de r{\'e}currence,

       \begin{align*}
            F_{n+1}&:=
\biggl (1-\Bigl (1-\bigl (1-\ldots (1-t_{n+1})^{-1}t_{n}\bigr )^{-1}
\ldots t_2\Bigr )^{-1}t_1 \biggr )^{-1}\\
             &=
\biggl (1-\Bigl (1-\bigl (1-\ldots (1-v_n)^{-1}v_{n-1}\bigr )^{-1}
\ldots v_2\Bigr )^{-1}v_1 \biggr )^{-1}\\
             &=\displaystyle\sum_{\alpha_1,\ldots,\alpha_n}
F_q(\alpha_1,\ldots,\alpha_n)v_n^{\alpha_n}\ldots v_1^{\alpha_1}.
       \end{align*}

On a 
$$
v_n^{\alpha_n}=\bigl [(1-t_{n+1})^{-1}t_n\bigr ]^{\alpha_n},
$$
et, par r{\'e}currence sur $k$, on montre que 
$$
\bigl [(1-t_{n+1})^{-1}t_n\bigr ]^{k}=
(1-t_{n+1})^{-1}\ldots (1-q^{k-1}t_{n+1})^{-1}t_n^{k}.$$
Ceci vient de la relation $$(1-t_{n+1})\, t_n^{-k}
=t_n^{-k}(1-q^k t_{n+1}).
$$
Donc, 
 
      \begin{align*}
          F_{n+1}&=
\biggl (1-\Bigl (1-\bigl (1-\ldots (1-t_{n+1})^{-1}t_{n}\bigr )^{-1}
\ldots t_2\Bigr )^{-1}t_1 \biggr )^{-1}\\
           &=\displaystyle\sum_{\alpha_1,\ldots,\alpha_n}
F_q(\alpha_1,\ldots,\alpha_n)
\bigl [(1-t_{n+1})^{-1}t_n\bigr ]^{\alpha_n}\ldots t_1^{\alpha_1}\\
           &=\displaystyle\sum_{\alpha_1,\ldots,\alpha_n}
F_q(\alpha_1,\ldots,\alpha_n)
(1-t_{n+1})^{-1}\ldots (1-q^{\alpha_n -1}t_{n+1})^{-1}t_n^{\alpha_n}
\ldots t_1^{\alpha_1}.
      \end{align*}

On utilise ensuite la relation classique~:
$$
\displaystyle\prod_{s=0}^{n-1}(1-q^st)^{-1}
=\displaystyle\sum_{k\geq 0}\C{n+k-1}{k}t^{k},
$$
pour obtenir~:
  
      \begin{align*}
      F_{n+1}&=\displaystyle\sum_{\alpha_1,\ldots,\alpha_n}
F_q(\alpha_1,\ldots,\alpha_n)
\Bigl
(\displaystyle\sum_{\alpha_{n+1}}\C{\alpha_n+\alpha_{n+1}-1}{\alpha_n}
t_{n+1}^{\alpha_{n+1}}\Bigr )t_n^{\alpha_n}\ldots t_1^{\alpha_1}\\
       &=\displaystyle\sum_{\alpha_1,\ldots,\alpha_{n+1}}
F_q(\alpha_1,\ldots,\alpha_{n+1})t_{n+1}^{\alpha_{n+1}}\ldots
t_1^{\alpha_1}.
      \end{align*}
  \end{dem}

De la m{\^e}me fa\c con, on montrerait la proposition suivante~:
    \begin{prop}
\label{pasdemprop}
Soit $q$ une ind{\'e}termin{\'e}e, et
$A=\CC[q,q^{-1}]$.
Alors, dans ${A\lbrace\lbrace t_1,\ldots,t_n\rbrace\rbrace}/J_{n}$,
o{\`u} $J_{n}$ est l'id{\'e}al engendr{\'e} par les {\'e}l{\'e}ments
$t_{i+1} t_{i}-qt_{i}t_{i+1}$ pour $i\in\lbrace 1,\ldots,n\rbrace$
et les $t_it_j-t_jt_i$ pour 
$\arrowvert i-j\arrowvert\geq 2$, on a~:

        \begin{multline*}
\biggl (1-t_1\Bigl (1-t_2\bigl (\ldots (1-t_n)^{-1}\bigr )^{-1}
\ldots\Bigr )^{-1}\biggr )^{-1}\\
=\displaystyle\sum_{\alpha_1,\ldots,\alpha_n}
F_q(\alpha_1,\ldots,\alpha_n)t_1^{\alpha_1}\ldots t_n^{\alpha_n}. 
        \end{multline*}
    \end{prop}

\section{Action de $U_q(\bm)$ sur les $u_i,1\leq i\leq 2n$}

\noindent
Tout est r{\'e}sum{\'e} dans la proposition suivante~:
   \begin{prop}
   \label{calculs+}
Pour tout $i\in\lbrace 1,\ldots,2n\rbrace$,
on a~:
       \begin{align}
       \label{equabui1} k.u_i&=q u_i,\\
       \label{equabui2} e_{+}.u_i&=\displaystyle\sum\limits_{r+s=i}u_r u_s,\\
       \label{equabui3} e_{-}.u_i&=q^{-1}\delta_{i}^{1}.
       \end{align}
   \end{prop}
\begin{dem}

La premi{\`e}re {\'e}galit{\'e} est triviale puisque $\deg{u_i}=1$.

\noindent
Pour la d{\'e}monstration des deux autres {\'e}galit{\'e}s nous aurons
besoin en particulier du lemme suivant~:

     \begin{lem}
     \label{lem1}
Pour tous $N,a_1,\ldots,a_{2N},k$ tels que 
$k\in\lbrace 1,\ldots ,N+1\rbrace$, on pose~:
$$
Y_N(k,a_1,\ldots,a_{2N})=k-1+\displaystyle\sum\limits_{j=1}^{2k-2}
(-1)^j a_j
$$ 
et
   \begin{multline}
   \label{defphi} \Phi_{N} (k,a_1,\ldots,a_{2N})=
q^{Y_N(k,a_1,\ldots,a_{2N})} 
(q^{-(a_{2k-2} +1)}-q^{-a_{2k-1}})\\
\times\, F_q(a_1 -1,a_2 +1,\ldots,a_{2k-2}+1,a_{2k-1},\ldots,a_{2N}).
   \end{multline}
avec la convention que $a_0=a_{2N+1}=0$.
Alors, 
$$\displaystyle\sum\limits_{k=1}^{N+1}
\Phi_{N}(k,a_1,\ldots,a_{2N})=0.
$$

     \end{lem}

\noindent
\begin{pr}
La d{\'e}monstration se fait par r{\'e}currence sur $N$.

\medskip
\noindent
$\bullet$ Si $N=1$, et si $(a_1,a_2)\in{\NN}^{2}-\lbrace 0\rbrace$,

   \begin{multline*}
\displaystyle\sum\limits_{k=1}^{2}\Phi(k,a_1,a_2)=
(q^{-1}-q^{-a_1})F_q(a_1,a_2)\\
+q^{-a_1+a_2+1}(q^{-(a_2+1)}-1)F_q(a_1-1,a_2+1)
   \end{multline*}
soit

    \begin{multline*}
\displaystyle\sum\limits_{k=1}^{2}\Phi(k,a_1,a_2)=
(q-1)q^{-a_1}\times\\
\times\left( \q{a_1-1}F_q(a_1,a_2)-\q{a_2+1}
F_q(a_1-1,a_2+1)\right).
    \end{multline*}

L'assertion d{\'e}coule alors de l'identit{\'e} suivante:
$$
\forall\; (x,y)\in{\NN}^2-\lbrace 0\rbrace,\quad
\q{x-1}F_q(x,y)=\q{y+1}F_q(x-1,y+1).
$$

\bigskip
\noindent
$\bullet$ Supposons la propri{\'e}t{\'e} vraie 
jusqu'au rang $N-1$, avec $N\ge 2$.\\
Soit $(a_1,\cdots,a_{2N})\in{\NN}^{2N}-\lbrace 0\rbrace$.
On a 
\begin{multline}\label{4.12}
\displaystyle\sum\limits_{k=1}^{N+1}
\Phi_{N}(k,a_1,\ldots,a_{2N})=\\
\displaystyle\sum\limits_{k=1}^{N-1}
\Phi_{N}(k,a_1,\ldots,a_{2N})
+\Phi_{N}(N,a_1,\ldots,a_{2N})+\Phi_{N}(N+1,a_1,\ldots,a_{2N})
\end{multline}
De plus, on note que pour $k\in\lbrace 1,\ldots ,N-1\rbrace$, on a
$$
\Phi_{N}(k,a_1,\ldots,a_{2N})=\Phi_{N-1}(k,a_1,\ldots,a_{2N-2})
F_q(a_{2N-2},a_{2N-1},a_{2N})
$$
Donc, l'hypoth{\`e}se de r{\'e}currence nous montre que
$$
\displaystyle\sum\limits_{k=1}^{N-1}
\Phi_{N}(k,a_1,\ldots,a_{2N})=
-\Phi_{N-1}(N,a_1,\ldots,a_{2N-2})F_q(a_{2N-2},a_{2N-1},a_{2N})
$$
d'o{\`u} 
\begin{equation}\label{4.13}
\begin{split}
\displaystyle\sum\limits_{k=1}^{N-1}
\Phi_{N}(k,a_1,\ldots,a_{2N})&=q^{Y_N}F_q(a_1 -1,a_2 +1,\ldots,a_{2N-2}+1)\\
&\times\, (1-q^{-(a_{2N-2}+1)})\,F_q(a_{2N-2},a_{2N-1},a_{2N})
\end{split}
\end{equation}
Par ailleurs, on a {\'e}galement
\begin{multline}\label{4.14}
\Phi_{N}(N,a_1,\ldots,a_{2N})=q^{Y_N} F_q(a_1-1,\cdots,a_{2N-2}+1)\\
\times\,(q^{-(a_{2N-2}+1)}-q^{-a_{2N-1}})\,
F_q(a_{2N-2}+1,a_{2N-1},a_{2N})
\end{multline}
et 
\begin{multline}\label{4.15}
\Phi_{N}(N+1,a_1,\ldots,a_{2N})=q^{Y_N} F_q(a_1 -1,a_2 +1,\ldots,a_{2N-2}+1)\\
\times\, q^{-a_{2N-1}+a_{2N}+1}(q^{-(a_{2N}+1)}-1)\,
F_q(a_{2N-2}+1,a_{2N-1}-1,a_{2N}+1).
\end{multline}
Donc, en rassemblant~(\ref{4.12}), (\ref{4.13}), (\ref{4.14})
et~(\ref{4.15}), on obtient
$$
\displaystyle\sum\limits_{k=1}^{N+1}
\Phi_{N}(k,a_1,\ldots,a_{2N})=q^{Y_N} F_q(a_1 -1,a_2
+1,\ldots,a_{2N-2}+1)\times A
$$
avec 
\begin{multline*}
A=(1-q^{-(a_{2N-2}+1)})\,F_q(a_{2N-2},a_{2N-1},a_{2N})\\
+(q^{-(a_{2N-2}+1)}-q^{-a_{2N-1}})
F_q(a_{2N-2}+1,a_{2N-1},a_{2N})\\
+q^{-a_{2N-1}+a_{2N}+1}(q^{-(a_{2N}+1)}-1)\,
F_q(a_{2N-2}+1,a_{2N-1}-1,a_{2N}+1)\\
\end{multline*}
Si $a_{2N-1}=0$, le r{\'e}sultat est clair.\\
Sinon, on a
$$
q^{-(a_{2N}+1)}-1=-q^{-(a_{2N+1})}\q{a_{2N}+1}(q-1)
$$
et
$$
\q{a_{2N}+1}F_q(a_{2N-1}-1,a_{2N}+1)=\q{a_{2N-1}-1}F_q(a_{2N-1},a_{2N}).
$$
Posons $x=a_{2N-2}$, $y=a_{2N-1}$ et $z=a_{2N}$.
Alors, $A=F_q(y,z)B$, avec
\begin{equation*}
\begin{array}{r@{}c@{}l}
B=&\,(1-q^{-(x+1)})\C{x+y-1}{y}&+\,(q^{-(x+1)}-q^{-y})\C{x+y}{y}\\
&&-\, q^{-y}(q^{y-1}-1)\C{x+y-1}{y-1}
\end{array}
\end{equation*}
En utilisant la relation
$$
\C{x+y}{y}=q^{y}\C{x+y-1}{y}+\C{x+y-1}{y-1},
$$
on obtient
$$
B=q^{-(x+1)}(q^y-1)\C{x+y-1}{y}-q^{-(x+1)}(q^x-1)\C{x+y-1}{y-1}.
$$
Le lemme d{\'e}coule alors de la simple {\'e}galit{\'e}:
$$
\q{y}\C{x+y-1}{y}=\q{x}\C{x+y-1}{y-1}.
$$
\end{pr}

\medskip
\noindent
{\bf D{\'e}monstration de l'{\'e}galit{\'e} (\ref{equabui3}):}

\medskip
\noindent
\begin{pr}
\smallskip

\noindent
Soit $k\in\lbrace 1,\ldots,2n\rbrace$.
Nous allons calculer $e_{-}.u_k$.

\noindent
$\bullet$ Si $k=1$, le calcul est imm{\'e}diat.\\
$\bullet$ Supposons $k\geq 2$. On a~: 
$$
e_{-} .u_k=\displaystyle{1\over q-1}
\displaystyle\sum\limits_{\arrowvert\alpha\arrowvert =k-1}
F_q(\alpha_{2n-1},\ldots,\alpha_1) U({\underline \alpha}),
$$
avec 
\begin{multline*}
U({\underline \alpha})={\Sigma}^-
(x_1 y_1)^{-\alpha_1}\ldots (x_n y_n)^{-\alpha_{2n-1}}y_n^{-1}\\
- q^{-1}(x_1 y_1)^{-\alpha_1}\ldots (x_n y_n)^{-\alpha_{2n-1}}
y_n^{-1}{\Sigma}^- .
\end{multline*}
En adoptant la convention $\alpha_{2n}=1$, et en utilisant les r{\`e}gles
de commutation entre les $x_i$ et les $y_j$, on a d'abord~:
\begin{multline*}
U({\underline \alpha})=
\displaystyle\sum\limits_{j=1}^{n}
\bigl ( q^{-\alpha_{2j-1}}-q^{-\alpha_{2j}}\bigr )
(x_1 y_1)^{-\alpha_1}\ldots (x_j y_j)^{-\alpha_{2j-1}}\times\\
y_j (y_j x_{j+1})^{-\alpha_{2j}}\ldots (x_n y_n)^{-\alpha_{2n-1}} y_n^{-1}.
\end{multline*}
Puis, toujours {\`a} l'aide des relations de commutation,
\begin{multline*}
U({\underline{\alpha}})=\sum\limits_{j=1}^{n}
q^{X_{j,{\underline{\alpha}}}}
(q^{-\alpha_{2j-1}}-q^{-\alpha_{2j}})
(x_1 y_1)^{-\alpha_1}\ldots (x_j y_j )^{-\alpha_{2j-1}}\times\\
(y_j x_{j+1})^{-(\alpha_{2j}-1)} (x_{j+1} y_{j+1})^{-(\alpha_{2j+1}+1)}
\ldots (x_n y_n)^{-(\alpha_{2n-1}+1)}
\end{multline*}
avec pour notation $\lbrace x\rbrace =1+2+\ldots+x$, et~:
\begin{align*}
X_{j,{\underline{\alpha}}}&=
\lbrace\alpha_1\rbrace
+\ldots+\lbrace\alpha_{2n-1}\rbrace\\
&-\bigl (
\lbrace\alpha_1\rbrace+\ldots+\lbrace\alpha_{2j-1}\rbrace
+\lbrace\alpha_{2j} -1\rbrace+\lbrace\alpha_{2j+1} +1\rbrace
+\ldots\lbrace\alpha_{2n-1} +1\rbrace\bigr )\\
&=\alpha_{2j}-(\alpha_{2j+1}+1)+\ldots-(\alpha_{2n-1} +1)\\
&=j-n+\sum\limits_{k=2j}^{2n-j} (-1)^k \alpha_k
\end{align*}
Donc, $e_{-} .u_k$ est une combinaison lin{\'e}aire de termes
de la forme
$$(x_1 y_1)^{-\beta_1}\ldots (x_n y_n)^{-\beta_{2n-1}}$$
tels que 
$$\beta_1+\ldots+\beta_{2n-1}=k-1.$$
Le coefficient correspondant est obtenu en regroupant 
tous les entiers\break
$j,\alpha_1,\ldots,\alpha_{2n-1}$ tels que~:
$$
\left\{
\begin{array}{l}
\alpha_1=\beta_1\\
\vdots\\
\alpha_{2j-1}=\beta_{2j-1}\\
\alpha_{2j}-1=\beta_{2j}\\
\vdots\\
\alpha_{2n-1}+1=\beta_{2n-1}\\
\end{array}
\right.
$$
On a~: 
\begin{align*}
j-n+\sum\limits_{k=2j}^{2n-1}
(-1)^k \alpha_k&=
j-n+\sum\limits_{k=2j}^{2n-1}
(-1)^k (\beta_k +(-1)^k)\\
&=n-j+\sum\limits_{k=2j}^{2n-1}
(-1)^k \beta_k .
\end{align*}
Par suite, ce coefficient $A$ vaut
\begin{multline*}
A=\displaystyle{1\over q-1}
\sum\limits_{j=1}^{n}
q^{X'_j}
(q^{-\beta_{2j-1}}-q^{-(\beta_{2j}+1)})\times\\
F_q(\beta_{2n-1}-1,\beta_{2n-2}+1,\ldots,\beta_{2j+1}-1,\beta_{2j}+1,
\beta_{2j-1},\ldots,\beta_1),
\end{multline*}
avec 
$$
X'_j=n-j+\sum\limits_{k=2j}^{2n-1}(-1)^k \beta_k .
$$
Posons $\gamma_j=\beta_{2n-j}$ pour $j\in\lbrace
1,\ldots,2n-1\rbrace$,
et $\gamma_0=0$.
Alors,
$$
X'_j=n-j+\sum\limits_{k=1}^{2(n-j)}(-1)^k \gamma_k 
$$
et
\begin{multline*}
A=\displaystyle{1\over q-1}
\sum\limits_{j=1}^n 
q^{X'_j}
(q^{-\gamma_{2(n-j)+1}}-q^{-(\gamma_{2(n-j)} +1)})\times\\
F_q(\gamma_1-1,\gamma_2 +1,\ldots,\gamma_{2(n-j)-1}-1,
\gamma_{2(n-j)}+1,
\gamma_{2(n-j)+1},\ldots,\gamma_1).
\end{multline*}
Soit, apr{\`e}s le changement de variables $s=n-j$,
\begin{multline*}
A=-\displaystyle{1\over q-1}
\sum\limits_{s=0}^{n-1}
q^{Y_{s+1}}\bigl ( q^{-(\gamma_{2s}+1)}-q^{-(\gamma_{2s+1})}\bigr )\times\\
F_q(\gamma_1 -1,\gamma_2 +1,\ldots,\gamma_{2s-1} -1,\gamma_{2s}
+1,\gamma_{2s+1},\ldots,\gamma_{2n-1}),
\end{multline*}
avec 
$$
Y_{s+1}=s+\sum\limits_{t=1}^{2s} (-1)^t \gamma_t .
$$
Supposons que $\gamma_{2n-1}\geq 2$.
Alors, pour tout $s\in\lbrace 0,\ldots,n-1\rbrace$, on a~:
$$
F_q(\gamma_1 -1,\gamma_2 +1,\ldots,\gamma_{2s-1} -1,\gamma_{2s} +1,
\gamma_{2s+1},\ldots,\gamma_{2n-1})=0.$$
En effet, sinon, il existe $s\in\lbrace 0,\ldots,n-1\rbrace$ tel que
$$
\left\{
\begin{array}{c}
\gamma_1 -1\geq 1\\
\gamma_3 -1\geq 1\\
\vdots\\
\gamma_{2s-1} -1\geq 1\\
\gamma_{2s+1}\geq 1\\
\vdots\\
\gamma_{2n-2}\geq 1
\end{array}
\right.
$$
Donc, 
$$
\left\{
\begin{array}{l}
\gamma_1+\ldots+\gamma_{2s-1}\geq 2s\\
\gamma_{2s+1}+\ldots+\gamma_{2n-2}\geq 2n-2-(2s+1)+1=2n-2s-2\\
\gamma_{2n-1}\geq 2
\end{array}
\right.
$$ D'o{\`u} 
$$
k-1=\gamma_1+\ldots+\gamma_{2n-1}\geq 2n
$$
ce qui est impossible, vu que $k\leq 2n$.
Par suite, $\gamma_{2n-1}\leq 1$, et, en posant $\gamma_{2n}=0$,
on a~: 
$$F_q(\gamma_1 -1,\gamma_2 +1,\ldots,\gamma_{2n-1}-1,\gamma_{2n} +1)=0
$$
et 
\begin{multline*}
F_q(\gamma_1 -1,\gamma_2 +1,\ldots,\gamma_{2s-1} -1,\gamma_{2s} +1,
\gamma_{2s+1},\ldots,\gamma_{2n-1})\\
=F_q(\gamma_1 -1,\gamma_2 +1,\ldots,\gamma_{2s-1} -1,\gamma_{2s} +1,
\gamma_{2s+1},\ldots,\gamma_{2n-1},\gamma_{2n}).
\end{multline*}
Donc,
\begin{multline*}
A=-\displaystyle{1\over q-1}
\sum\limits_{s=0}^n 
q^{Y_{s+1}}\bigl ( q^{-(\gamma_{2s}+1)}-q^{-(\gamma_{2s+1})}\bigr )\times\\
F_q(\gamma_1 -1,\gamma_2 +1,\ldots,\gamma_{2s-1} -1,\gamma_{2s}
+1,\gamma_{2s+1},\ldots,\gamma_{2n-1},\gamma_{2n}).
\end{multline*}
On utilise ensuite le lemme \ref{lem1} pour conclure: $A=0$.
\end{pr}

\bigskip
\noindent
Pour la d{\'e}monstration de l'{\'e}galit{\'e} (\ref{equabui2}),
nous aurons besoin des deux lemmes suppl{\'e}mentaires suivants~:

          \begin{lem}
          \label{lem2}
Soient $t_1,\ldots,t_n$ les variables $q$-commutatives
d{\'e}finis dans la proposition \ref{pasdemprop} et $B$ l'alg{\`e}bre
qu'elles engendrent. Alors,
\begin{enumerate}
     \item   L'alg{\`e}bre $B$ est un $\CC[q,q^{-1}]$-module libre
             dont une base est donn{\'e}e par la famille
             $t_1^{\a_1}\ldots t_n^{\a_n}$, avec $\forall i,\,
             \a_i\in\NN$.
     \item   Pour tout $k\in\lbrace 1,\ldots,n\rbrace$, l'application
             $$
                 \begin{array}{rcl}
                 \phi_k:\quad B&\longrightarrow&B\\
                 t_i&\longmapsto&
                   \begin{cases}
                    t_i,&\textnormal{ si }i\not=k\,;\\
                    qt_k,&\textnormal{ si }i=k\,. 
                   \end{cases}
                 \end{array}
             $$
             est un morphisme d'alg{\`e}bres graudu{\'e}es. 
     \item   Pour tout $u\in B$, on posera
             $$
               \phi_k(u)=u(t_1,\ldots,q t_k,\ldots,t_n),
             $$
             et
             $$
               \partial_k u ={{1}\over {q-1}}\displaystyle\lbrack
        u(t_1,\ldots,q
        t_k,\ldots,t_n)-u(t_1,\ldots,t_k,\ldots,t_n)\rbrack .
             $$
                Alors, pour tous  $u$ et $v$ dans $B$, on a~:
                \begin{align}
                \label {r1} \partial_k 1    &= 0\\
                \label {r2} \partial_k (t_1^{\alpha_1}\ldots
                t_n^{\alpha_n})
                                            &=[\alpha_k]
                t_1^{\alpha_1},\ldots,t_n^{\alpha_n}\\
                \label {r3} \partial_k(u+v) &=\partial_k(u)+\partial_k(v)\\
                \label {r4} \partial_k(uv)  &=(\partial_k u)v+
                u\bigl (t_1,\ldots,qt_k,\ldots,t_n\bigr )\partial_k v.
                \end{align}
            De plus, si $u$ est inversible, alors
                \begin{equation}
                \label{r5}
                  \partial_k(u^{-1}) = 
                  -u^{-1}\bigl ( t_1,\ldots,qt_k,\ldots,t_n\bigr )
                  (\partial_k u)u^{-1}\bigl ( t_1,\ldots, t_k, \ldots,t_n).
                \end{equation}
\end{enumerate}

\end{lem}

\noindent
\begin{dem}
Pour le point 1, on utilise l'appendice avec $q$ chang{\'e} en
$q^{-1}$. Pour le reste, il suffit d'effectuer les calculs.
En particulier, les {\'e}galit{\'e}s \ref{r1} et \ref{r2}
montrent que l'expression $\partial_k u$ a bien un sens, i.e.,
que l'on peut l{\'e}gitimement diviser par $q-1$.
\end{dem}

\noindent
\begin{lem}\label{lem3}
 Pour toute suite d'entiers ${\underline \alpha}=
(\alpha_1,\ldots,\alpha_n)$, on a~:
$$
[\alpha_1+1]
F_q(\alpha_1,\ldots,\alpha_n)=\displaystyle\sum\limits_{{\underline \beta},
{\underline \gamma}
\atop{
{\underline \beta}+{\underline \gamma}={\underline \alpha}
}}
q^{\beta_1+\beta_2\gamma_1+\ldots +\beta_n\gamma_{n-1}}
F_q({\underline \beta})F_q({\underline \gamma}).
$$
\end{lem}

\noindent
\begin{pr}

\noindent
Dans le cas classique, en posant
$$
f(t_1,\ldots,t_n)=
\cfrac{1}{1-\cfrac{t_1}{1-\cfrac{t_2}{\ddots\cfrac{t_{n-1}}{1-t_n}}}}
$$
la relation signifie simplement que $\partial_1 (t_1 f)=f^2.$

\noindent
Dans le cas g{\'e}n{\'e}ral, posons
$$
f:=\biggl (1-t_1\Bigl (1-t_2\bigl (\ldots (1-t_n)^{-1}\bigr )^{-1}
\ldots\Bigr )^{-1}\biggr )^{-1}=(1-t_1g)^{-1}=h^{-1},
$$
avec $h=1-t_1g$.
Alors, par~(\ref{r5}), 
$$
\partial_1 f=-f\bigl (qt_1,t_2,\ldots, t_n\bigr )(\partial_1 h)
f(t_1,\ldots,t_n).
$$
D'apr{\`e}s le lemme et du fait que $g$ ne d{\'e}pend
pas de $t_1$, on a~:
\begin{align*}
\partial_1 h&=-(\partial_1 t_1)g-(qt_1)(\partial_1 g)\\
&=-t_1 g\\
&=h-1.
\end{align*}
Donc, 
\begin{align*}
\partial_1 f&=-f\bigl (qt_1,t_2,\ldots, t_n\bigr )(h-1)f\\
&=-f\bigl (qt_1,t_2,\ldots, t_n\bigr )
\bigl ( 1-f(t_1,\ldots,t_n)\bigr )\\
&=-f\bigl (qt_1,t_2,\ldots, t_n\bigr )+
f\bigl (qt_1,t_2,\ldots, t_n\bigr )f(t_1,\ldots,t_n).
\end{align*}
D'o{\`u},
\begin{equation}\label{r6}
t_1(\partial_1 f)+t_1 f\bigl (qt_1,t_2,\ldots, t_n\bigr )=
t_1 f\bigl (qt_1,t_2,\ldots, t_n\bigr )f(t_1,\ldots,t_n).
\end{equation}
Par ailleurs, d'apr{\`e}s la proposition~{\ref{pasdemprop}}, on a~:
$$
f(t_1,\ldots,t_n)=\displaystyle\sum_{\alpha_1,\ldots,\alpha_n}
F_q(\alpha_1,\ldots,\alpha_n)t_1^{\alpha_1}\ldots t_n^{\alpha_n}. 
$$
Donc, l'application de d{\'e}rivation $\partial_1$ {\'e}tant continue
pour la graduation d{\'e}finie par la valuation sur $B$, on a~:
\begin{align}
\notag t_1(\partial_1 f)+t_1 f\bigl (qt_1,t_2,\ldots, t_n\bigr )&=
t_1 \displaystyle\sum_{\alpha_1,\ldots,\alpha_n}[\alpha_1]
F_q(\alpha_1,\ldots,\alpha_n)t_1^{\alpha_1}\ldots t_n^{\alpha_n}\\
\notag &+\displaystyle\sum_{\alpha_1,\ldots,\alpha_n}q^{\alpha_1}
F_q(\alpha_1,\ldots,\alpha_n)t_1^{\alpha_1 +1}\ldots t_n^{\alpha_n}\\
\label{r7} &=\displaystyle\sum_{\alpha_1,\ldots,\alpha_n}[\alpha_1 +1]
F_q(\alpha_1,\ldots,\alpha_n)t_1^{\alpha_1 +1}\ldots t_n^{\alpha_n}.
\end{align}
D'autre part, 
\begin{multline*}
f\bigl (qt_1,t_2,\ldots, t_n\bigr )f(t_1,\ldots,t_n)=
\displaystyle\Bigl\lbrace\displaystyle\sum\limits_{{\underline
\beta}}F_q({\underline{\beta}})q^{\beta_1}t_1^{\beta_1}\ldots
t_n^{\beta_n}\Bigr\rbrace\times\\
\times\displaystyle\Bigl\lbrace\displaystyle\sum\limits_{{\underline
\gamma}}F_q({\underline{\gamma}})t_1^{\gamma_1}\ldots
t_n^{\gamma_n}\Bigr\rbrace
\end{multline*}
Soit,
$$
f\bigl (qt_1,t_2,\ldots, t_n\bigr )f(t_1,\ldots,t_n)=
\displaystyle\sum\limits_{{\underline
\beta,\underline \gamma}}F_q({\underline \beta}) F_q({\underline \gamma})
q^{\beta_1}t_1^{\beta_1}\ldots t_n^{\beta_n}
t_1^{\gamma_1}\ldots t_n^{\gamma_n}.
$$
Mais, vu les relations de commutativit{\'e} entre les $t_i$, on a~: 
$$
t_{i+1}^m t_i^n=q^{mn}t_i^{n} t_{i+1}^m.
$$
Donc, 
\begin{multline*}
f\bigl (qt_1,t_2,,\ldots, t_n\bigr )f(t_1,\ldots,t_n)\\
=\displaystyle\sum\limits_{({\underline
\beta,\underline \gamma})}F_q({\underline \beta}) F_q({\underline \gamma})
q^{\beta_1}\times q^{\beta_2 \gamma_1}\times \ldots\times 
q^{\beta_n \gamma_{n-1}}
t_1^{\beta_1+\gamma_1}\ldots t_n^{\beta_n+\gamma_n}.
\end{multline*}
Soit,
\begin{multline}\label{r8}
t_1 f\bigl (qt_1,t_2,\ldots, t_n\bigr )f(t_1,\ldots,t_n)
=\\
\displaystyle\sum\limits_{{\underline \alpha}}\Bigl\lbrace
\displaystyle\sum\limits_{({\underline
\beta,\underline \gamma})\atop{\underline\beta
+\underline\gamma=\alpha}}q^{\beta_1+\beta_2\gamma_1+
\ldots \beta_n\gamma_{n-1}}
F_q({\underline \beta})F_q({\underline \gamma})\Bigr\rbrace
t_1^{\alpha_1 +1}\ldots t_n^{\alpha_n}.
\end{multline}
Il suffit ensuite d'identifier~(\ref{r7}) et~(\ref{r8}) {\`a} l'aide 
de~(\ref{r6}) pour obtenir le r{\'e}sultat.
\end{pr}

\medskip
\noindent
{\bf D{\'e}monstration de l'{\'e}galit{\'e} (\ref{equabui2}):}

\medskip
\noindent
Soit $k\in \lbrace 1,\ldots, 2n\rbrace$. On a~;
\begin{multline*}
e_{+} .u_k=\displaystyle{1\over q-1}\displaystyle\sum\limits_{
(\alpha_1,\ldots,\alpha_{2n-1})\atop{
\alpha_1+\ldots+\alpha_{2n-1}=k-1}}F_q(\alpha_{2n-1},\ldots,\alpha_1)
\times\\
\displaystyle\sum\limits_{i=1}^n 
x_i(x_1 y_1)^{-\alpha_1}\ldots (x_n y_n)^{-\alpha_{2n-1}}y_n^{-1}
-q(x_1 y_1)^{-\alpha_1}\ldots (x_n y_n)^{-\alpha_{2n-1}}y_n^{-1}x_i .
\end{multline*}
Les relations de commutation entre les $x_i$ et $y_i$ montrent que 
$x_i$ commute avec $(y_{j-1} x_j)^{-\alpha_{2j-2}}$
et $(x_j y_j)^{-\alpha_{2j-1}}$, pour $j\not= i$.
De plus, 
$$
\left\{
\begin{array}{l}
q y_n^{-1} x_i =x_i y_n^{-1}\\
x_i (y_{i-1} x_i)^{-\alpha_{2i-2}}=q^{-\alpha_{2i-2}} 
(y_{i-1} x_i)^{-\alpha_{2i-2}}x_i\\
(x_i y_i)^{-\alpha_{2i-1}} x_i =q^{-\alpha_{2i-1}}
x_i (x_i y_i)^{-\alpha_{2i-1}}
\end{array}
\right.
$$
Donc, 
\begin{multline*}
e_{+} .u_k=\displaystyle{1\over q-1}\displaystyle\sum\limits_{
(\alpha_1,\ldots,\alpha_{2n-1})\atop{
\alpha_1+\ldots+\alpha_{2n-1}=k-1}}F_q(\alpha_{2n-1},\ldots,\alpha_1)\times\\
\displaystyle\sum\limits_{i=1}^n 
( q^{-\alpha_{2i-2}}-q^{-\alpha_{2i-1}})
(x_1 y_1)^{-\alpha_1}\ldots (y_{i-1} x_i)^{-\alpha_{2i-2}} x_i
(x_i y_i)^{-\alpha_{2i-1}}\ldots\\
\ldots (x_n y_n)^{-\alpha_{2n-1}}y_n^{-1}
\end{multline*}
En notant $\lbrace x\rbrace =\displaystyle{x(x+1)\over 2}
=\sum\limits_{t=1}^x t$,
et en utilisant plusieurs fois les relations
$$
\left\{
\begin{array}{l}
(x_j y_j)^{-n} =q^{\lbrace n\rbrace} x_j^{-n} y_j^{-n}\\
(y_{j-1} x_j)^{-n}=q^{\lbrace n\rbrace} y_{j-1}^{-n} x_j^{-n}
\end{array}
\right.
$$ valables pour tous $j$ et $n$, on obtient~:
\begin{multline*}
e_{+} .u_k=\displaystyle{1\over q-1}\displaystyle\sum\limits_{
(\alpha_1,\ldots,\alpha_{2n-1})\atop{
\alpha_1+\ldots+\alpha_{2n-1}=k-1}}F_q(\alpha_{2n-1},\ldots,\alpha_1)
\times\\
\displaystyle\sum\limits_{i=1}^n 
q^{\lbrace\alpha_1\rbrace +\ldots +\lbrace\alpha_{2n-1}\rbrace
-\bigl (
\lbrace\alpha_1\rbrace +\ldots +\lbrace\alpha_{2i-2}\rbrace
+\lbrace\alpha_{2i-1} -1\rbrace +\lbrace\alpha_{2i} +1\rbrace
+\ldots +\lbrace\alpha_{2n-1} -1\rbrace\bigr )}\times\\
(q^{-\alpha_{2i-2}}-q^{-\alpha_{2i-1}})
(x_1 y_1)^{-\alpha_1}\ldots (y_{i-1} x_i)^{-\alpha_{2i-2}}\times\\
(x_i y_i)^{-(\alpha_{2i-1} -1)}
(y_i x_{i+1})^{-(\alpha_{2i} +1)}\ldots
(x_n y_n)^{-(\alpha_{2n-1}-1)} y_n^{-2}
\end{multline*}
Il s'ensuit que $e_{+}.u_k$ est une combinaison lin{\'e}aire de
termes de la forme 
$$(x_1 y_1)^{-\beta_1}\ldots 
(x_n y_n)^{-\beta_{2n-1}}y_n^{-2},$$
avec, 
$$\beta_1+\ldots+\beta_{2n-1}=\alpha_1+\ldots+\alpha_{2n-1}-1=k-2.$$
Le coefficient de $(x_1 y_1)^{-\beta_1}\ldots 
(x_n y_n)^{-\beta_{2n-1}}y_n^{-2}$ est obtenu pour chaque\break 
multi-indice
${\underline \alpha}=(\alpha_1,\ldots,\alpha_{2n-1})$ et chaque $i\in\lbrace
1,\ldots,n\rbrace $ tels que~:
$$
\left\{
\begin{array}{l}
\alpha_1=\beta_1\\
\vdots\\
\alpha_{2i-2}=\beta_{2i-2}\\
\alpha_{2i-1}=\beta_{2i-1} +1\\
\alpha_{2i}=\beta_{2i} -1\\
\vdots\\
\alpha_{2n-2}=\beta_{2n-2}-1\\
\alpha_{2n-1}=\beta_{2n-1}+1
\end{array}
\right.
$$
Par suite, 
\begin{multline*}
e_{+}.u_k=\displaystyle{1\over q-1}
\displaystyle\sum\limits_{\arrowvert {\underline\beta}\arrowvert=k-2}
\displaystyle\sum\limits_{i=1}^n
q^{n-i+1+\sum\limits_{t=2i-1}^{2n-1}
(-1)^{t+1} \beta_t}\bigl
(q^{-\beta_{2i-2}}-q^{-(\beta_{2i-1}+1)}\bigr )\times\\
F_q(\beta_{2n-1}+1,
\beta_{2n-2}-1,\ldots,\beta_{2i}-1,\beta_{2i-1}+1,
\beta_{2i-2},\ldots,\beta_1)\\
(x_1 y_1)^{-\beta_1}\ldots (x_n y_n)^{-\beta_{2n-1}}y_n^{-2},
\end{multline*}
avec la notation que 
$\arrowvert {\underline\alpha}\arrowvert =\sum\alpha_i$.

\noindent
Le lemme~{\ref{lem1}} appliqu{\'e} {\`a} $a_1=2$ et
$\forall i\in\lbrace 2,\ldots,2n\rbrace,\, a_i=\a_{2n-i+1}$
entra{\^\i}ne~:
\begin{multline*}
(q^{\alpha_{2n-1}+1}-1)F_q(\alpha_{2n-1},\ldots,\alpha_1)\\
=\displaystyle\sum\limits_{i=1}^n
q^{n-i+1+\sum\limits_{t=2i-1}^{2n-1}
(-1)^{t+1} \alpha_t}
\bigl
(q^{-\alpha_{2i-2}}-q^{-(\alpha_{2i-1}+1)}\bigr )\times\\
F_q(\alpha_{2n-1}+1,
\alpha_{2n-2}-1,\ldots,\alpha_{2i}-1,\alpha_{2i-1}+1,
\alpha_{2i-2},\ldots,\alpha_1)
\end{multline*}
Donc, 
\begin{multline*}
e_{+}.u_k=
\displaystyle\sum\limits_{
({\underline\alpha})\atop{\arrowvert
{\underline\alpha}\arrowvert =k-2}}
[\alpha_{2n-1} +1] F_q(\alpha_{2n-1},\ldots,\alpha_1)\times\\
(x_1 y_1)^{-\alpha_1}\ldots (x_n y_n)^{-\alpha_{2n-1}}y_n^{-2}.
\end{multline*}
D'o{\`u} {\`a} l'aide du lemme~{\ref{lem3}},
\begin{multline*}
e_{+}.u_k=
\displaystyle\sum\limits_{i+j=k}
\displaystyle\sum\limits_{({\underline \beta},
{\underline \gamma})
\atop{
\arrowvert{\underline \beta}\arrowvert =i-1;
\arrowvert{\underline \gamma}\arrowvert=j-1
}}
q^{\beta_2\gamma_1+\ldots+\beta_{2n-1}\gamma_{2n-2}+\gamma_{2n-1}}\times\\
F_q(\beta_{2n-1},\ldots,\beta_1)
F_q(\gamma_{2n-1},\ldots,\gamma_1)\times\\
(x_1 y_1)^{-(\beta_1+\gamma_1)}\ldots 
(x_n y_n)^{-(\beta_{2n-1}+\gamma_{2n-1})}y_n^{-2}.
\end{multline*}
Or, d'apr{\`e}s les relations de commutation
entre les $x_i$ et $y_j$, on v{\'e}rifie sans probl{\`e}me que
\begin{multline*}
(x_1 y_1)^{-\beta_1}\ldots (x_n y_n)^{-\beta_{2n-1}}y_n^{-1}
(x_1 y_1)^{-\gamma_1}\ldots (x_n y_n)^{-\gamma_{2n-1}}y_n^{-1}\\
=q^{\beta_2\gamma_1+\ldots+\beta_{2n-1}\gamma_{2n-2}+\gamma_{2n-1}}\times\\
(x_1 y_1)^{-(\beta_1+\gamma_1)}\ldots 
(x_n y_n)^{-(\beta_{2n-1}+\gamma_{2n-1})}y_n^{-2}
\end{multline*}
Donc, 
\begin{multline*}
e_{+}.u_k=
\displaystyle\sum\limits_{i+j=k}
\Bigl [
\displaystyle\sum\limits_{{\underline \beta}
\atop{\arrowvert{\underline \beta}\arrowvert =i-1}}
F_q(\beta_{2n-1},\ldots,\beta_1)
(x_1 y_1)^{-\beta_1}\ldots (x_n y_n)^{-\beta_{2n-1}}y_n^{-1}
\Bigr ]\\
\times\Bigl [\displaystyle\sum\limits_{{\underline \gamma}
\atop{\arrowvert{\underline \gamma}\arrowvert =j-1}}
F_q(\gamma_{2n-1},\ldots,\gamma_1)
(x_1 y_1)^{-\gamma_1}\ldots (x_n y_n)^{-\gamma_{2n-1}}y_n^{-1}\Bigr ]
\end{multline*}
Soit, finalement~:
$$
e_{+}.u_k=\displaystyle\sum\limits_{i+j=k}u_i^{(n)}u_j^{(n)}.
$$

\end{dem}

\section{Le morphisme d'alg{\`e}bres de
${\cal A}_q^{(n)}$ dans $A_q^{(n)}$}\label{E4}

\noindent
Bien que cette section soit courte, nous voulons mettre en {\'e}vidence
le morphisme de ${\cal A}_q^{(n)}$ dans $A_q^{(n)}$
dont la restriction donne le plongement de
${\cal U}_q^{(n)}$ dans $A_q^{(n)}$.
\smallskip

Vu l'isomorphisme (\ref{isonp}) entre $U_q(\nm)$ et ${\cal A}_q$,
de mani{\`e}re groupique, les morphismes $\varphi_i$
et $\psi_i$ de \ref{action} sont les suivants~:
$$
\begin{array}{rcl}
\varphi_i:\quad {\cal A}_q&\longrightarrow&A_q^{(n)}\\
{\cal L}(\lambda)&\longmapsto&\begin{pmatrix}
1&0\\
\lambda x_i&1
\end{pmatrix},
\end{array}
$$
et~:
$$
\begin{array}{rcl}
\psi_i:\quad {\cal A}_q&\longrightarrow&A_q^{(n)}\\
{\cal L}(\lambda)&\longmapsto&\begin{pmatrix}
1&y_i\\
0&1
\end{pmatrix}.
\end{array}
$$
Donc, par composition, et en utilisant (\ref{debiso}),
le morphisme $\varphi$ de (\ref{ecran}) n'est autre que~:
\begin{equation}\label{lemorq}
\begin{array}{rcl}
\varphi:\quad {\cal A}_q&\longrightarrow&A_q^{(n)}\\
{\cal L}(\lambda)&\longmapsto&
\displaystyle\prod\limits_{i=1}^{n}
\begin{pmatrix}
1&0\\
\lambda x_i&1
\end{pmatrix}
\begin{pmatrix}
1&y_i\\
0&1
\end{pmatrix}.
\end{array}
\end{equation}
Il appara{\^\i}t que $\varphi$ est l'analogue quantique
de l'application $\varphi_{\textnormal{cl}}$ classique
de \ref{phicl}, rendant naturel son introduction.

Il est clair que pour $k>n,\, a_{1,1}^{(k-1)},\,
a_{1,2}^{(k-1)},\, a_{2,1}^{(k)},\,  a_{2,2}^{(k)}$ appartiennent {\`a} 
$\ker\varphi$.
De plus,
  $$
  \begin{array}{rcl}
  \varphi\bigl( a_{1,1}^{(n-1)}\bigr)&=&
  y_1 x_2\ldots y_{n-1} x_n,\\
  \varphi\bigl( a_{1,2}^{(n-1)}\bigr)&=&
  y_1 x_2\ldots y_{n-1} x_n y_n,\\
  \varphi\bigl( a_{2,1}^{(n)}\bigr)&=&
  x_1 y_1\ldots x_{n-1} y_{n-1} x_n,\\
  \varphi\bigl( a_{2,2}^{(n)}\bigr)&=&
  x_1 y_1\ldots x_{n-1} y_{n-1} x_n y_n,
  \end{array}
  $$
et chacun des termes de droite est inversible dans $A_q^{(n)}$.
Par suite, d'apr{\`e}s \ref{defaqn}, on en d{\'e}duit l'existence d'un
morphisme que l'on note encore $\varphi$ par abus de notation~:
$$
\varphi:\quad {\cal A}_q^{(n)}\longrightarrow A_q^{(n)}
$$
Contrairement au cas classique, on ne sait pas si cette application
est injective.
\section{Relations de commutations entre les $u_i$}\label{E5}

\noindent
Dans le cas quantique, par restriction du morphisme $\varphi$ pr{\'e}c{\'e}dent
{\`a} ${\cal U}_q^{(n)}$, on a un diagramme commutatif {\'e}quivalent {\`a}
(\ref{diagclfonc}), sauf que l'on ne sait pas que les fl{\`e}ches 
en diagonales sont
injectives (mais nous montrerons que celle de gauche l'est
effectivement)~:
\begin{equation}\label{diagqfonc}
\begin{array}{rcl}
{\cal U}_q^{(n)}&\longrightarrow&\quad {\cal A}_q^{(n)}\\
\searrow&&\swarrow\\
&A_q^{(n)}&
\end{array}
\end{equation}
Soit $\lambda$ une ind{\'e}termin{\'e}e, et 
$v(\lambda^{-1})=a_{2,2}(\lambda)^{-1}a_{2,1}(\lambda)$.
Par analogie avec le cas classique,
on posera~:
\begin{equation}\label{somuiq}
v(\lambda^{-1})=\displaystyle\sum\limits_{i=0}^{+\infty}
(-1)^{i}\a_{i+1}\lambda^{-i}
\end{equation}
avec $\forall i\in{\NN}^{*},\, \a_{i}\in{\cal A}_q^{(n)}$
et $\forall i\in\lbrace 1,\ldots,2n\rbrace,\, 
\a_{i}\in{\cal U}_q^{(n)}$.
L'alg{\`e}bre ${\cal U}_q^{(n)}$ est par d{\'e}finition
l'alg{\`e}bre engendr{\'e}e par les $\a_{i},\, 1\leq i\leq 2n$
(voir \ref{defuqng} et la fin de \ref{phicl} pour le cas classique).

Nous allons {\'e}tablir une relation de commutation entre 
$v(\lambda^{-1})$ et $v(\mu^{-1})$ qui n'est autre que l'analogue
quantique de (\ref{vadon}), pour $\lambda$
et $\mu$ deux ind{\'e}termin{\'e}es.
Celle-ci nous fournira les relations de commutations
voulues entre les $u_i$ pour $1\leq i\leq 2n$,
en utilisant le morphisme $\varphi$.

\begin{prop}\label{vadonq}
L'{\'e}galit{\'e} suivante est v{\'e}rifi{\'e}e dans
${\cal A}_q^{(n)}((\lambda^{-1},\mu^{-1}))$~:
$$
q(\lambda -\mu)\bigl[
v(\lambda^{-1}),v(\mu^{-1})\bigr]
=(q-1)\bigl(
\mu v(\lambda^{-1})-\lambda v(\mu^{-1})\bigr)
\bigl( v(\lambda^{-1})-v(\mu^{-1})\bigr).
$$
\end{prop}

\begin{dem}
On a~:
\begin{align*}
\bigl[ 
v(\lambda^{-1})- v(\mu^{-1})
\bigr]&=
\bigl[ 
{a_{2,2}(\lambda)}^{-1}
a_{2,1}(\lambda),{a_{2,2}(\mu)}^{-1}
a_{2,1}(\mu)
\bigr]\\
&=
\bigl[ 
{a_{2,2}(\lambda)}^{-1}
a_{2,1}(\lambda),{a_{2,2}(\mu)}^{-1}
\bigr]\,
a_{2,1}(\mu)\\
&\phantom{=}+
{a_{2,2}(\mu)}^{-1}
\bigl[ 
{a_{2,2}(\lambda)}^{-1}
a_{2,1}(\lambda),a_{2,1}(\mu)
\bigr]\\
&={a_{2,2}(\lambda)}^{-1}
\bigl[ 
a_{2,1}(\lambda),{a_{2,2}(\mu)}^{-1}
\bigr]
a_{2,1}(\mu)\\
&\phantom{=}+
\bigl[
{a_{2,2}(\lambda)}^{-1},{a_{2,2}(\mu)}^{-1}
\bigr]\,
a_{2,1}(\lambda)a_{2,1}(\mu)\\
&\phantom{=}+{a_{2,2}(\mu)}^{-1}{a_{2,2}(\lambda)}^{-1}
\bigl[
a_{2,1}(\lambda),a_{2,1}(\mu)
\bigr]\\
&\phantom{=}+{a_{2,2}(\mu)}^{-1}
\bigl[
{a_{2,2}(\lambda)}^{-1},a_{2,1}(\mu)
\bigr]\,
a_{2,1}(\lambda)
\end{align*}
Mais, 
$\bigl[ 
a_{2,1}(\lambda),a_{2,1}(\mu)
\bigr]=\bigl[
{a_{2,2}(\lambda)}^{-1},a_{2,2}(\mu^{-1})
\bigr]=0$ d'apr{\`e}s (\ref{RR1}).
Donc,
\begin{align*}
\bigl[ 
v(\lambda^{-1})- v(\mu^{-1})
\bigr]&={a_{2,2}(\lambda)}^{-1}\bigl[
a_{2,1}(\lambda),{a_{2,2}(\mu)}^{-1}\bigr]\,
a_{2,1}(\mu)\\
&\phantom{=}+{a_{2,2}(\mu)}^{-1}\bigl[
{a_{2,2}(\lambda)}^{-1},a_{2,1}(\mu)\bigr]\,
a_{2,1}(\lambda)\\
&=-{a_{2,2}(\lambda)}^{-1}{a_{2,2}(\mu)}^{-1}
\bigl[
a_{2,1}(\lambda),a_{2,2}(\mu)
\bigr]\,
{a_{2,2}(\mu)}^{-1}a_{2,1}(\mu)\\
&\phantom{=}
+{a_{2,2}(\lambda)}^{-1}{a_{2,2}(\mu)}^{-1}
\bigl[
a_{2,1}(\mu),a_{2,2}(\lambda)
\bigr]\,
{a_{2,2}(\lambda)}^{-1}a_{2,1}(\lambda)
\end{align*}
Or, on peut montrer que (\ref{RR1})
conduit aux relations~:
\begin{align*}
[ a_{2,1}(\lambda),a_{2,2}(\mu)]&=
[ a_{2,1}(\mu),a_{2,2}(\lambda)]\\
q(\lambda -\mu)[ a_{2,1}(\lambda),a_{2,2}(\mu)]&=
(q-1)\bigl(
\mu a_{2,2}(\mu)a_{2,1}(\lambda)
-\lambda a_{2,2}(\lambda)a_{2,1}(\mu)\bigr)
\end{align*}
Par suite, on obtient~:
\begin{multline*}
\bigl[ 
v(\lambda^{-1})- v(\mu^{-1})
\bigr]
={a_{2,2}(\lambda)}^{-1}{a_{2,2}(\mu)}^{-1}
\bigl[
a_{2,1}(\lambda),a_{2,2}(\mu)
\bigr]\\
\times\bigl(
{a_{2,2}(\lambda)}^{-1}a_{2,1}(\lambda)-
{a_{2,2}(\mu)}^{-1}a_{2,1}(\mu)
\bigr)
\end{multline*}
En utilisant la relation~:
$$
a_{2,2}(\lambda)a_{2,2}(\mu)=a_{2,2}(\mu)a_{2,2}(\lambda),
$$
qui provient de (\ref{RR1}),
on obtient le r{\'e}sultat.
\end{dem}

De la proposition pr{\'e}c{\'e}dente, nous pouvons en d{\'e}duire les lois
de commutation dans ${\cal U}_{q}^{(n)}$.
Le corollaire suivant est l'analogue quantique de (\ref{ahoha})~:
\begin{cor}\label{lecora}
Soient $1\leq i<j\leq 2n$ deux entiers. On a~:
$$
q\, [\a_{i},\a_{j}]=(q-1)\displaystyle\sum\limits_{k=i}^{j-1}
\a_{k}\a_{i+j-k}.
$$
\end{cor}
\begin{dem}
Soient $k$ et $l$ deux entiers positifs.
En prenant
le coefficient de $\lambda^{k-1}\mu^{l-1}$ dans chacun des membres de
l'{\'e}galit{\'e} exprim{\'e}e par la proposition \ref{vadonq}, on obtient~:
\begin{equation}\label{eblea}
q\bigl(
[\a_{k+1},\a_{l}]-
[\a_{k},\a_{l+1}]
\bigr)
=(q-1)( -\a_{k}\a_{l+1}
-\a_{l}\a_{k+1}).
\end{equation}
En faisant $k=l$, on voit que~:
\begin{equation}\label{derui}
\forall i\in{\NN}^*,\,
q[\a_{i},\a_{i+1}]=(q-1)
\a_{i}\a_{i+1}.
\end{equation}
ce qui est le r{\'e}sultat demand{\'e} avec
$j-i=1$.
La formule g{\'e}n{\'e}rale se prouve alors par r{\'e}currence sur
$j-i$ {\`a} l'aide de (\ref{eblea}).
\end{dem}

Relions maintenant les $\a_{i}$ aux $u_{i}$.
Le lemme suivant g{\'e}n{\'e}ralise l'{\'e}galit{\'e} classique (\ref{somuicl})~:
\begin{lem}\label{lienui}
Pour tout $i\in{\NN}^{*},\, \varphi(\a_{i})=u_{i}$.
Autrement dit, 
$$
\varphi\bigl( v(\lambda^{-1})\bigr)
=\displaystyle\sum\limits_{k=0}^{+\infty}
(-1)^{k}u_{k+1}\lambda^{-k}.
$$
\end{lem}
\begin{dem}
Posons~:
$$
L_{n}=\begin{pmatrix}
a_n&b_n\\
c_n&d_n
\end{pmatrix}:=\displaystyle\prod\limits_{i=1}^{n}
\begin{pmatrix}
1&0\\
\lambda x_i&1
\end{pmatrix}
\begin{pmatrix}
1&y_i\\
0&1
\end{pmatrix},
$$
et $U_n:=d_{n}^{-1}c_{n}$.
\smallskip

\noindent
$\bullet$ Dans le cas o{\`u} $n=1$, on v{\'e}rifie sans probl{\`e}me que
$d_1=\lambda (x_1 y_1)$ est inversible, et~: 
$$U_1=\bigl ( 1+(\lambda x_1 y_1)^{-1}\bigr )^{-1}y_1^{-1}.$$
Le r{\'e}sultat en d{\'e}coule.

\noindent
$\bullet$ Supposons $n>1$.
On a~: 
$$ L_n=L_{n-1}
\begin{pmatrix}
1&y_n\\
\lambda x_n& \lambda x_n y_n +1 
\end{pmatrix}.$$
Donc, 
$$\left\{ \begin{array}{l}
c_n=c_{n-1}+\lambda d_{n-1} x_n\\
d_n=c_{n-1}y_n +d_{n-1}(\lambda x_n y_n +1)
\end{array}
\right.$$
Donc, $d_n\in{\mathbb C}[x_i^{\pm 1},y_i^{\pm 1}]_q
[\lambda]$ et son coefficient dominant $\omega_n$ v{\'e}rifie~:
\begin{displaymath}
\omega_n=\omega_{n-1}(x_n y_n).
\end{displaymath}
Une r{\'e}currence imm{\'e}diate nous prouve alors que 
$\omega_n\in{\mathbb C}[x_i^{\pm 1},y_i^{\pm 1}]_q$ est inversible.
Par suite, $d_n$ est inversible.
De plus,
\begin{align*}
U_n&=\bigl (c_{n-1}y_n+d_{n-1}(\lambda x_n y_n +1)
\bigr )^{-1}(c_{n-1}+\lambda d_{n-1} x_n)\\
&=\Bigl ( (c_{n-1} +\lambda d_{n-1} x_n )^{-1} 
\bigl ( c_{n-1}y_n +d_{n-1} (\lambda x_n y_n +1)\bigr )\Bigr )^{-1}\\
&=\Bigl ( (c_{n-1} +\lambda d_{n-1} x_n )^{-1} 
(d_{n-1}^{-1})^{-1}d_{n-1}^{-1}
\bigl ( c_{n-1}y_n +d_{n-1} (\lambda x_n y_n +1)\bigr )\Bigr )^{-1}\\
&=\Bigl ( (d_{n-1}^{-1}c_{n-1} +\lambda x_n )^{-1} 
\bigl (d_{n-1}^{-1} c_{n-1}y_n + \lambda x_n y_n +1\bigr )\Bigr )^{-1}\\
&=\Bigl ( y_n \bigl ( d_{n-1}^{-1}c_{n-1}y_n +\lambda x_n y_n \bigr )^{-1}
\bigl ( d_{n-1}^{-1} c_{n-1} y_n +\lambda x_n y_n +1\bigr )\Bigr)^{-1}\\
&=\Bigl ( \bigl ( d_{n-1}^{-1}c_{n-1}y_n +\lambda x_n y_n \bigr )^{-1}
\bigl ( d_{n-1}^{-1} c_{n-1} y_n +\lambda x_n y_n +1\bigr )\Bigr)^{-1}
y_n^{-1}\\
&=\Bigl ( 1+\bigl ( d_{n-1}^{-1}c_{n-1}y_n +\lambda x_n y_n \bigr )^{-1}
\Bigr )^{-1} y_n^{-1}\\
&=\Bigl ( 1+\bigl ( d_{n-1}^{-1}c_{n-1}\lambda^{-1}x_n^{-1}
(\lambda x_n y_n) +(\lambda x_n y_n)\bigr )^{-1}\Bigr )^{-1}
y_n^{-1}
\end{align*}
\newpage
\begin{align*}
&=\Bigl ( 1+ (\lambda x_n y_n)^{-1}\bigl ( 
1+d_{n-1}^{-1}c_{n-1}\lambda^{-1}x_n^{-1}
\bigr )^{-1}\Bigr )^{-1}y_n^{-1}\\
&=\Bigl ( 1+ (\lambda x_n y_n)^{-1}\bigl ( 
1+d_{n-1}^{-1}c_{n-1}y_{n-1}\lambda^{-1}y_{n-1}^{-1}x_n^{-1}
\bigr )^{-1}\Bigr )^{-1}y_n^{-1}\\
&=\Bigl ( 1+ (\lambda x_n y_n)^{-1}\bigl ( 
1+q^{-1}(d_{n-1}^{-1}c_{n-1}y_{n-1})(\lambda y_{n-1}x_n)^{-1}
\bigr )^{-1}\Bigr )^{-1}y_n^{-1}
\end{align*}
car $$ y_{n-1}^{-1}x_n^{-1}=q^{-1}(y_{n-1} x_n)^{-1}.$$
Donc, en posant $V_n=U_n y_n$, et 
$W_n^{(k)}=\bigl ( 1+q^{k} V_{n-1}(\lambda y_{n-1} x_n)^{-1}
\bigr )^{-1}$,
on voit que
$$V_n=\bigl ( 1+ (\lambda x_n y_n)^{-1}W_n^{(-1)}\bigr )^{-1}.$$
De l{\`a}, on en d{\'e}duit par r{\'e}currence que $\forall k>n,$ $(x_k y_k)$
et $V_n$ commutent.
Puis, 
$$
\forall p\in{\mathbb N},\,
(\lambda x_n y_n)^{-p}W_n^{(k)}=W_n^{(k+p)}(\lambda x_n y_n)^{-p}.
$$
D'o{\`u}~:
\begin{align*} 
V_n&=\bigl ( 1+W_n^{(0)}(\lambda x_n y_n)^{-1}\bigr )^{-1}\\
&=\Bigl (1+\bigl ( 1+V_{n-1}(\lambda y_{n-1} x_n)^{-1}\bigr )^{-1}
(\lambda x_n y_n)^{-1}\Bigr )^{-1}.
\end{align*}
Soit, {\`a} la suite d'une nouvelle r{\'e}currence,
\begin{multline*}
V_n=\Biggl ( 1+\biggl ( 1+\Bigl ( 1+\ldots\bigl (
1+(\lambda x_1 y_1)^{-1}\bigr )^{-1}(\lambda y_1 x_2)^{-1}
\Bigr )^{-1}\ldots\\
\ldots (\lambda y_{n-1} x_n)^{-1}\biggr )^{-1}
(\lambda x_n y_n)^{-1}\Biggr )^{-1}
\end{multline*}
Le r{\'e}sultat d{\'e}coule alors de T1 et de (\ref{lemorq}).
\end{dem}
\begin{cor}
Soient $1\leq i<j\leq 2n$ deux entiers. On a~:
$$
q\, [u_{i},u_{j}]=(q-1)\displaystyle\sum\limits_{k=i}^{j-1}
u_{k}u_{i+j-k}.
$$
\end{cor}\label{lecorq}
\begin{dem}
Il suffit d'appliquer le corollaire \ref{lecora},
le lemme \ref{lienui} et le morphisme d'alg{\`e}bres $\varphi$.

\end{dem}

\section{D{\'e}monstration de T2 et de T3}

\noindent
D'ores et d{\'e}ja, d'apr{\`e}s la proposition \ref{calculs+}
et le corollaire \ref{lecorq}, nous savons que
$U_q^{(n)}$ est une sous $U_q(\bm)$-alg{\`e}bre-module
de $A_q^{(n)}$.
Il faut voir que les relations du corollaire \ref{lecorq}
sont les seules entre les $u_i,\, 1\leq i\leq 2n$.

Soit $N\in{\NN}^{*}$, et consid{\'e}rons 
${\cal M}_q^{(N)}$ l'alg{\`e}bre engendr{\'e}e sur
$\CC[q,q^{-1}]$
par g{\'e}n{\'e}rateurs: $t_1,\ldots,t_{N}$
et relations~:
$$
\forall i<j,\quad q\, [t_i,t_j]=(q-1)
\displaystyle\sum\limits_{k=i}^{j-1}
t_k t_{i+j-k}.
$$
Alors, d'apr{\`e}s le corollaire \ref{lecorq},
il existe un morphisme d'alg{\`e}bres not{\'e} $\phi_{2n}$~:
$$
\begin{array}{rcl}
\phi_{2n}:\quad {\cal M}_q^{(2n)}&\longrightarrow&
A_q^{(n)}\\
t_i&\longmapsto&u_i
\end{array}
$$
Consid{\'e}rons {\'e}galement la simple application lin{\'e}aire entre
espace vectoriels~:
$$
\begin{array}{rcl}
\psi_{N}:\quad\CC[X_1,\ldots,X_N]
\quad &\longrightarrow&{\cal M}_q^{(N)}\\
\sum\limits_{{\underline{\a}}}
a_{{\underline{\a}}} X_1^{\a_1}\ldots X_N^{\a_N}
&\longmapsto&
\sum\limits_{{\underline{\a}}}a_{{\underline{\a}}}
\textnormal{cl}\bigl(
t_1^{\a_1}\ldots t_N^{\a_N}\bigr)
\end{array}
$$
Posons $i_{2n}=\phi_{2n}\circ \psi_{2n}$.
Le diagramme suivant est {\'e}videmment commutatif~:
\begin{equation}\label{bof1}
\begin{array}{rcll}
\CC[X_1,\ldots,X_{2n}]&\overset{\psi_{2n}}{\longrightarrow}&
{\cal M}_q^{(2n)}\\
i_{2n}\searrow&&\swarrow\phi_{2n}\\
&A_q^{(n)}&
\end{array}
\end{equation}
\begin{lem}\label{psin}
Pour tout $N$, l'application $\psi_N$ est surjective.
\end{lem}
\begin{dem}
Par r{\'e}currence sur $j-i$,
en utilisant plusieurs fois la formule~:
$$
\forall k<l,\quad
q\, [t_k,t_l]=
(q-1)\sum_{\a=k}^{l-1}t_{\a}t_{k+l-\a},
$$
on montre que si $i<j$, alors $t_jt_i$ est
une combinaison lin{\'e}aire
de termes de la forme $t_u t_v$ avec $i\leq u\leq v\leq j$.
D'o{\`u} le r{\'e}sultat.
En particulier, remarquons que l'on a l'{\'e}galit{\'e}~:
\begin{equation}\label{titii}
\forall i,\quad t_i t_{i+1}=qt_{i+1}t_i
\end{equation}
\end{dem}
\begin{prop}\label{in}
L'application $i_{2n}$ est injective.
\end{prop}
\begin{dem}
Pour $k\in\lbrace 1,\ldots,2n\rbrace$,
notons
$\CC[x_k^{\pm 1},\ldots,y_n^{\pm 1}]_q$
la sous-alg{\`e}bre de $A_q^{(n)}$ engendr{\'e}e par
$x_k^{\pm 1},\ldots,y_n^{\pm 1}$. De m{\^e}me,
$\CC[y_k^{\pm 1},\ldots,y_n^{\pm 1}]_q$
d{\'e}signe la sous-alg{\`e}bre de $A_q^{(n)}$ engendr{\'e}e par
$y_k^{\pm 1},\ldots,y_n^{\pm 1}$.
Alors, d'apr{\`e}s \ref{modlib2}, la sous-alg{\`e}bre
$\CC[x_k^{\pm 1},\ldots,y_n^{\pm 1}]_q$
est un $\CC[y_k^{\pm 1},\ldots,y_n^{\pm 1}]_q$-module libre
({\`a} gauche ou {\`a} droite) dont une base est donn{\'e}e
par la famille ${(x_{k}^{p})}_{p\in\NN}$.
De m{\^e}me, la sous-alg{\`e}bre $\CC[y_k^{\pm 1},\ldots,y_n^{\pm 1}]_q$
est un $\CC[x_{k+1}^{\pm 1},\ldots,y_n^{\pm 1}]_q$-module
libre ({\`a} gauche ou {\`a} droite) dont une base est donn{\'e}e
par la famille ${(y_{k}^{p})}_{p\in\NN}$.

Toujours pour $k\in\lbrace 1,\ldots,2n\rbrace$,
notons $j_k$ la restriction de l'application 
$i_{2n}$ {\`a} 
$\CC[X_1,\ldots,X_k]$, et montrons par r{\'e}currence sur $k$
que $j_k$ est injective.

\noindent
$\bullet$ Si $k=1$, c'est clair car $i_1(X_1)=u_1=y_n^{-1}$
et la famille $(y_n^{-k})_{k\in\NN}$ est libre dans 
$A_{q}^{(n)}$.

\smallskip
\noindent
$\bullet$ Soit $k\geq 2$, et supposons le r{\'e}sultat vrai jusqu'au
rang $k-1$. Pour fixer les id{\'e}es, supposons $k$ pair, $k=2r$.
Prenons~: 
$$
x\in\CC[X_1,\ldots,X_k]\setminus \CC[X_1,\ldots,X_{k-1}],
$$
et montrons que $j_k(x)\not= 0$. Posons $\a =\partial_{k}^{o}x$.
Il existe des polyn{\^o}mes $P_0,\ldots,P_{\a}$ dans
$\CC[X_1,\ldots,X_{k-1}]$ tels que
$P_{\a}(X_1,\ldots,X_{k-1})\not=0$ et~:
\begin{equation}\label{hihan2}
x=P_0(X_1,\ldots,X_{k-1})+\ldots+
P_{\a}(X_1,\ldots,X_{k-1})X_k^{\a}.
\end{equation}
Par ailleurs, en examinant la forme des $u_i$,
on constate que~:
$$
\forall i\in\lbrace 1,\ldots, k-1\rbrace,\quad
u_i\in\CC[y_{n-r+1}^{-1},\ldots,x_n^{-1},y_n^{-1}]_q,
$$
et,
\begin{equation}\label{lala+}
u_k=v+(x_{n-r+1}y_{n-r+1})^{-1}\ldots (x_n y_n)^{-1}y_n^{-1}
\end{equation}
avec $v\in\CC[y_{n-r+1}^{-1},\ldots,x_n^{-1},y_n^{-1}]_q$.
Donc, par (\ref{hihan2}), 
$$
j_k(x)=i_{2n}(x)\in\CC[x_{n-r+1}^{-1},\ldots,x_n^{-1},y_n^{-1}]_q,
$$
et le coefficient de $j_k(x)$ sur $x_{n-r+1}^{-\a}$ 
dans $\CC[x_{n-r+1}^{\pm 1},\ldots,y_n^{\pm 1}]_q$
consid{\'e}r{\'e} comme
$\CC[y_{n-r+1}^{\pm 1},\ldots,y_n^{\pm 1}]_q$-module
est
$j_{k-1}\bigl( P_{\a}(X_1,\ldots,X_{k-1})\bigr)$ qui est non nul
d'apr{\`e}s l'hypoth{\`e}se de r{\'e}currence.
Donc, $j_k(x)\not= 0$.
Le m{\^e}me raisonnement s'adaptant au cas $k$ impair,
on en d{\'e}duit que $j_{k}$ est injective.
Puisque $j_{2n}=i_{2n}$, la proposition est d{\'e}montr{\'e}e.
\end{dem}
\begin{cor}\label{bof3}
L'application $\phi_{2n}$ est injective.
\end{cor}
\begin{dem}
La surjectivit{\'e} de $\psi_{2n}$,
l'injectivit{\'e} de $i_{2n}$ et la commutativit{\'e} du
diagramme (\ref{bof1}) entra{\^\i}ne imm{\'e}diatement
l'injectivit{\'e} de $\phi_{2n}$.
\end{dem}
\begin{cor}\label{bof2}
Pour tout $N$, l'application $\psi_{N}$ est un isomorphisme de 
$\CC$-espace vectoriel.
\end{cor}
\begin{dem}
Le m{\^e}me raisonnement que pr{\'e}c{\'e}demment montre que
$\psi_{2n}$ est bijectif.

Reste {\`a} montrer que
$\psi_{2n-1}$ est bijectif pour $n\geq 1$.
Soit $n\geq 1$. On a un morphisme d'alg{\`e}bres~:
$$
\begin{array}{rcl}
r:\quad {\cal M}_q^{(2n-1)}&\longrightarrow&{\cal M}_q^{(2n)}\\
t_i&\longmapsto&t_i
\end{array}
$$
ce qui nous permet de d{\'e}finir le morphisme d'alg{\`e}bres~:
$$
\phi_{2n-1}:\quad {\cal M}_q^{(2n-1)}\longrightarrow
\CC[x_1^{\pm 1},\ldots,y_n^{\pm 1}]_q
$$
par $\phi_{2n-1}:=\phi_{2n}\circ r$.
Alors, clairement, $r\circ\psi_{2n-1}$ est la restriction de
$\psi_{2n}$ {\`a} $\CC[X_1,\ldots,X_{2n-1}]$. De plus,
en notant comme dans la d{\'e}monstration pr{\'e}c{\'e}dente
$j_{2n-1}$ la restriction de $i_{2n}$ {\`a}
$\CC[X_1,\ldots,X_{2n-1}]$, on a le diagramme commutatif suivant~:
\begin{equation}\label{bof4}
\begin{array}{rcccl}
\CC[X_1,\ldots,X_{2n-1}]&\overset{\psi_{2n-1}}{\longrightarrow}&
{\cal M}_q^{(2n-1)}&\overset{r}{\longrightarrow}&
{\cal M}_q^{(2n)}\\
j_{2n-1}&\searrow&\downarrow\phi_{2n-1}&\swarrow&\phi_{2n}\\
&&\CC[x_1^{\pm 1},\ldots,y_n^{\pm 1}]_q&&
\end{array}
\end{equation}
Par suite, l'injectivit{\'e} de $j_{2n-1}$ entra{\^\i}ne celle
de $\psi_{2n-1}$. D'o{\`u} le r{\'e}sultat.
\end{dem}

Etant donn{\'e} que l'image de ${\cal M}_q^{(2n)}$
par $\phi_{2n}$ n'est autre que $U_{q}^{(n)}$,
nous avons presque d{\'e}montr{\'e} T2.
Il nous reste {\`a} voir que $U_{q}^{(n)}$
a le m{\^e}me corps de fractions que $A_{q}^{(n)}$.
Ceci d{\'e}coule du lemme suivant~:
\begin{lem}
Pour tout $k\in\lbrace 1,\ldots,2n\rbrace$,
le corps de fractions de l'alg{\`e}bre engendr{\'e}e
par $u_1,\ldots,u_k$ est le m{\^e}me que celui de
$\CC[x_{n-r+1}^{\pm 1},\ldots,y_n^{\pm 1}]_q$
si $k=2r$ et $\CC[y_{n-r}^{\pm 1},\ldots,y_n^{\pm 1}]_q$
si $k=2r+1$.
\end{lem}
\begin{dem}
Cela d{\'e}coule de la formule (\ref{lala+}) si $k=2r$ et d'une relation
analogue si $k=2r+1$.
\end{dem}

Par suite, T2 est d{\'e}montr{\'e}. De plus,
d'apr{\`e}s le corollaire \ref{bof2}, on voit qu'une base de
$U_{q}^{(n)}$ est donn{\'e}e par la famille
$\displaystyle\prod\limits_{i=1}^{2n}u_{i}^{\a_i},\, \a_i\in\NN$.

Il nous reste {\`a} v{\'e}rifier T3.
\begin{prop}
La restriction de $\varphi$ {\`a} ${\cal U}_{q}^{(n)}$
est un isomorphisme de ${\cal U}_{q}^{(n)}$
sur $U_{q}^{(n)}$.
\end{prop}
\begin{dem}
D'apr{\`e}s le corollaire \ref{lecora}
et le fait que $\phi_{2n}$ soit un isomorphisme de 
${\cal M}_q^{(2n)}$ sur $U_q^{(n)}$,
on a un morphisme d'alg{\`e}bres~:
$$
\begin{array}{rcl}
\psi:\quad U_q^{(n)}&\longrightarrow&{\cal U}_q^{(n)}\\
u_i&\longmapsto&\a_i
\end{array}
$$
Ce morphisme $\psi$ n'est autre que l'inverse recherch{\'e}.
\end{dem}

Ainsi, une base de ${\cal U}_q^{(n)}$
est donn{\'e}e par la famille $\displaystyle\prod\limits_{i=1}^{2n}
\a_i^{a_i}$, avec $a_i\in\NN$.
Le corollaire \ref{lecora}
et la formule (\ref{ahoha})
donnant le crochet de Poisson sur ${\cal U}^{(n)}$
ach{\`e}ve la preuve de T3.

\section{Fin de la d{\'e}monstration}

\noindent
Dor{\'e}navant, $q$ d{\'e}signera un nombre complexe non nul,
non racine de l'unit{\'e}, et $U_q(\bm)$ d{\'e}signera l'alg{\`e}bre
obtenue {\`a} partir de l'``ancienne" alg{\`e}bre  de Hopf $U_q(\bm)$
en adjoignant une racine carr{\'e}e de $k:\, k^{{1\over 2}}$,
ainsi que son inverse $k^{-{1\over 2}}$.

Alors, $U_q(\bm)$ reste une alg{\`e}bre de Hopf en posant~:
$$
k^{{1\over 2}} e_{\pm} k^{-{1\over 2}}=q^{\pm {1\over 2}}e_{\pm},
$$
et $A_{q}^{(n)}$ reste une $U_q(\bm)$-alg{\`e}bre de Hopf en posant~:
$$
k^{\pm {1\over 2}}.P=q^{\pm {1\over 2}\deg{P}}P,
$$
pour $P$ homog{\`e}ne.

Par ailleurs, pour $(i,j)\in {\lbrace 1,\ldots,2n\rbrace}^{2}$,
avec $i<j$, nous noterons par 
$\CC[u_i,\ldots,u_j]_q$ la sous-alg{\`e}bre de $A_{q}^{(n)}$
engendr{\'e}e par $u_k$ pour $i\leq k\leq j$.
C'est une $U_q(\bm)$-sous-alg{\`e}bre-module de $U_q^{(n)}$.

\subsection{Graduation principale sur $A_q^{(n)}$ et $U_q(\bm)$.}

\noindent
En plus de la graduation utilis{\'e}e jusqu'ici, nous disposons d'une
graduation principale sur $A_q^{(n)}$ et $U_q(\bm)$ d{\'e}finie de la 
fa\c con suivante~:
$$
\degp{x_i}=\degp{y_i}=1,
$$
et~:
$$
\degp{e_{\pm}}=1,\quad \degp{k^{\pm {1\over 2}}}=0.
$$
Cette graduation (principale) est int{\'e}ressante car l'action 
de $U_q(\bm)$ sur $A_q^{(n)}$ est gradu{\'e}e (pour cette graduation 
comme pour l'autre).

Il est facile de voir que~: 
$$
\forall i\in\lbrace 1,\ldots,2n\rbrace,
\quad \deg{u_i}=1\textnormal{ et }
\degp{u_i}=1-2i.
$$
Par suite, en posant $\dega{}=\deg{}$,
et $\degb{}=\displaystyle{1\over 2}(\deg{}-\degp{} )$,
on obtient deux graduations sur $U_q^{(n)}$ et $U_q(\bm)$ telles que~:
$$
\forall i\in\lbrace 1,\ldots,2n\rbrace,
\quad \dega{u_i}=1\, \degb{u_i}=i,
$$
et~:
$$
\dega{e_{\pm}}=\pm 1,\; \degb{e_+}=0,\;
\degb{e_-}=-1,\; \dega{k^{\pm {1\over 2}}}=
\degb{k^{\pm {1\over 2}}}=0.
$$
\begin{lem}\label{atlem0}
Soient $k\in\lbrace 1,\ldots, 2n\rbrace$ et 
$P\in\CC[u_1,\ldots,u_k]_q$ homog{\`e}ne.
Alors, 
$$
\dega{P}\leq\degb{P}\leq k\, \dega{P}.
$$
\end{lem}
\begin{dem}
C'est clair, car si $P=u_1^{\a_1}\ldots u_k^{\a_k}$, alors~:
\begin{align*}
\dega{P}&=\a_1+\a_2+\ldots+\a_k,\\
\textnormal{ et }\quad\quad\degb{P}&=\a_1+2\a_2+\ldots+k\a_k.
\end{align*}
\end{dem}
\subsection{L'id{\'e}al $I_{q}^{(n)}$.}
\noindent
Nous cherchons {\`a} g{\'e}n{\'e}raliser l'isomorphisme 
(\ref{coindclass}) qui provient du couplage non-d{\'e}g{\'e}n{\'e}r{\'e}e (\ref{coupla2}).
Nous sommes donc amener {\`a} consid{\'e}rer l'id{\'e}al {\`a} droite
$I_{q}^{(n)}$ de $U_q(\bm)$ d{\'e}fini par~:
$$
I_{q}^{(n)}:=
\bigl\lbrace
x\in U_q(\bm) / \, \varepsilon(x.P)=0\quad 
\forall P\in U_{q}^{(n)}
\bigr\rbrace,
$$
o{\`u} $\varepsilon$ d{\'e}signe l'{\'e}valuation du terme constant dans
$A_{q}^{(n)}$ comme dans $U_{q}^{(n)}$ (c'est la m{\^e}me chose).
Notons que la fonction $\varepsilon$ est un morphisme
d'alg{\`e}bres et que si $P\in U_{q}^{(n)}$ est homog{\`e}ne
pour $\dega{}$ ou $\degb{}$ de degr{\'e} strictement positif,
alors $\varepsilon (P)=0$.

L'id{\'e}al {\`a} droite $I_{q}^{(n)}$ se d{\'e}crit difficilement dans l'alg{\`e}bre 
$U_q(\bm)$.
Par contre, nous allons voir qu'il s'exprime ais{\'e}ment comme id{\'e}al
de $\CC[\BP]_q$ ainsi qu'en termes de nouvelles r{\'e}alisations
de Drinfeld.

\subsection{D{\'e}finition de $\CC[B+]_q$.}
\noindent
C'est l'alg{\`e}bre engendr{\'e}e sur 
$\CC[q^{{1\over 2}},q^{{-1\over 2}}]$
par les g{\'e}n{\'e}rateurs~: $a_{i,j}^{(k)}$
et les relations~:
\begin{equation}\label{RLL}
R({\lambda,\mu}){\cal L}_1(\lambda){\cal L}_2(\mu)=
{\cal L}_2(\mu){\cal L}_1(\lambda)R({\lambda,\mu}),
\end{equation}
o{\`u} $R({\lambda,\mu})$ est la $R$-matrice de (\ref{mat1}),
et~:
\begin{equation}
a_{2,1}^{(0)}=0,
\end{equation}
ainsi que~:
\begin{equation}\label{detq}
a_{1,1}(q\lambda)\bigl[
a_{2,2}(\lambda)-
a_{2,1}(\lambda)a_{1,1}(\lambda)^{-1}a_{1,2}(\lambda)\bigr]=1.
\end{equation}
La derni{\`e}re relation est celle du d{\'e}terminant quantique.
On a en particulier~:
\begin{align}
a_{1,1}^{(0)}a_{2,2}^{(0)}&=a_{2,2}^{(0)}a_{1,1}^{(0)}=1\\
\label{aii}
\forall\, i,j,\a,\b,\quad
a_{i,j}^{(\a)}a_{i,j}^{(\b)}&=a_{i,j}^{(\b)}a_{i,j}^{(\a)}
\end{align}
On peut {\'e}galement d{\'e}finir une antipode qui ne nous sera pas utile
ici.
On montre que l'on obtient une alg{\`e}bre de Hopf.
La comultiplication est donn{\'e}e par~:
$$
\begin{array}{rcl}
\Delta:\quad \CC[B+]_q&\longrightarrow&\CC[B+]_q\\
{\cal L}(\lambda)&\longmapsto&
{\cal L}(\lambda)\otimes {\cal L}(\lambda)
\end{array}
$$
L'alg{\`e}bre $\CC[\BP]_q$ est bigradu{\'e}e naturellement~:

\noindent
- par $\deg{\bigl( a_{i,j}^{(k)}\bigr)}=k$ correspondant 
{\`a} la d{\'e}rivation~:
$$
d\bigl( L(\lambda)\bigr)=\lambda\displaystyle{\partial\over
\partial\lambda}L(\lambda)
$$

\noindent
- par $\dego{\bigl( a_{i,j}^{(k)}\bigr)}=2(j-i)$ 
correspondant {\`a} la d{\'e}rivation~:
$$
d'\bigl( L(\lambda)\bigr) ={\textnormal{ad}}(h)\bigl(
L(\lambda)\bigr).
$$ 

\noindent
Nous poserons $\dega{}=-\displaystyle{1\over 2}\dego{}$
et $\degb{}=-\displaystyle{1\over 2}\dego{}-\deg{}$.
Donc,
\begin{equation}\label{degaijk}
\forall\, i,j,k,\quad 
\dega{\bigl(a_{i,j}^{(k)}\bigr)}=i-j\quad\textnormal{ et }\quad
\degb{\bigl(a_{i,j}^{(k)}\bigr)}=i-j-k.
\end{equation}
Nous admettrons la proposition~:
\begin{prop}\label{isofren}
L'application suivante entre $U_q(\bm)$ et $\CC[\BP]_q$ 
est un isomorphisme d'alg{\`e}bres de Hopf
bigradu{\'e} pour $\dega{}$ et $\degb{}$~:
$$
\begin{array}{rcl}
U_q(\bm)&\longrightarrow&\CC[\BP]_q\\
e_{+}&\longrightarrow&a_{2,1}^{(1)}a_{2,2}^{(0)}\\
e_{-}&\longrightarrow&a_{1,2}^{(0)}a_{1,1}^{(0)}\\
k^{{1\over 2}}&\longrightarrow&a_{2,2}^{(0)}\\
k^{-{1\over 2}}&\longrightarrow&a_{1,1}^{(0)}
\end{array}
$$
\end{prop}
Nous allons maintenant nous int{\'e}resser {\`a} l'action de
$\CC[\BP]_q$ sur $U_q^{(n)}$ via l'isomorphisme de la proposition 
pr{\'e}c{\'e}dente. 
Commen\c cons par calculer l'action des $a_{i,j}^{(k)}$
sur $U_{q}^{(n)}$.
\subsection{Action de $a_{2,1}^{(k)}$ sur $U_{q}^{(n)}$.}
\noindent
Soit $k\geq 1$. On a 
$\dega{\bigl(a_{2,1}^{(k)}\bigr)}=1$.
Donc, si $P\in A_{q}^{(n)}$ est homog{\`e}ne pour $\dega{}$,
$a_{2,1}^{(k)}.P$ est homog{\`e}ne de degr{\'e} 
strictement positif pour $\dega{}$.
Par suite, $a_{2,1}^{(k)}\in I_{q}^{(n)}$.

\subsection{ Action de $a_{1,1}^{(1)}$ et 
$a_{2,2}^{(1)}$ sur les $u_i,\, 1\leq i\leq 2n$.}
\noindent
On a $\dega{\bigl(a_{1,1}^{(1)}\bigr)}=0$ et
$\degb{\bigl(a_{1,1}^{(1)}\bigr)}=-1$.
Donc $,a_{1,1}^{(1)}\in I_{q}^{(n)}$ et~: 
$$
\forall i,\,
\exists\, \a_i,\quad
a_{1,1}^{(1)}.u_i=\a_i u_{i-1}.
$$
Notons $\phi$ l'isomorphisme de la proposition \ref{isofren}.
Il est bigradu{\'e}. Donc,
$$
\exists\, \lambda,\mu, r,s,\quad
a_{1,1}^{(1)}=
\phi\bigl[
\lambda k^{r}e_{+}e_{-}+\mu k^{s}e_{-}e_{+}
\bigr] .
$$
Par suite, $\a_i$ est ind{\'e}pendant de $i$.
Notons que (\ref{derui}) peut se r{\'e}ecrire~:
\begin{equation}\label{derui2}
\forall k,\quad u_k u_{k+1}=q u_{k+1} u_k.
\end{equation}
En appliquant $a_{1,1}^{(1)}$ aux deux membres de l'{\'e}galit{\'e}
pr{\'e}c{\'e}dente pour $k=1$, et en utilisant la structure de
$\CC[\BP]_q$ alg{\`e}bre-module, on constate que $\a_i=\a\not= 0$.

\noindent
De la m{\^e}me fa\c con, on prouve que $a_{2,2}^{(1)}\in I_{q}^{(n)}$ 
et que~:
$$
\exists\,\b\not= 0\,\,\forall i,\quad 
a_{2,2}^{(1)}.u_i=\b u_{i-1}.
$$
Le calcul montre que 
$\a=\displaystyle{q^{-{1\over 2}}\over q-1}$
et $\b=-\displaystyle{q^{-{1\over 2}}\over q-1}$.
\subsection{ Action de $a_{1,2}^{(k-1)},\, k>2n$ sur $U_{q}^{(n)}$.}
\noindent
Soit $k>2n$. On a $\dega{a_{1,2}^{(k-1)}}=-1$ et
$\degb{a_{1,2}^{(k-1)}}=-k$.
Donc,
$$
\forall\, j\in\lbrace 1,\ldots,2n\rbrace,
\quad a_{1,2}^{(k-1)}.u_j=0.
$$
Par ailleurs, si $P$ et $Q$ sont homog{\`e}nes
dans $U_{q}^{(n)}$, alors~:
$$
\varepsilon\Bigl[a_{1,2}^{(k-1)}.(PQ)\Bigr]=\sum_{\a+\b=k-1}
\varepsilon\bigl[ a_{1,1}^{(\a)}.P\bigr]\varepsilon\bigl[
a_{1,2}^{(\b)}.Q\bigr]+\sum_{\a+\b=k-1}
\varepsilon\bigl[ a_{1,2}^{(\a)}.P\bigr]\varepsilon\bigl[
a_{2,2}^{(\b)}.Q\bigr].
$$
Or, si $k>0,\, 
a_{1,1}^{(k)}\in I_{q}^{(n)}$ et $a_{2,2}^{(k)}\in I_{q}^{(n)}$.
Ceci vient du fait que $\dega{\bigl( a_{1,1}^{(k)}\bigr)}=0$ et
$\degb{\bigl( a_{1,1}^{(k)}\bigr)}=-k$.

\noindent
Donc,
$$
\varepsilon\Bigl[a_{1,2}^{(k-1)}.(PQ)\Bigr]=
\varepsilon\bigr[ a_{1,1}^{(0)}.P\bigr]
\varepsilon\bigr[ a_{1,2}^{(k-1)}.Q\bigr]
+\varepsilon\bigr[ a_{1,2}^{(k-1)}.P\bigr]
\varepsilon\bigr[ a_{2,2}^{(0)}.Q\bigr].
$$
Par suite, on montre par r{\'e}currence sur $\dega{P}$, que si $P$ est
homog{\`e}ne, alors $\varepsilon\bigr[ a_{1,2}^{(k-1)}.P\bigr]=0$
d'ou l'on d{\'e}duit que 
$a_{1,2}^{(k-1)}\in I_{q}^{(n)}$.
\subsection{Action de $a_{1,2}^{(k-1)},\, k\leq 2n$ sur
$U_{q}^{(n)}$.}\label{inco}
\noindent
Soit $k\leq 2n$.
On a toujours $\dega{a_{1,2}^{(k-1)}}=-1$ et
$\degb{a_{1,2}^{(k-1)}}=-k$.
Donc, 
\begin{equation}\label{saisp}
\exists C_k\, \forall i,\quad
a_{1,2}^{(k-1)}.u_i=C_k\delta_i^k.
\end{equation}
D'apr{\`e}s les propositions \ref{calculs+} et \ref{isofren},
on a $C_1\not= 0$.
De plus, en appliquant $a_{1,2}^{(k)}$ {\`a} (\ref{derui2}), on constate
que~:
$$
\exists r\not= 0\, \forall k,\quad C_{k+1}=rC_k.
$$
Par suite, $C_k\not= 0$. 
En particulier, $a_{1,2}^{(k-1)}\notin I_{q}^{(n)}$.
Le calcul montre que 
$$
r=-\displaystyle{1\over (q-1)^2}\quad\textnormal{ et }\quad 
C_1=q^{-{1\over 2}}.
$$
\subsection{ Action de $a_{1,1}^{(k)}$ et $a_{2,2}^{(k)},\,
k\geq 1$ sur $U_{q}^{(n)}$.}\label{bachi}
\noindent
Soit $k\geq 1$. Comme on l'a d{\'e}ja remarqu{\'e} pr{\'e}c{\'e}demment, 
on a $\dega{\bigl( a_{1,1}^{(k)}\bigr)}=0$ et
$\degb{\bigl( a_{1,1}^{(k)}\bigr)}=-k$.
Donc, $a_{1,1}^{(k)}\in I_{q}^{(n)}$ et~:
$$
\forall\, i\, \exists\, \lambda_{i,k},\quad
a_{1,1}^{(k)}.u_i=\lambda_{i,k}u_{i-k}.
$$
En appliquant $a_{1,1}^{(1)}$ aux deux membres de l'{\'e}quation
et en utilisant (\ref{aii}), on constate que $\lambda_{i,k}$
est une constante $\lambda_k$ ind{\'e}pendante de $i$.
Par ailleurs, en appliquant $a_{1,1}^{(k)}$ {\`a} (\ref{derui2}),
et en regardant le coefficient devant $u_1u_k$ suivant la base
$\displaystyle\prod\limits_{i=1}^{2n}u_{i}^{\a_{i}},\, \a_{i}\in\NN$
de $U_{q}^{(n)}$,
on peut calculer cette constante, et on constate que 
$\lambda_k=\displaystyle{C_k\over q-1}\not= 0$.

\noindent
De la m{\^e}me fa\c con, on montre que $a_{2,2}^{(k)}\in I_{q}^{(n)}$ et que~:
$$
\exists\, \mu_k\not= 0\, \forall\, i,\quad
a_{2,2}^{(k)}.u_i=\mu_k u_{i-k},
$$
avec $\mu_k=-\displaystyle{C_k\over q-1}\not= 0$.
\subsection{ Calcul de $I_q^{(n)}$.}
\noindent
Notons $I_0$ l'id{\'e}al {\`a} droite engendr{\'e} par les {\'e}l{\'e}ments~:
$$
a_{1,1}^{(0)}-1,\,
a_{2,2}^{(0)}-1,\, a_{2,1}^{(\a)},\, \a>0,\,
a_{1,2}^{(\b)},\, \b\geq 2n.
$$
D'apr{\`e}s ce qui pr{\'e}c{\`e}de, on a $I_0\subset I_q^{(n)}$.
Nous allons montrer qu'en fait, $I_0=I_q^{(n)}$.
Pour cela, nous aurons besoin de la s{\'e}rie de lemmes suivants~:
\begin{lem}\label{leplem}
Pour $r\in\lbrace 1,\ldots,2n\rbrace\,
(i,j)\in{\lbrace 1,2\rbrace}^2$, et $k\in \NN$, on a~:
$$
\forall\, P\in\CC[u_1,\ldots,u_r]_q,\,
a_{i,j}^{(k)}.P\in\CC[u_1,\ldots,u_r]_q.
$$
\end{lem}
\begin{dem}
Cela veut juste dire que $\CC[u_1,\ldots,u_r]_q$
est un $U_q(\bm)$-module. 
\end{dem}
\begin{lem}\label{atlem1}
$\forall k\in\lbrace 1,\ldots,2n\rbrace,\,
\forall P\in\CC[u_{k+1},\ldots,u_{2n}]_q,\,
a_{1,2}^{(k-1)}.P=0$.
\end{lem}
\begin{dem}
Par r{\'e}curence sur $\dega{P}$, et en utilisant la formule~:
$$
a_{1,2}^{(k-1)}.(PQ)=\displaystyle\sum\limits_{u=1}^{k}
\bigl( a_{1,1}^{(k-u)}.P\bigr)
\bigl( a_{1,2}^{(u-1)}.Q\bigr)
+\displaystyle\sum\limits_{u=1}^{k}
\bigl( a_{1,2}^{(u-1)}.P\bigr)
\bigl( a_{2,2}^{(k-u)}.Q\bigr).
$$
\end{dem}

\begin{lem}\label{atlem2}
Pour $r\in\NN,\, k\in\lbrace 1,\ldots,2n\rbrace,\,
\a\in\lbrace 1,\ldots,k-1\rbrace$ et\break
$i\in\lbrace 0,\ldots,r-1\rbrace$, il existe
$P_{i,\a}\in\CC[u_1,\ldots,u_{k-1}]_q$ tel que~:
$$
a_{2,2}^{(\a)}.u_{k}^{r}=\displaystyle\sum\limits_{i=0}^{r-1}
P_{i,\a} u_{k}^{i}.
$$
\end{lem}
\begin{dem}
La d{\'e}monstration se fait par r{\'e}currence sur $r$.
Soit\break $k\in\lbrace 1,\ldots,2n\rbrace$, et 
$\a\in\lbrace 1,\ldots,k-1\rbrace$.
On a~:
\begin{align*}
a_{2,2}^{(\a)}.u_{k}^{r+1}&=a_{2,2}^{(\a)}.(u_k u_{k}^{r})\\
&=\displaystyle\sum\limits_{u+v=\a}
\bigl( a_{2,1}^{(u)}.u_{k}\bigr)
\bigl( a_{1,2}^{(v)}.u_{k}^{r} \bigr)
+\displaystyle\sum\limits_{u+v=\a}
\bigl( a_{2,2}^{(u)}.u_{k}\bigr)
\bigl( a_{2,2}^{(v)}.u_{k}^{r} \bigr)
\end{align*}
Donc, d'apr{\`e}s le lemme pr{\'e}c{\'e}dent, le fait que 
$a_{2,1}^{(0)}=0$, et l'hypoth{\`e}se de r{\'e}currence,
on a avec les notations de la section \ref{bachi},
$$
a_{2,2}^{(\a)}.u_{k}^{r+1}=
q^{{r\over 2}}\mu_{\a}u_{k-\a}u_{k}^{r}+
\displaystyle\sum\limits_{x=1}^{\a}
\mu_{\a-x}u_{k+x-\a}\displaystyle\sum\limits_{i=0}^{r-1}
P_{i,x}u_k^i.
$$
Il suffit ensuite d'utiliser les relations de commutation
entre $u_k$ et $u_i,\, i<k$ (lorsque $x=\a$) pour obtenir
le r{\'e}sultat.
\end{dem}

\begin{lem}\label{atlem3}
Pour $k\in\lbrace 1,\ldots,2n\rbrace,\, j\in\lbrace 1,\ldots,k-1
\rbrace,\, r\in\NN$, et\break
$P\in\CC[u_1,\ldots,u_{k-1}]_q$,
ils existent des $P_i\in\CC[u_1,\ldots,u_{k-1}]_q$ tels que~:
$$
a_{1,2}^{(j-1)}.(P u_{k}^{r})=
\displaystyle\sum\limits_{i=0}^{r} P_i u_{k}^{i}.
$$
Ces $P_i$ sont uniques. De plus, 
$P_r=q^{{r\over 2}}a_{1,2}^{(j-1)}.P$.
\end{lem}
\begin{dem}
On a~:
$$
a_{1,2}^{(j-1)}.(P u_{k}^{r})=
\displaystyle\sum\limits_{\a+\b=j-1}
\bigl( a_{1,2}^{(\a)}.P\bigr)
\bigl( a_{2,2}^{(\b)}.u_{k}^{r}\bigr)+
\displaystyle\sum\limits_{\a+\b=j-1}
\bigl( a_{1,1}^{(\a)}.P\bigr)
\bigl( a_{1,2}^{(\b)}.u_{k}^{r}\bigr).
$$
La derniere somme est nulle en vertu du lemme \ref{atlem1}.
Donc,
$$
a_{1,2}^{(j-1)}.(P u_{k}^{r})=
\bigl( a_{1,2}^{(j-1)}.P\bigr)q^{{r\over 2}}u_k^r
+\displaystyle\sum\limits_{p=1}^{j-1}
\bigl( a_{1,2}^{(j-1-p)}.P\bigr)
\bigl( a_{2,2}^{(p)}.u_{k}^{r}\bigr).
$$
Il suffit ensuite d'appliquer les lemmes \ref{atlem2} 
et \ref{leplem} pour obtenir
l'existence des $P_i$.
L'unicit{\'e} est claire car
$\CC[u_1,\ldots,u_k]_q$ est un $\CC[u_1,\ldots,u_{k-1}]_q$-module
libre. 
\end{dem}

\begin{lem}\label{exder}
$\forall k\in{1,\ldots,2n}\, \forall l\in{\NN}^*\,\exists\,
C_{k,l}\not=0\quad
a_{1,2}^{(k-1)}. u_{k}^{l}=C_{k,l} u_{k}^{l-1}$.
\end{lem}
\begin{dem}
\noindent
On a $\Delta^{(l-1)}\bigl( a_{1,2}^{(k-1)}\bigr)=
\sum a_{i_1,i_2}^{(\a_1)}\otimes\ldots\otimes a_{i_l,i_{l+1}}^{(\a_l)}
$, la somme portant sur les indices $i_1,\ldots,i_{l+1},
\a_1,\ldots,\a_l$ tels que $i_1=1,i_{l+1}=2$ et
$\sum\a_i=k-1$.

\noindent
Soient $i_1,\ldots,i_{l+1},
\a_1,\ldots,\a_l$ de tels indices, et notons $r$ le plus petit entier
v{\'e}rifiant $i_r=2$, de telle sorte que $i_u=1$ si $u<r$.
D'apr{\`e}s (\ref{saisp}), on a $a_{i_{r}-1,i_r}^{(\a_{r-1})}.u_k=
a_{1,2}^{(\a_{r-1})}.u_k=0$ sauf si
$\a_{r-1}=k-1$, auquel cas on a 
$\forall j\not=r-1,\, \a_j=0$ car $\sum_v\a_v =k-1$.
Faisons l'hypoth{\`e}se qu'il existe $s\geq r$ tel que
$i_s=1$, et soit $s_0$ le plus petit entier $\geq r$ v{\'e}rifiant cette
propri{\'e}t{\'e}. Alors, on a $s_0> r$, et
$a_{i_{s_0-1},i_{s_0}}^{(\a_{s_0-1})}=a_{2,1}^{(0)}=0$.
Par suite, pour calculer 
$a_{1,2}^{(k-1)}.\bigl( u_{k}^{l}\bigr)$, on peut se restreindre {\`a}
sommer sur les indices $i_1,\ldots,i_{l+1},
\a_1,\ldots,\a_l$ tels qu'il existe $r\in\lbrace 1,\ldots,l\rbrace$
v{\'e}rifiant~:
$$
i_1=\ldots=i_r=1;\,\, i_{r+1}=\ldots=i_l=2;\,\, \a_r=k-1\,
\textnormal { et }
\forall j\not= r,\quad \a_j=0.
$$
Donc, avec les notations de la sous-section \ref{inco},

\begin{align}
\notag
a_{1,2}^{(k-1)}.u_k^l&=\sum_{r=1}^{l}
\bigl( a_{1,1}^{(0)}.u_k\bigr)\ldots
\bigl( a_{1,1}^{(0)}.u_k\bigr)
\underbrace{a_{1,2}^{(k-1)}.u_k}_{r^{\textnormal{{\`e}me}}\textnormal
{ t{\`e}rme}}
\bigl( a_{2,2}^{(0)}.u_k\bigr)\ldots
\bigl( a_{2,2}^{(0)}.u_k\bigr)\\
\notag
&=\sum_{r=1}^{l}\Bigl(
q^{-{r-1\over 2}}u_k^{r-1}\Bigr)C_k\Bigl(
q^{{1\over 2}(l-r)}u_k^{l-r}\Bigr)\\
\label{jor3}
&=q^{{(l+1)\over 2}}\Bigl(
\sum_{r=1}^{l}q^{-r}\Bigr) C_k u_k^{l-1}
\end{align}
D'o{\`u} le r{\'e}sultat, car $q$ n'est pas une racine de l'unit{\'e}.
\end{dem}

\begin{lem}\label{atlem4}
Soit $k\in\lbrace 1,\ldots,2n\rbrace$. Alors,
$\forall \a_k\, \exists\, \Lambda_{{\a}_k}\not= 0$ tel que~: 
$$
\forall\, \b_k\, \forall P\in\CC[u_1,\ldots,u_{k-1}]_q,\quad
\varepsilon\Bigl[
{\bigl( a_{1,2}^{(k-1)}\bigr)}^{{\a}_k}.\bigl( P u_k^{{\b}_k}\bigr)\Bigr]
=\Lambda_{{\a}_k}\varepsilon(P)\delta_{{\a}_k}^{{\b}_k}.
$$
\end{lem}
\begin{dem}
On peut supposer que $P$ est un mon{\^o}me en les $u_i$,\break
$1\leq i\leq k-1$.
Posons $n_i=\degi{P},\, i\in\lbrace 1,2\rbrace$.
Si $\varepsilon\Bigl[
{\bigl( a_{1,2}^{(k-1)}\bigr)}^{{\a}_k}.\bigl( P u_k^{{\b}_k}\bigr)\Bigr]
\not= 0$, alors~: 
 $$
  \degi{{a_{1,2}^{(k-1)}}^{{\a}_k}}
   +\degi{P u_k^{{\b}_k}}=0,
 $$ 
pour $i\in\lbrace 1,2\rbrace$.
Donc,
\begin{align*}
n_1+{\b}_k&={\a}_k\\
n_2+k{\b}_k&=k{\a}_k
\end{align*}
D'o{\`u} $n_2=k n_1$.
Mais d'apr{\`e}s le lemme \ref{atlem0}, ceci entraine $n_1=n_2=0$.
Donc $P$ est une constante et ${\a}_k={\b}_k$.
R{\'e}ciproquement, si ${\a}_k={\b}_k$ et $P=\varepsilon(P)$,
en utilisant plusieurs fois le lemme \ref{exder}, on a~:
$$
\varepsilon\Bigl[
{\bigl( a_{1,2}^{(k-1)}\bigr) }^{{\a}_k}.\bigl( P u_k^{{\a}_k}\bigr)\Bigr]=
\varepsilon(P)\varepsilon\Bigl[
{\bigl( a_{1,2}^{(k-1)}\bigr) }^{{\a}_k}. 
u_k^{{\a}_k}\Bigr]=\varepsilon(P)\Lambda_{{\a}_k},
$$
avec $\Lambda_{{\a}_k}\not= 0$.
\end{dem}

\noindent
Id{\'e}alement, on aimerait obtenir une base duale
dans l'alg{\`e}bre engendr{\'e}e par les $a_{1,2}^{(k-1)}$, 
(avec $k\in\lbrace 1,\ldots,2n\rbrace$) 
aux $\displaystyle\prod\limits_{i=1}^{2n}u_i^{\a_i}$
pour la forme bilin{\'e}aire~:
$$
\begin{array}{rcl}
\CC[\BP]_q/I_0\times U_q^{(n)}&\longrightarrow&\CC\\
(\bar{x},P)&\longmapsto&\varepsilon(x.P).
\end{array}
$$
Mais ceci n'est pas tout {\`a} fait {\'e}vident car
on a par exemple les {\'e}galit{\'e}s~:
$$
\varepsilon\Bigl[
{\bigl( a_{1,2}^{(2)}\bigr)}^2.\bigl( u_2 u_4\bigr)
\Bigr]\not= 0
$$
et
$$
\varepsilon\Bigl[\bigl( a_{1,2}^{(1)} a_{1,2}^{(3)}\bigr). 
\bigl( u_2 u_4\bigr)\Bigr]\not= 0.
$$
Donc, il serait illusoire de croire que les 
$\displaystyle\prod\limits_{i=1}^{2n}
{\bigl( a_{1,2}^{(i-1)}\bigr)}^{\a_i}$
est la base duale des 
$\displaystyle\prod\limits_{i=1}^{2n}u_i^{\a_i}$
comme pourrait nous le laisser supposer la formule
(\ref{saisp}).

Pour se tirer d'affaires, consid{\'e}rons $\leq$
la relation d'ordre lexicographique naturelle sur ${\NN}^{2n}$,
$R$ le ``renversement'' de suites~:
  $$
  \begin{array}{rcl}
  R:\quad {\NN}^{2n}&\longrightarrow&{\NN}^{2n}\\
  (a_i)&\longmapsto&(a_{2n-i+1})
  \end{array}
  $$
et $\preceq$ la relation d'ordre lexicographique invers{\'e}e d{\'e}finie
par~:
$$
\forall \bigl( {\underline a}\, {\underline b}\bigr)
\in {\bigl( {\NN}^{2n}\bigr)}^{2},\quad
{\underline a}\preceq {\underline b}\,\Longleftrightarrow\,
R({\underline a})\leq R({\underline b}).
$$

Prenons l'exemple $n=2$. La suite des coefficients 
correspondant {\`a} 
$u_2 u_4$ est $\b=(0,1,0,1):\, u_2 u_4=u_1^0u_2^1u_3^0u_4^1$.
La suite des coefficients correspondant {\`a} 
${\bigl( a_{1,2}^{(2)}\bigr)}^2$ est 
$\a=(0,0,2,0):\, {\bigl( a_{1,2}^{(2)}\bigr)}^2=
{\bigl( a_{1,2}^{(0)}\bigr)}^0
{\bigl( a_{1,2}^{(1)}\bigr)}^0
{\bigl( a_{1,2}^{(2)}\bigr)}^2
{\bigl( a_{1,2}^{(3)}\bigr)}^0$.
On a $\a\prec\b$ et
$\varepsilon\Bigl[ \displaystyle\prod
{\bigl( a_{1,2}^{(i-1)}\bigr)}^{\a_i}.
\displaystyle\prod u_i^{\b_i}\Bigr]\not= 0$

Le lemme suivant montre que ce genre de choses ne peut pas se produire
si $\b\preceq\a$.

 \begin{lem}\label{esper}
 Soient $k\in\lbrace 1,\ldots,2n\rbrace$ et 
 $\a\in{\NN}^{k}$. Alors, il existe $M_{\a}\in{\CC}^{*}$ tel que~:
 $$
  \forall \b\preceq\a,\quad
  \varepsilon\Bigl[ \displaystyle\prod\limits_{i=1}^{k}
  {\bigl( a_{1,2}^{(i-1)}\bigr)}^{\a_i}.
  \displaystyle\prod\limits_{i=1}^{k} 
  u_i^{\b_i}\Bigr]=M_{\a}\delta_{\a}^{\b}.
 $$
 \end{lem}

 \begin{dem}
 Par r{\'e}currence sur $k$. Si $k=1$, on applique le lemme \ref{atlem4}
 avec $P=1$ et $k=1$. Supposons le r{\'e}sultat vrai jusqu'{\`a} l'ordre
 $k-1\geq 1$. Soient $\a$ et $\b$ deux suites finies dans
 ${\NN}^{k}$, tels que $\b\preceq\a$.
 Posons~:
 $$
  P=u_1^{\b_1}\ldots u_{k-1}^{\b_{k-1}}.
 $$
 En vertu du lemme \ref{atlem3}, on a~:
  $$
  a_{1,2}^{(0)}.(P u_k^{\b_k})=\displaystyle\sum\limits_{i=0}^{\b_k}
  P_{i,0} u_{k}^{i},
  $$
 avec $P_{i,0}\in\CC[u_1,\ldots,u_{k-1}]_q$ et
 $P_{\b_k,0}=q^{{\b_k\over 2}} a_{1,2}^{(0)}.P$.
 
 En appliquant de mani{\`e}re r{\'e}p{\'e}t{\'e}e $a_{1,2}^{(0)},\ldots,
 a_{1,2}^{(k-2)}$ autant de fois qu'il le
 faut aux deux membres de l'{\'e}quation pr{\'e}c{\'e}dente, 
 et en utilisant les lemmes
 \ref{leplem} et \ref{atlem3}, on voit qu'ils existent des $Q_i$
 uniques tels que~:
  \begin{align*}
   Q_i&\in\CC[u_1,\ldots,u_{k-1}]_q,\\
   Q_{\b_k}&=q^{{\b_k\over 2}(\a_1+\ldots+\a_{k-1}})
   {\bigl( a_{1,2}^{(0)}\bigr)}^{\a_1}\ldots
   {\bigl( a_{1,2}^{(k-2)}\bigr)}^{\a_{k-1}}.P
  \end{align*}
 et~:
  $$
   {\bigl( a_{1,2}^{(0)}\bigr)}^{\a_1}\ldots
   {\bigl( a_{1,2}^{(k-2)}\bigr)}^{\a_{k-1}}.(P u_k^{\b_k})
   =\displaystyle\sum\limits_{i=0}^{\b_k}
   Q_i u_k^i.
  $$
 D'apr{\`e}s (\ref{aii}), les $a_{1,2}^{(k-1)}$ commutent entre eux.
 Donc, en reprenant la notation de la constante $\Lambda_{\a_k}\not= 0$
 introduite dans le lemme \ref{atlem4}, on a~:
 \begin{align*}
  \varepsilon\Bigl[ \displaystyle\prod\limits_{i=1}^{k}
  {\bigl( a_{1,2}^{(i-1)}\bigr)}^{\a_i}.
  \displaystyle\prod\limits_{i=1}^{k} 
  u_i^{\b_i}\Bigr]&=
  \varepsilon\Bigl[
  {\bigl( a_{1,2}^{(k-1)}\bigr)}^{\a_k}
  {\bigl( a_{1,2}^{(0)}\bigr)}^{\a_1}
  \ldots
  {\bigl( a_{1,2}^{(k-2)}\bigr)}^{\a_{k-1}}.(P u_k^{\b_k})
  \Bigr]\\
  &=\displaystyle\sum\limits_{i=0}^{\b_k}
  \varepsilon\Bigl[ {\bigl( a_{1,2}^{(k-1)}\bigr)}^{\a_k}.
  (Q_i u_k^i)\Bigr]\\
  &=\Lambda_{\a_k}\varepsilon (Q_i)\delta_i^{\a_k}.
  \end{align*}
 Or, $\b\preceq\a$. Donc, $\b_k\leq\a_k$.
 Donc, 
  \begin{align*}
  \varepsilon\Bigl[ \displaystyle\prod\limits_{i=1}^{k}
  {\bigl( a_{1,2}^{(i-1)}\bigr)}^{\a_i}.
  \displaystyle\prod\limits_{i=1}^{k} 
  u_i^{\b_i}\Bigr]&=
  \Lambda_{\a_k}\varepsilon (Q_{\b_k})\delta_{\a_k}^{\b_k}\\
  &=q^{{\b_k\over 2}(\a_1+\ldots+\a_{k-1}})
  \Lambda_{\a_k}\varepsilon\Bigl[ 
  \displaystyle\prod\limits_{i=1}^{k-1}
  {\bigl( a_{1,2}^{(i-1)}\bigr)}^{\a_i}.
  \displaystyle\prod\limits_{i=1}^{k-1} 
  u_i^{\b_i}\Bigr]\delta_{\a_k}^{\b_k}
  \end{align*}
 Si $\a_k\not=\b_k$, cette quantit{\'e} est nulle comme attendue.
 Sinon, on a~:
  $$
   (\b_1,\ldots,\b_{k-1})\preceq  (\a_1,\ldots,\a_{k-1}),
  $$
 et on applique l'hypoth{\`e}se de r{\'e}currence.
 \end{dem}
 
\noindent
Nous pouvons {\`a} pr\'sent en d{\'e}duire la proposition~:
\begin{prop}\label{lapropx}
L'id{\'e}al {\`a} droite $I_q^{(n)}$ est engendr{\'e} ({\`a} droite)
par les {\'e}l{\'e}ments~:
$a_{1,1}^{(\a)},a_{2,2}^{(\a)},
a_{2,1}^{(\a)},a_{1,2}^{(2n-1+\a)},
a_{1,1}^{(0)}-1,a_{2,2}^{(0)}-1$, avec $\a>0$.
Autrement dit, $I_q^{(n)}=I_0$.
\end{prop}
\begin{dem}
 Soit $x\in I_q^{(n)}$. Montrons que $x\in I_0$.
 Pour tous entiers $(i,j,k,l)\in{\NN}^4$,
 et toutes suites d'entiers
 $\a=(\a_r),\, 0\leq r\leq i,\, \b=(\b_s),\, 0\leq s\leq j,\,
 \gamma=(\gamma_t),\, 0\leq t\leq k,\,
 \delta=(\delta_u),\, 0\leq u\leq l$, on pose~:
  $$
   x_{\a,\b,\gamma,\delta}^{i,j,k,l}=
   \displaystyle\prod\limits_{r=1}^{i}
   {\bigl( a_{2,1}^{(r)}\bigr)}^{\a_r}
   \displaystyle\prod\limits_{s=0}^{j}
   {\bigl( a_{1,1}^{(s)}\bigr)}^{\b_s}
   \displaystyle\prod\limits_{t=0}^{k}
   {\bigl( a_{2,2}^{(t)}\bigr)}^{\gamma_t}
   \displaystyle\prod\limits_{u=0}^{l}
   {\bigl( a_{1,2}^{(u)}\bigr)}^{\delta_u}.
  $$
 D'apr{\`e}s la d{\'e}finition de $\CC[\BP]_q$, on peut voir que
 $x$ est combinaison lin{\'e}aire de termes de cette forme.
 Etant donn{\'e} $i,j,k,l,\a,\b,\gamma,\delta$, chacun des cas suivants
 entraine $x_{\a,\b,\gamma,\delta}^{i,j,k,l}\in I_0$~:
  \begin{itemize}
   \item Il existe $r$ tel que $\a_r>0$,\\
   \item Il existe $s$ tel que $s\b_s>0$,\\
   \item Il existe $t$ tel que $t\gamma_t>0$,\\
   \item Il existe $u\in\lbrace 2n,\ldots,l\rbrace$ tel que
         $\delta_u\not= 0$.
  \end{itemize}
 Par ailleurs, si $i\in\lbrace 1,2\rbrace,\, a_{i,i}^{(0)}-1\in I_0$.
 Donc, la formule du bin{\^o}me montre que l'on peut supposer
 que $x$ est en fait une
 combinaison lin{\'e}aire de termes de la forme
 $\displaystyle\prod\limits_{i=1}^{2n}
   {\bigl( a_{1,2}^{(i-1)}\bigr)}^{\a_i}$.
 Pour conclure, il suffit de montrer que $x=0$.
 Supposons le contraire, et posons~:
  $$
   x=\displaystyle\sum\limits_{\a\in E}\lambda_{\a}
     \displaystyle\prod\limits_{i=1}^{2n}
     {\bigl( a_{1,2}^{(i-1)}\bigr)}^{\a_i},
  $$
 o{\`u} $E$ est un ensemble fini non vide tel que $\forall\a\in E\,
 \lambda_{\a}\not= 0$.
 Soit $\a$ le plus petit {\'e}l{\'e}ment de l'ensemble $E$ pour la
 relation d'ordre $\preceq$.
 Alors, en utilisant le lemme \ref{esper}, et en regardant l'action de
 $x$ sur $\displaystyle\prod\limits_{i=1}^{2n} u_i^{\a_i}$,
 on constate que $\lambda_{\a}=0$. Contradiction.
\end{dem}

\medskip
\noindent
{\bf Remarques:}

\smallskip
\noindent
\begin{enumerate}

\item
Sous forme ramass{\'e}e, l'action de $U_q(\bm)$ sur les $U_i$
est la suivante~:
\begin{align*}
e_{+}.u(\lambda) &=-{\lambda}^{-1}u^{2}(\lambda),\\
e_{-}.u(\lambda) &=q^{-1},\\
k.u(\lambda) &=q u(\lambda).
\end{align*}
avec $u(\lambda)=\displaystyle\sum\limits_{i=0}^{+\infty}
(-1)^{i}u_{i+1}\lambda^{-i}$.

\item
Soit 
$$
\begin{array}{rcl}
T:\quad \CC[u_1,\ldots,u_{2n-1}]_q&\longrightarrow&
\CC[u_1,\ldots,u_{2n}]_q=U_q^{(n)}\\
u_i&\longrightarrow&u_{i+1}
\end{array}
$$
l'application de translation sur $U_q^{(n)}$. On peut v{\'e}rifier que $T$
est bien un morphisme d'alg{\`e}bres.

\noindent
Alors, pour
$i\in\lbrace 1,2\rbrace,\,
k\in\lbrace 1,\ldots,2n\rbrace$ et 
$x\in\CC[u_{k},\ldots,u_{2n-1}]_q$ on a 
$$
a_{i,i}^{(k-1)}\circ T.x=T\circ a_{i,i}^{(k-1)}.x
$$
De plus, en reprenant la constante $r$ de \ref{inco},
on montre que si $k\in\lbrace
1,\ldots,2n-1\rbrace$ et $x\in\CC[u_{k},\ldots,u_{2n-1}]_q$,
alors~:
$$
a_{1,2}^{(k)}\circ T.x=r\,T\circ a_{1,2}^{(k-1)}.x
$$

\item On doit pouvoir remplacer $\preceq$ par $\leq$
dans le lemme \ref{esper}

\item En reprenant la formule (\ref{jor3}), il semblerait que lorsque
$q$ est une racine primitive $l-$i{\`e}me de l'unit{\'e}, l'id{\'e}al $I$
soit engendr{\'e}
par $I_0$ et les ${\bigl( a_{1,2}^{(k)}\bigr)}^{l},$ avec $k\in\NN$.
\end{enumerate}
\subsection{Le couplage non d{\'e}g{\'e}n{\'e}r{\'e}.}

\noindent
Nous allons nous int{\'e}resser {\`a} pr{\'e}sent au couplage~:
  $$
   \begin{array}{rcl}
   B:\quad {U_q(\bm)}/I_q^{(n)} \times U_q^{(n)}&\longrightarrow&\CC\\
                           (\bar{x},P)&\longmapsto&
                                       \varepsilon(x.P)
   \end{array}
  $$
D'apr{\`e}s la forme de $I_0=I_q^{(n)}$,
il est clair que $I_q^{(n)}$ est un id{\'e}al de Hopf bigradu{\'e}.
Or, l'action de $U_q(\bm)$ sur $U_q^{(n)}$ est bigradu{\'e}e.
Donc, ${U_q(\bm)}/I_q^{(n)}$ est bigradu{\'e}, ainsi que le couplage 
$B$.
De plus,
la structure de $U_q(\bm)$-alg{\`e}bre-module sur
$U_q^{(n)}$, et le fait que $\varepsilon$ soit un morphisme 
d'alg{\`e}bres entra{\^\i}ne que l'application provenant de $B$~:
$$
\begin{array}{rcl}
\xi:\quad U_q^{(n)}&\longrightarrow&{\Bigl( {U_q(\bm)}/I_q^{(n)}\Bigr)}^*\\
P&\longmapsto&\bigl[
x\rightarrow\varepsilon(x.P)\bigr]
\end{array}
$$
est un morphisme d'alg{\`e}bres.
\begin{prop}\label{injec}
L'application $\xi$ est injective.
\end{prop}
 \begin{dem}
  Soit $P\in  U_q^{(n)},\, P\not= 0$. Alors, nous avons vu 
  qu'il existe $F$ un ensemble fini non vide tel que~:
   $$
    P=\displaystyle\sum\limits_{\b\in F}\lambda_{\b}
      \displaystyle\prod\limits_{i=1}^{2n} u_i^{\b_i},
   $$ 
 avec $\forall \b\in F,\, \lambda_{\b}\not= 0$.
 Soit $\a\in F$ le plus grand {\'e}l{\'e}ment de $F$ pour la relation
 $\preceq$,
 et posons~:
  $$
   x=\displaystyle\prod\limits_{i=1}^{2n} 
       {\bigl( a_{1,2}^{(i-1)}\bigr)}^{\a_i}.
  $$
 Alors, le lemme \ref{esper} montre que~:
   \begin{align*}
    \xi(P).x&=\varepsilon(x.P)\\
            &=\displaystyle\sum\limits_{\b\in F}\lambda_{\b}
              \varepsilon\Bigl[ 
              x.\displaystyle\prod\limits_{i=1}^{2n} 
              u_i^{\b_i}\Bigr]\\
            &=\lambda_{\a}M_{\a}\not= 0
   \end{align*}
 ce qui prouve l'injectivit{\'e}.
 \end{dem}
\begin{prop}\label{noyau}
L'application lin{\'e}aire~:
 $$
 \begin{array}{rcl}
 l:\quad {U_q(\bm)}/I_q^{(n)}&\longrightarrow&{U_q^{(n)}}^*\\
            \bar{x}&\longmapsto&\varepsilon(x.\,)
 \end{array}
 $$
est injective.
\end{prop}
\begin{dem}
Cela d{\'e}coule de la d{\'e}finition m{\^e}me de $I_q^{(n)}$~!
\end{dem}

\noindent
Le couplage $B$ {\'e}tant (bi)gradu{\'e}, on peut remplacer
le dual normal utilis{\'e} jusqu'ici par le dual restreint 
(ou gradu{\'e}). De plus, toutes les composantes gradu{\'e}es
sont de dimension finie.
Les deux propositions pr{\'e}c{\'e}dentes entrainent alors imm{\'e}diatement
que~:
  $$
   U_q^{(n)}={\bigl( {U_q(\bm)}/I_q^{(n)}\bigr)}'
  $$ 
et~: 
$$
{U_q(\bm)}/I_q^{(n)}={\bigl( U_q^{(n)}\bigr)}',
$$ 
en notant par ' le dual
restreint.
Notons qu'une base du $\CC$ espace vectoriel ${U_q(\bm)}/I_q^{(n)}$
est donn{\'e}e par les 
$\textnormal{cl}\Bigl[\displaystyle\prod\limits_{i=1}^{2n} 
{\bigl( a_{1,2}^{(i-1)}\bigr)}^{\a_i}\Bigr]$.
Comme espace-vectoriel, $U_q^{(n)}$ et ${U_q(\bm)}/I_q^{(n)}$ 
sont tous les deux
isomorphes {\`a} $\CC[X_1,\ldots,X_{2n}]$.

\noindent
Ceci ach{\`e}ve la d{\'e}monstration de T4.
En ce qui concerne T5, il suffit d'appliquer la proposition
\ref{lapropx} ainsi que l'isomorphisme entre $\CC[\BP]_q$
et les nouvelles r{\'e}alisations mis en {\'e}vidence dans [13].

\bigskip
\noindent
{\bf RETOUR SUR ${\cal A}_q^{(n)}$}

\bigskip
\noindent
De la m{\^e}me fa\c con que ${\cal U}_q^{(n)}$ est une d{\'e}formation
quantique de ${\cal U}^{(n)}$, on aimerait pouvoir affirmer que
${\cal A}_q^{(n)}$ est une d{\'e}formation quantique de
${\cal A}^{(n)}$. Comme nous l'avons d{\'e}ja remarqu{\'e}, il faudrait
imposer en plus les relations entre les $a_{i,j}^{(k)}$
et ${a_{1,1}^{(n-1)}}',\,{a_{1,2}^{(n-1)}}',\,  
{a_{2,1}^{(n)}}',\, {a_{2,2}^{(n)}}'$.
Mais m{\^e}me apr{\`e}s cela, on ne voit pas trop comment on pourrait
d{\'e}montrer que ${\cal A}_q^{(n)}$ est une d{\'e}formation quantique de
${\cal A}^{(n)}$, et que l'application $\varphi$ de \ref{E4}
est injective (comme c'est le cas dans le cas classique).

Cela dit, par analogie avec le cas classique, il est plus naturel
d'introduire ${\cal A}_q^{(n)}$ comme un localis{\'e} de
${\bar{{\cal A}}}_q^{(n)}$~: le localis{\'e} suivant la partie
multiplicative $S$ engendr{\'e}e par
$a_{1,1}^{(n-1)},\ldots,a_{2,2}^{(n)}$.
Mais, pour pouvoir localiser dans un anneau non commutatif, il faut
r{\'e}unir trois conditions exprim{\'e}es dans [14].
Les deux premi{\`e}res sont facilement v{\'e}rifi{\'e}es sauf la troisi{\`e}me
qui exprime que $S$ ne contient pas de diviseurs de $0$.

Par ailleurs, on peut voir que l'application~:
$$
\varphi:\quad {\bar{{\cal A}}}_q^{(n)}\longrightarrow A_q^{(n)}
$$
d{\'e}finie en \ref{E4} est injective si 
${\bar{{\cal A}}}_q^{(n)}$ est int{\`e}gre.
Il est donc naturel de conjecturer que ${\bar{{\cal A}}}_q^{(n)}$ 
est int{\`e}gre.

Etant donn{\'e} l'isomorphisme entre $U_q(\nm)$ et 
${\cal A}_q$, dire que ${\bar{{\cal A}}}_q$ est int{\`e}gre,
c'est dire que l'alg{\`e}bre engendr{\'e}e par les op{\'e}rateurs
d'{\'e}crants $\Sigma^{+}$ et $\Sigma^{-}$ est isomorphe {\`a}
${\bar{{\cal A}}}_q^{(n)}$.

D'autre part, notons que l'application~:
$$
\begin{array}{rcl}
\iota:\quad {\cal A}^{(n)}&\longrightarrow&{\CC}^{2n-2}\times
{({\CC}^*)}^2\\
\begin{pmatrix}
a&b\\
c&d
\end{pmatrix}&\longmapsto&\bigl( c_1,\ldots,c_{n-1},d_1,\ldots,
d_{n-1},c_n,d_n\bigr)
\end{array}
$$
est injective car $(a,b)$ est le couple de Bezout unique tel que
$du-cv=1$ et $\deg{u}<\deg{c},\, \deg{v}<\deg{d}$.
De plus, l'image de ${\cal A}^{(n)}$ par $\iota$ est 
$\textnormal{Res}^{<-1>}({\CC}^*)$ o{\`u} Res est l'application 
R{\'e}sultant~:
$$
\begin{array}{rcl}
\textnormal{Res}:\quad {\CC}^{2n-2}\times {({\CC}^*)}^2
&\longrightarrow&\CC\\
\bigl( c_1,\ldots,c_{n-1},d_1,\ldots,d_{n-1},c_n,d_n\bigr)
&\longmapsto&\textnormal{Res}(c,d),
\end{array}
$$
avec $c=\displaystyle\sum\limits_{i=1}^{n}c_i\lambda^{i}$,
et $d=1+\displaystyle\sum\limits_{i=1}^{n}d_i\lambda^{i}$.
Par suite, $\CC[{\cal A}^{(n)}]$ s'identifie au localis{\'e}~:
\begin{equation}\label{res}
\CC[{\cal A}^{(n)}]={\CC[c_i,d_i,\, 1\leq i\leq n-1][c_n^{\pm
  1},d_n^{\pm 1}]}_{(\textnormal{Res})}
\end{equation}
o{\`u} $(\textnormal{Res})$ est l'id{\'e}al engendr{\'e} par le r{\'e}sultant 
"g{\'e}n{\'e}rique'' de $c$ et de $d$.
En outre, il est facile de voir que la restriction de 
$\varphi$ {\`a} la sous alg{\`e}bre de ${\bar{{\cal A}}}_q^{(n)}$
engendr{\'e} par les $a_{2,1}^{(k)},\, i\in\lbrace 1,2\rbrace,\,
k\in\NN$ est injective. Donc l'analogue quantique
${\CC[c_i,d_i,\, 1\leq i\leq n]}_q$ est int{\`e}gre.
Il serait donc satisfaisant de pouvoir introduire un r{\'e}sultant
quantique pour g{\'e}n{\'e}raliser (\ref{res}), et
d{\'e}crire la cellule de Schubert quantique.

De mani{\`e}re g{\'e}n{\'e}rale, on conjecture que si
$\NP$ est la partie nilpotente positive de $\widehat{\textnormal{SL}_n}$,
et si $p$ est un entier positif, alors le quotient
de $\CC[\NP]_q$ par l'id{\'e}al engendr{\'e}
par les $a_{i,j}^{(k)}$ pour $k\geq p$ est int{\`e}gre.

\bigskip
\bigskip
\bigskip
\noindent
\centerline{\bf{REMERCIEMENTS}}

\bigskip
\bigskip
\noindent
Tous mes remerciements vont {\`a} mon directeur de th{\`e}se 
B. Enriquez.

\newpage
\appendix
\noindent
{\centerline {\bf {\huge APPENDICE}}}
\section{Sur l'alg{\`e}bre $A_n$ de la proposition \ref{fracq2}}
\noindent
Rigoureusement les alg{\`e}bres ${A\{ t_1,\ldots,t_n \}}/I$
et ${A\{\{ t_1,\ldots,t_n \}\}}/I$ sont les suivantes~:

\medskip
\noindent
$\bullet$
  L'espace $A\{ t_1,\ldots,t_n \}$
  est par d{\'e}finition le $A$-module libre
  dont une base est donn{\'e}e par l'ensemble $S$ des suites finies
  {\`a} valeurs dans $\{ 1,\ldots,n\}$. Si $s$ et $s'$ sont dans $S$,
  on d{\'e}finit de mani{\`e}re {\'e}vidente le produit de $s$ avec $s'$
  que l'on note $ss'$, et qui n'est rien d'autre que la juxtaposition
  de $s$ avec $s'$.
  
  \smallskip
  \textnormal{\bf Fait 1:}\quad
  $\forall t\in S,\; \bigl\lbrace (s,s')\in S^2\, /\, t=ss'\bigr\rbrace\;
  \textnormal{est un ensemble fini.}$
  \smallskip
  
  On peut donc {\'e}tendre sans probl{\`e}me le produit dans $S$ {\`a} un
  produit dans $A\{ t_1,\ldots,t_n \}$ par~:
  $$
  \forall x=\displaystyle\sum\limits_{s\in S} a_s s,\; 
  \forall x'=\displaystyle\sum\limits_{s'\in S}a'_{s'},\quad   
  xx':=\displaystyle\sum\limits_{t\in S}\Bigl(
  \displaystyle\sum\limits_{ss'=t} a_s a'_{s'}\Bigr)t.
  $$ 
  L'unit{\'e} de $A\{ t_1,\ldots,t_n \}$ est la suite vide 
  $\emptyset$. Le $A$-module $A\{ t_1,\ldots,t_n \}$ muni de ce
  produit
  est une alg{\`e}bre. C'est l'alg{\`e}bre
  (libre) engendr{\'e}e par la suite $\emptyset\, (=1)$
  et les suites {\'e}l{\'e}mentaires $t_i=(i)$.
  
  L'ensemble des graduations sur $A\{ t_1,\ldots,t_n \}$ est en
  bijection avec ${\ZZ}^n$. Nous consid{\`e}rerons celle d{\'e}finie par~:
  $$
  \forall i,\quad \deg{t_i}=1.
  $$
  Elle munit $A\{ t_1,\ldots,t_n \}$ d'une structure d'alg{\`e}bre 
  valu{\'e}e. On note $v$ la valuation.

  \medskip
  \noindent
  $\bullet$ 
  On montre facilement que l'alg{\`e}bre quotiente ${A\{ t_1,\ldots,t_n
  \}}/I$ est un $A$-module libre dont une base est donn{\'e}e par la
  famille $t_1^{\a_1}\ldots t_n^{\a_n}$ avec $\a_i\in\NN$. C'est donc
  le $A$-module libre dont une base est donn{\'e}e par l'ensemble $N$
  des $n$-uplets {\`a} valeurs dans $\NN$.
  
  Dans cette identification, la suite {\'e}l{\'e}mentaire $t_i$
  correspond au $n$-uplet~: 
   $$
    (0,\ldots,0,1,0,\ldots,0),
   $$
  o{\`u} le $1$ est en $i^{\textnormal{{\`e}me}}$ position, et le produit
  de deux $n$-uplets $\underline{\a}=(\a_1,\ldots,\a_n)$
  avec $\underline{\b}=(\b_1,\ldots,\b_n)$ est d{\'e}fini par~:
  $$
  \underline{\a} . \underline{\b}=q^{-(\a_2\b_1+\ldots+\a_n\b_{n-1})}
  \underline{\a+\b},
  $$ 
  avec 
  $$
  \underline{\a+\b}=(\a_1+\b_1,\ldots,\a_n+\b_n).
  $$
  La graduation introduite pr{\'e}c{\'e}demment sur $A\{ t_1,\ldots,t_n
  \}$
  passe au quotient, faisant de ${A\{ t_1,\ldots,t_n\}}/I$ une
  alg{\`e}bre valu{\'e}e. De plus, la surjection canonique~:
  $$
  \pi:\quad A\{ t_1,\ldots,t_n \}\longrightarrow {A\{ t_1,\ldots,t_n\}}/I
  $$ 
  est telle que si $x\in A\{ t_1,\ldots,t_n \}$, alors, 
  $v\bigl(\pi(x)\bigr)\geq v(x)$. Donc, $\pi$ est continue.
 
  \medskip
  \noindent
  $\bullet$ 
  L'espace $A\{\{ t_1,\ldots,t_n \}\}$ est par d{\'e}finition, le
  $A$-module form{\'e} par  les fonctions de $S$ dans $A$.
  Si $f$ et $g$ sont dans $A\{\{ t_1,\ldots,t_n \}\}$, on d{\'e}finit
  le produit 
  $fg$ par~:
  $$
  \begin{array}{rcl}
  fg:\quad S&\longrightarrow&A\\
  t&\longmapsto&\displaystyle\sum\limits_{ss'=t}f(s)g(s').
  \end{array}
  $$
  Ceci d{\'e}finit une structure d'alg{\`e}bre sur $A\{\{ t_1,\ldots,t_n \}\}$.
  De plus, l'application~:
  $$
    \begin{array}{rcl}
      S&\longrightarrow&A\{\{ t_1,\ldots,t_n \}\}\\
      s&\longmapsto&\delta_s=\left\{\begin{array}{rcl}
                               S&\longrightarrow&A\\ 
                               t&\longmapsto&\delta_s^{\, t}
                              \end{array}
                              \right.
    \end{array}
  $$
  o{\`u} $\delta_s^{\, t}$ est le symbole de Kronecker
  s'{\'e}tend en un monomorphisme d'alg{\`e}bres de $A\{ t_1,\ldots,t_n \}$
  dans $A\{\{ t_1,\ldots,t_n \}\}$  car 
  $\forall s,s'\in S,\; \delta_{s}\delta_{s'}=\delta_{ss'}$ et de
  plus,
  l'unit{\'e} de  
  $A\{\{ t_1,\ldots,t_n \}\}$ est clairement $\delta_{\emptyset}$.
  Ceci nous permet d'identifier $A\{ t_1,\ldots,t_n \}$ {\`a} la
  sous-alg{\`e}bre de $A\{\{ t_1,\ldots,t_n \}\}$ form{\'e}e par les
  fonctions {\`a} support fini.

  On peut {\'e}tendre {\`a} $A\{\{ t_1,\ldots,t_n \}\}$ la valuation
  d{\'e}finie sur $A\{ t_1,\ldots,t_n \}$ par~:
  $$
  \forall f\in A\{\{ t_1,\ldots,t_n \}\},\quad
  v(f)=\textnormal{Inf}\bigl\lbrace
  v(s)/\, s\in S\textnormal{ et }f(s)\not=0\bigr\rbrace,
  $$
  avec la convention que $\textnormal{Inf}\,\emptyset=0$. 
  Le plongement pr{\'e}c{\'e}dent est valu{\'e}. 
  On montre que $A\{\{ t_1,\ldots,t_n \}\}$ est la compl{\'e}tion
  formelle de $A\{ t_1,\ldots,t_n \}$.

  Notons les r{\'e}sultats suivants~:
  
  \smallskip
  \textnormal{\bf Fait 2:}\quad
  Soit $f\in A\{\{ t_1,\ldots,t_n \}\},\quad v(f)\geq 1$. Alors,
  $1-f$ est inversible, et $(1-f)^{-1}=
  1+f+\ldots+f^n+\ldots$
  \smallskip

  La s{\'e}rie est bien convergente, car $\forall n,\, v(f^n)\geq n$. 

  \smallskip
  \textnormal{\bf Fait 3:}\quad
  $\forall s\in S\,\, \exists\, !\, \lambda_s\in A\,\exists\, !\, n_s=:
  \underline{\a}\in N/\, s=\lambda_s t_1^{\a_1}\ldots t_n^{\a_n} \pmod I$. 
  \smallskip

  \textnormal{\bf Fait 4:}\quad
  $\forall\underline{\a}\in N,\, \lbrace s\in S,\, n_s=\underline{\a}
  \rbrace$ est fini.
  \smallskip

  \textnormal{\bf Fait 5:}\quad
  L'application~:
   $$
   \begin{array}{rcl}
   N&\longrightarrow&S\\
   \underline{\a}&\longmapsto&(
   \underset{\a_1\textnormal{ fois}}{\underbrace{1,\ldots,1}}
   ,\ldots,
   \underset{\a_i\textnormal{ fois}}{\underbrace{i,\ldots,i}}
   ,\ldots,
   \underset{\a_n\textnormal{ fois}}{\underbrace{n,\ldots,n}})
   \end{array}
   $$
  est injective.
  \smallskip

  \textnormal{\bf Fait 6:}\quad
  $\forall\, \underline{\a}\in N,\, 
  \bigl\lbrace (\underline{\b},\underline{\gamma})\in N^2\, /\,
  \underline{\b +\gamma}=\underline{\a}
  \bigr\rbrace
  $ est fini.
  \smallskip

  Ceux-ci sont utiles pour le point suivant.

  \medskip
  \noindent
  $\bullet$ 
  L'alg{\`e}bre $A\{\{ t_1,\ldots,t_n \}\}/I$ est par d{\'e}finition
  le quotient de l'alg{\`e}bre\break
  $A\{\{ t_1,\ldots,t_n \}\}$ par l'id{\'e}al engendr{\'e}
  par $I$ dans $A\{\{ t_1,\ldots,t_n \}\}$, via le plongement~:
   $$
    A\{ t_1,\ldots,t_n \}\hookrightarrow
    A\{\{ t_1,\ldots,t_n \}\}.
   $$
  L'id{\'e}al
  engendr{\'e} par $I$ est ferm{\'e} dans $A\{\{ t_1,\ldots,t_n \}\}$.
  Ceci nous permet de d{\'e}finir une valuation quotiente (not{\'e}e
  encore $v$) par~:
  $$
  \forall x\in A\{\{ t_1,\ldots,t_n \}\},\quad
  v(\bar{x})=\textnormal{sup }\bigl\lbrace v(x+a),\, a\in I
  \bigr\rbrace.
  $$
  On montre que $A\{\{ t_1,\ldots,t_n \}\}/I$ est la compl{\'e}tion
  formelle de l'alg{\`e}bre\break
  $A\{ t_1,\ldots,t_n \}/I$.
  \smallskip
  
  Par ailleurs, consid{\'e}rons ${\cal F}(N,A)$, le $A$-module
  des fonctions de $N$ dans $A$. Si $f$ et $g$ sont dans 
  ${\cal F}(N,A)$, on d{\'e}finit le produit de $f$ avec $g$ 
  comme {\'e}tant~:
   $$
   \begin{array}{rcl}
   f\, .\, g:\quad N&\longrightarrow&A\\
   \underline{\a}&\longmapsto&
   \displaystyle\sum\limits_{{\underline \beta},
   {\underline \gamma}
   \atop{
   {\underline {\beta+\gamma}}={\underline \alpha}}}
   q^{-(\b_{2}\gamma_{1}+\ldots+\b_n\gamma_{n-1})}
   f({\underline \beta})\, g({\underline \gamma}).
   \end{array}  
   $$
  L'expression a bien un sens d'apr{\`e}s le fait 6.
  Muni de ce produit, ${\cal F}(N,A)$ est une $A$-alg{\`e}bre.

  D'apr{\`e}s le fait 5, $N$ se voit comme une partie de $S$. 
  Notons $\chi_{N}$, la fonction caract{\'e}ristique de $N$ dans $S$.
  Alors, l'application~:
   $$
    \begin{array}{rcl}
    {\cal F}(N,A)&\longrightarrow&A\{\{ t_1,\ldots,t_n \}\}/I\\
    f&\longmapsto&\textnormal{cl}(g)
    \end{array}
   $$
   o{\`u} $g$ est la fonction~:
       $$
        \begin{array}{rcl}
             g:\quad S&\longrightarrow&A\\
                     s&\longmapsto&\chi_{N}(s)f(s)
        \end{array}
       $$
  est un isomorphisme d'alg{\`e}bres. Son inverse est l'application~:
  $$
    \begin{array}{rcl}
       A\{\{ t_1,\ldots,t_n \}\}/I&\longrightarrow&{\cal F}(N,A)\\
       \textnormal{cl}(f)&\longmapsto&g
    \end{array}
   $$
  o{\`u} $g$ est la fonction~:
  $$
   \begin{array}{rcl}
   g:\quad N&\longrightarrow&A\\
   {\underline \a}&\longmapsto&
   \displaystyle\sum\limits_{n_s={\underline \a}}\lambda_{s}f(s).
   \end{array}
  $$
  D'apr{\`e}s le fait 4, ceci a bien un sens. Ainsi, 
  $A\{\{ t_1,\ldots,t_n \}\}/I$ s'identifie {\`a} ${\cal F}(N,A)$
  comme alg{\`e}bre.
 
\newpage
\section{Comultiplication tordue}
\noindent
Nous rappelons la d{\'e}finition du produit tensoriel
tordu d'alg{\`e}bres gradu{\'e}es. Dans toute la suite, $k$ d{\'e}signera un
anneau commutatif unitaire et $k^*$ l'ensemble de ses {\'e}l{\'e}ments
inversibles.

\medskip
\noindent
{\bf Rappel:}

\medskip
\noindent
{\it Soient $A_1$, $A_2$, deux $k-$alg{\`e}bres gradu{\'e}es et $q\in k^*$.
Notons $i_1$ (resp. $i_2$) l'injection canonique de $A_1$
(resp. $A_2$) dans $A_1\times A_2$ d{\'e}fini par
$\forall a_1\in A_1$ (resp. $\forall a_2\in A_2$),
$i_1(a_1)=(a_1,1)$ (resp. $i_2(a_2)=(1,a_2)$).
Il existe {\`a} isomorphisme gradu{\'e} pr{\`e}s, une et une seule
alg{\`e}bre gradu{\'e}e $(A,\delta)$ et une unique application
bilin{\'e}aire
$f:\, A_1\times A_2\longrightarrow A$ qui v{\'e}rifient la propri{\'e}t{\'e}
universelle suivante~:\\
Soient $B$ une $k-$alg{\`e}bre gradu{\'e}e et 
$\varphi:A_1\times A_2\longrightarrow B$ une application bilin{\'e}aire
balanc{\'e}e
(i.e., $\forall\lambda\in k\,\forall (a_1,a_2)\in A_1\times A_2,\,
\varphi (\lambda a_1,a_2)=\varphi (a_1,\lambda a_2)=
\lambda\varphi (a_1,a_2)$)
qui satisfait les deux conditions suivantes~:\\
(1) Pour tout $k\in\lbrace 1,2\rbrace$, l'application
$\varphi\circ i_k$ est un morphisme d'alg{\`e}bres gradu{\'e}es
de $A_k$ vers $B$.\\
(2) Pour tout $(a_1,a_2)\in A_1\times A_2$, avec $a_i$ homog{\`e}ne
de degr{\'e} $\alpha_i$ ($i\in\lbrace 1,2\rbrace$),
on a 
$$
\begin{array}{rl}
&(\varphi\circ i_1)(a_1)(\varphi\circ i_2)(a_2)=\varphi(a_1,a_2)\\
{\text et\quad }&
(\varphi\circ i_2)(a_2)(\varphi\circ i_1)(a_1)=q^{\alpha_1\alpha_2}
\varphi(a_1,a_2)
\end{array}
$$
Alors, il existe un unique morphisme gradu{\'e} 
$\widehat{\varphi}:A\longrightarrow B$ tel que
le diagramme suivant soit commutatif~:
$$
\begin{array}{rcll}
A&\overset{\widehat{\varphi}}{\longrightarrow}&B\\
f\searrow&&\nearrow\varphi\\
&A_1\times A_2&
\end{array}
$$}

\smallskip
\noindent
Comme espace vectoriel, $A$ est isomorphe
{\`a} $A_1\otimes_k A_2$.
La structure d'alg{\`e}bre sur $A$ se d{\'e}finit par~:
pour tous couples $(a_1,a_2)$ et $(a'_1,a'_2)$ dans $A_1\times A_2$
constitu{\'e}s d'{\'e}l{\'e}ments homog{\`e}nes, on pose~:
$$
(a_1\otimes a_2).(a'_1\otimes a'_2)=
q^{(\partial^{\circ}a'_1)(\partial^{\circ}a_2)}a_1 a'_1\otimes a_2 a'_2.
$$
La graduation sur $A$ est alors donn{\'e}e par l'endomorphisme
diagonalisable~:
$\delta=\delta_1\otimes Id_{A_1}+Id_{A_2}\otimes\delta_2$,
$\delta_i$ d{\'e}signant la d{\'e}rivation donnant la graduation sur
$A_i$.\\
On notera $A_1\otimes_q A_2$ le produit tensoriel $q-$tordu de $A_1$
et de $A_2$, et $A_1\ootimes A_2$ le produit tensoriel
$q^{-1}-$tordu de $A_1$ et $A_2$, c'est {\`a} dire 
$A_1\otimes_{q^{-1}}A_2$.\medskip\\
Par associativit{\'e}, on montre {\'e}galement que l'on peut d{\'e}finir
le produit tensoriel 
$q-$tordu de plusieurs alg{\`e}bres gradu{\'e}es.
Si $n\in{\NN}^*$, si $A_i,\, i\in\lbrace 1,\ldots n\rbrace$
est une alg{\`e}bre gradu{\'e}e, et si 
$a_i\in A_i$ est homog{\`e}ne de degr{\'e} $\alpha_i$, alors
$a_1\ootimes\ldots\ootimes a_n$ est homog{\`e}ne de degr{\'e} 
$\alpha_1+\ldots+\alpha_n$.

\smallskip
\noindent
Soient maintenant $A,B,C,D$ quatre alg{\`e}bres gradu{\'e}es,
$g:C\times D\longrightarrow C\otimes_q D$ l'application 
bilin{\'e}aire naturelle,
et $u:A\longrightarrow C$ et $v:B\longrightarrow D$ deux morphismes
gradu{\'e}s.
alors, $g\circ(u\times v):A\times B\longrightarrow C\otimes_q D$
est une application bilin{\'e}aire satisfaisant les deux propri{\'e}t{\'e}s
requises de la proposition pr{\'e}c{\'e}dente.
Par la propri{\'e}t{\'e} universelle du produit tensoriel $q-$tordu,
on en d{\'e}duit l'existence d'un morphisme gradu{\'e} not{\'e}
$$
u\otimes_q v:A\otimes_q B\longrightarrow C\otimes_q D
$$
tel que 
$$
\forall (a,b)\in A\times B,\, 
(u\otimes_q v)(a\otimes_q b)=u(a)\otimes_q v(b).
$$
Par associativit{\'e}, on peut {\'e}galement d{\'e}finir le produit
tensoriel $q-$tordu de plusieurs morphismes gradu{\'e}s.

\medskip
\noindent
{\bf Exemples:}

\medskip
\noindent
$\bullet$
Pour $i\in\lbrace 1,\ldots,n\rbrace$, on d{\'e}finit
une graduation sur $k[x_i,x_i^{-1}]$ et $k[y_i,y_i^{-1}]$, en d{\'e}cr{\'e}tant
que $\deg{x_i}=1$ et $\deg{y_i}=-1$).

\noindent
Pour $i\in\lbrace 1,\ldots,n\rbrace$, on a deux morphismes gradu{\'e}s~:
$$
\begin{array}{rcl}
u_i:\quad k[x_i,x_i^{-1}]&\longrightarrow&
k[x_1^{\pm 1},\ldots,y_n^{\pm 1}]_q\\
x_i^{\pm 1}&\longmapsto&x_i^{\pm 1}
\end{array}
$$
et~:
$$
\begin{array}{rcl}
v_i:\quad k[y_i,y_i^{- 1}]&
\longrightarrow&k[x_1^{\pm 1},\ldots,y_n^{\pm 1}]_q\\
y_i^{\pm 1}&\longmapsto&y_i^{\pm 1}
\end{array}
$$
\noindent
$\bullet$ Sur l'ensemble 
$E:=\bigl\lbrace u_i,v_i\, ;i,j\in\lbrace 1,\ldots n\rbrace\bigr\rbrace$,
on a une relation d'ordre total d{\'e}fini par:
Pour $(\a,\b)\in E^2,\, \a\prec\b$ si et seulement si on est dans l'un
des quatres cas suivants~:
\begin{enumerate}
\item $\a=u_i,\b=u_j$ et $i<j$
\item $\a=u_i,\b=v_j$ et $i\leq j$
\item $\a=v_i,\b=u_j$ et $i<j$
\item $\a=v_i,\b=v_j$ et $i<j$.
\end{enumerate}
On montre facilement que si $(\a,\b)\in E^2$, alors
$$
\a\prec\b\Longrightarrow \forall x,y,\, 
\b(y)\a(x)=q^{-(\deg{x})(\deg{y})}\a(x)\b(y).
$$
Par suite, il existe un morphisme gradu{\'e}~:
\begin{equation}
\begin{array}{rcl}
\phi:\quad k[x_1,x_1^{-1}]\ootimes\ldots\ootimes
k[y_n,y_n^{-1}]&\longrightarrow&
k[x_1^{\pm 1},\ldots,y_n^{\pm 1}]_q\\
u_1\ootimes\ldots\ootimes v_n&
\longmapsto&u_1\ldots v_n
\end{array}
\end{equation}
$\bullet$
Il existe sur $k\lbrace
x_1^{\pm 1},\ldots,y_n^{\pm 1}\rbrace$,
une graduation telle que~:
$$
\deg{x_i^{\pm 1}}=\pm 1\quad \textnormal{ et }\quad 
\deg{y_i^{\pm 1}}=\mp 1.
$$
Le morphisme ~:
$$
\begin{array}{rcl}
\widehat{\psi}:\quad k\lbrace x_1^{\pm 1},\ldots,y_n^{\pm 1}\rbrace&
\longrightarrow&k[x_1,x_1^{-1}]\ootimes\ldots\ootimes k[y_n,y_n^{-1}]\\
x_i^{\pm 1}&\longmapsto&1\ootimes\ldots\ootimes x_i^{\pm 1}\ootimes
\ldots\ootimes 1\\
y_j^{\pm 1}&\longmapsto&1\ootimes\ldots\ootimes y_j^{\pm 1}
\ootimes\ldots\ootimes 1
\end{array}
$$
est gradu{\'e} et se factorise en un morphisme gradu{\'e}~:
$$
\begin{array}{rcl}
\psi:\quad k[x_1^{\pm 1},\ldots,y_n^{\pm 1}]_q&
\longrightarrow&
k[x_1,x_1^{-1}]\ootimes\ldots\ootimes k[y_n,y_n^{-1}]\\
x_i^{\pm 1}&\longmapsto&1\ootimes\ldots\ootimes x_i^{\pm 1}
\ootimes\ldots\ootimes 1\\
y_j^{\pm 1}&\longmapsto&1\ootimes\ldots\ootimes y_j^{\pm 1}
\ootimes\ldots\ootimes 1
\end{array}
$$
$\bullet$
Il est facile de voir $\phi$ et $\psi$ sont inverses
l'une de l'autre. Par suite, $\phi$ et $\psi$ sont des 
isomporphismes gradu{\'e}s.

\medskip
\noindent
Revenons {\`a} l'alg{\`e}bre ${\cal A}_q$ de \ref{defaq}.

\medskip
\noindent
$\bullet$
On pose~:
${\cal L}(\lambda)=[a_{i,j}(\lambda)]\in M_2({\cal A}_q)$.
On a~: 
\begin{equation}\label{comt1}
R({\lambda\over \mu}){\cal L}^{1}(\lambda)H{\cal L}^{2}(\mu)=
{\cal L}^{2}(\mu)H{\cal L}^{1}(\lambda)R({\lambda\over \mu}).
\end{equation}
Les deux matrices ${\cal L}^{1}(\lambda)$ et ${\cal L}^{2}(\mu)$ ne commutent
pas dans $M_4({\cal A}_q)$.

\smallskip
\noindent
$\bullet$
L'alg{\`e}bre ${\cal A}_q$ {\'e}tant gradu{\'e}e, consid{\'e}rons les deux morphismes
d'alg{\`e}bres suivants~:
$$
\begin{array}{rclcrcl}
g:\quad {\cal A}_q&\longrightarrow&{\cal A}_q\ootimes {\cal A}_q&
\textnormal{ et\quad}&d:\quad
{\cal A}_q&\longrightarrow&{\cal A}_q\ootimes {\cal A}_q\\
x&\longrightarrow&x\ootimes 1&&x&\longrightarrow&1\ootimes x\\
\end{array}
$$
Pour tout $n\in{\NN}^{*}$, l'application $g$ (resp. $d$)
se prolonge en un morphisme d'alg{\`e}bres de $M_n({\cal A}_q)$
dans $M_n({\cal A}_q\ootimes {\cal A}_q)$. On notera $M^g$ (resp. $M^d$)
l'image de $M$ par ce morphisme.

\smallskip
\noindent
$\bullet$
Il est clair que
$$
\forall i\in\lbrace 1,2\rbrace\,
\forall x\in\lbrace g,d\rbrace\quad
\bigl( {\cal L}^i(\lambda)\bigr)^x=\bigl( {\cal L}^x(\lambda)\bigr)^{i}.
$$
Par suite, on notera ${\cal L}^{i,x}(\lambda)$ pour 
$\bigl( {\cal L}^i(\lambda)\bigr)^x$.

\noindent
La relation \ref{comt1} entra{\^\i}ne imm{\'e}diatement les 
deux relations suivantes~:
\begin{align}
\label{comt2}
R({\lambda\over \mu}){\cal L}^{1g}(\lambda)H{\cal L}^{2g}(\mu)&=
{\cal L}^{2g}(\mu)H{\cal L}^{1g}(\lambda)R({\lambda\over \mu})\\
\label{comt3}
R({\lambda\over \mu}){\cal L}^{1d}(\lambda)H{\cal L}^{2d}(\mu)&=
{\cal L}^{2d}(\mu)H{\cal L}^{1d}(\lambda)R({\lambda\over \mu})
\end{align}
Les matrices ${\cal L}^{1g}(\lambda)$ et ${\cal L}^{2d}(\mu)$ ne commutent pas dans 
$M_4({\cal A}_q)$, pas plus que ${\cal L}^{1d}(\lambda)$ et ${\cal L}^{2g}(\mu)$.
Cependant, on a le lemme suivant~:
\begin{lem}
Les relations suivantes sont v{\'e}rifi{\'e}es dans $M_4({\cal A}_q)$~:
\begin{align}
\label{comt4}
Ad(H)\bigl( {\cal L}^{1g}(\lambda)\bigr){\cal L}^{2d}(\mu)&=
{\cal L}^{2d}(\mu)Ad(H)\bigl( {\cal L}^{1g}(\lambda)\bigr)\\
\label{comt5}
Ad(H)\bigl( {\cal L}^{2g}(\lambda)\bigr){\cal L}^{1d}(\mu)&=
{\cal L}^{1d}(\mu)Ad(H)\bigl( {\cal L}^{2g}(\lambda)\bigr)\\
\label{comt6}
Ad(H)\bigl( R({\lambda\over \mu})\bigr)&=R({\lambda\over \mu})
\end{align}
\end{lem}
\begin{dem}
La derni{\`e}re {\'e}galit{\'e} est claire d'apr{\`e}s (\ref{mat1})
et (\ref{mat2}). Les deux autres relations se v{\'e}rifient
matrici{\`e}lement. Par exemple, 
\begin{multline*}
Ad(H)\bigl( {\cal L}^{2g}(\lambda)\bigr)\\
=\begin{pmatrix}
a_{1,1}(\lambda)\ootimes 1&q^{-{1\over 2}}a_{1,2}(\lambda)\ootimes 1&
0&0\\
q^{{1\over 2}}a_{2,1}(\lambda)\ootimes 1&a_{2,2}(\lambda)\ootimes 1&
0&0\\0&0&
a_{1,1}(\lambda)\ootimes 1&q^{{1\over 2}}a_{1,2}(\lambda)\ootimes 1\\
0&0&
q^{-{1\over 2}}a_{2,1}(\lambda)\ootimes 1&a_{2,2}(\lambda)\ootimes 1
\end{pmatrix}
\end{multline*}
et
$$
{\cal L}^{1d}(\mu)=\begin{pmatrix}
1\ootimes a_{1,1}(\mu)&0&1\ootimes a_{1,2}(\mu)&0\\
0&1\ootimes a_{1,1}(\mu)&0&1\ootimes a_{1,2}(\mu)\\
1\ootimes a_{2,1}(\mu)&0&1\ootimes a_{2,2}(\mu)&0\\
0&1\ootimes a_{2,1}(\mu)&0&1\ootimes a_{2,2}(\mu)
\end{pmatrix}.
$$
Le fait que ces deux matrices commutent est une cons{\'e}quence
de l'{\'e}galit{\'e}~:
$$\bigl( 1\ootimes a_{i,j}(\mu)\bigr)
\bigl( a_{k,l}(\lambda)\ootimes 1\bigr)=
q^{-(i-j)(k-l)}\bigl( a_{k,l}(\lambda)\ootimes 1\bigr)
\bigl( 1\ootimes a_{i,j}(\mu)\bigr),$$
pour $i,j,k,l\in\lbrace 1,2\rbrace$.
\end{dem}

\noindent
Nous pouvons {\`a} pr{\'e}sent en d{\'e}duire l'existence d'une 
comultiplication tordue sur ${\cal A}_q$.
\begin{prop}
Soit $\DDelta$ l'application suivante~:
$$
\begin{array}{rcl}
\DDelta:\quad {\cal A}&\longrightarrow&{\cal A}_q\ootimes {\cal A}_q\\
{\cal L}({\lambda})&\longrightarrow&{\cal L}({\lambda})\ootimes {\cal L}({\lambda})=
{\cal L}^{g}({\lambda}){\cal L}^{d}({\lambda})
\end{array}.
$$
Alors, $\DDelta$ d{\'e}finit un morphisme d'alg{\`e}bres gradu{\'e}es
tel que
$I_q\subset\ker\DDelta$.
D'o{\`u} l'on en d{\'e}duit un morphisme d'alg{\`e}bres gradu{\'e}es
not{\'e} encore
$\DDelta$~:
$$
\begin{array}{rcl}
\DDelta:\quad {\cal A}_q&\longrightarrow&{\cal A}_q\ootimes {\cal A}_q\\
{\cal L}({\lambda})&\longrightarrow&{\cal L}({\lambda})\ootimes 
{\cal L}({\lambda})
\end{array}
$$
\end{prop}
\begin{dem}
Le fait que $\DDelta$ soit gradu{\'e} est clair. Pour montrer qu'il
d{\'e}finit
un morphisme d'alg{\`e}bres de ${\cal A}_q$ dans
${\cal A}_q\ootimes {\cal A}_q$, il faut voir que 
$$
\DDelta\Bigl(
R({\lambda\over \mu}){\cal L}^{1}(\lambda)H{\cal L}^{2}(\mu)\Bigr)
=\DDelta\Bigl(
{\cal L}^{2}(\mu)H{\cal L}^{1}(\lambda)R({\lambda\over \mu})\Bigr).
$$
L'application $\DDelta$ {\'e}tant un morphisme d'alg{\`e}bres, on a~:
\begin{align*}
\DDelta\Bigl(R({\lambda\over \mu})
{\cal L}^{1}(\lambda)H{\cal L}^{2}(\mu)\Bigr)&=
R({\lambda\over \mu})\DDelta\bigl( {\cal L}^{1}(\lambda)\bigr)
H\DDelta\bigl( {\cal L}^{2}(\mu)\bigr)\\
&=R({\lambda\over \mu}){\cal L}^{1g}(\lambda){\cal L}^{1d}(\lambda)
H{\cal L}^{2g}(\mu){\cal L}^{2d}(\mu)\\
&=R({\lambda\over \mu}){\cal L}^{1g}(\lambda)H H^{-1}{\cal L}^{1d}(\lambda)
H{\cal L}^{2g}(\mu){\cal L}^{2d}(\mu)\\
&=R({\lambda\over \mu}){\cal L}^{1g}(\lambda)H\Bigl(
Ad(H^{-1})\bigl( {\cal L}^{1d}(\lambda)\bigr) 
{\cal L}^{2g}(\mu)\Bigr){\cal L}^{2d}(\mu)\\
\intertext{Soit, d'apr{\`e}s (\ref{comt5}),}
\DDelta\Bigl(R({\lambda\over \mu})
{\cal L}^{1}(\lambda)H{\cal L}^{2}(\mu)\Bigr)
&=R({\lambda\over \mu}){\cal L}^{1g}(\lambda)H\Bigl(
{\cal L}^{2g}(\mu)Ad(H^{-1})\bigl( {\cal L}^{1d}(\lambda)\bigr)
\Bigr){\cal L}^{2d}(\mu)\\
&=\bigl( R({\lambda\over \mu}) {\cal L}^{1g}(\lambda) H
{\cal L}^{2g}(\mu)\bigr) H^{-1}\bigl( 
{\cal L}^{1d}(\lambda) H {\cal L}^{2d}(\mu)\bigr)
\end{align*}
Puis, d'apr{\`e}s (\ref{comt2}) et (\ref{comt6}),
\begin{align*}
\DDelta\bigl(R({\lambda\over \mu})
{\cal L}^{1}(\lambda)H{\cal L}^{2}(\mu)\bigr)&=
\bigl( {\cal L}^{2g}(\mu)H{\cal L}^{1g}(\lambda)R({\lambda\over \mu})
\bigr) H^{-1}
{\cal L}^{1d}(\lambda) H {\cal L}^{2d}(\mu)\\
&={\cal L}^{2g}(\mu)H{\cal L}^{1g}(\lambda)H^{-1}\bigl( R({\lambda\over \mu})
{\cal L}^{1d}(\lambda) H {\cal L}^{2d}(\mu)\bigr)
\end{align*}
En utilisant ensuite (\ref{comt3}) et (\ref{comt4}),
on obtient~:
\begin{align*}
\DDelta\Bigl(R({\lambda\over \mu})
{\cal L}^{1}(\lambda)H{\cal L}^{2}(\mu)\Bigr)&=
{\cal L}^{2g}(\mu)H{\cal L}^{1g}(\lambda)H^{-1}\Bigl(
{\cal L}^{2d}(\mu) H {\cal L}^{1d}(\lambda) R({\lambda\over \mu})\Bigr)\\
&={\cal L}^{2g}(\mu)\bigl( HL^{1g}(\lambda)H^{-1}{\cal L}^{2d}(\mu)\bigr)
H {\cal L}^{1d}(\lambda) R({\lambda\over \mu})\\
&={\cal L}^{2g}(\mu)\Bigl( Ad(H)\bigl( {\cal L}^{1g}(\lambda)\bigr)
{\cal L}^{2d}(\mu)\Bigr) H {\cal L}^{1d}(\lambda) R({\lambda\over \mu})\\
&={\cal L}^{2g}(\mu)\Bigl( {\cal L}^{2d}(\mu)
Ad(H)\bigl( {\cal L}^{1g}(\lambda)\bigr)\Bigr) 
H {\cal L}^{1d}(\lambda) R({\lambda\over \mu})\\
&={\cal L}^{2g}(\mu)\bigl(
{\cal L}^{2d}(\mu) H {\cal L}^{1g}(\lambda) H^{-1}\bigr)
H {\cal L}^{1d}(\lambda) R({\lambda\over \mu})\\
&={\cal L}^{2g}(\mu){\cal L}^{2d}(\mu)H {\cal L}^{1g}(\lambda)
{\cal L}^{1d}(\lambda) R({\lambda\over \mu})\\
&=\DDelta\bigl( {\cal L}^{2}(\mu)\bigr)H
\DDelta\bigl( {\cal L}^{1}(\lambda)\bigr)
R({\lambda\over \mu})\\
&=\DDelta\Bigl( {\cal L}^{2}(\mu)H{\cal L}^{1}(\lambda)
R({\lambda\over \mu})\Bigr).
\end{align*}
\end{dem}
\newpage
\noindent
{\centerline{{\bf {\huge R{\'E}F{\'E}RENCES}}}}

\bigskip
\noindent
[1] P. Bouwknegt, J. McCarthy, K. Pilch,
\emph{Quantum group structure in the Fock space resolutions of
$\widehat{\textnormal{sl}_n}$ representations},
Comm. Math. Phys. {\bf 131}\break (1990),
${\textnormal{n}}^{\circ}1$, 125-155.

\smallskip
\noindent
[2] B. Enriquez, B. Feigin, \emph{Integrals of motion of classical lattice
sine-Gordon system}, Theor.Math.Phys. {\bf 103} (1995) 738-756.

\smallskip
\noindent
[3] A.G. Izergin, V.E. Korepin, \emph{The lattice sine-Gordon model}
Lett.Math.\break Phys. {\bf 5} (1981),199-205.

\smallskip
\noindent
[4] A.G. Izergin, V.E. Korepin, \emph{Lattice versions of quantum
field theory models in two dimensions}
Nucl.Phys. B {\bf 205} (1985),401-413.

\smallskip
\noindent
[5] L. Faddeev, A. Volkov, \emph{Quantum inverse scattering method on a
space-time lattice}, Theor.Math.Phys. {\bf 92} (1992), 207-214.

\smallskip
\noindent
[6] A. Bobenko, N. Kutz, U. Pinkall, \emph{The discrete quantum
pendulum}, Physics Letters A {\bf 177} (1993), 399-404.

\smallskip
\noindent
[7] N. Reshetikhin, \emph{Integrable Systems with Discrete Time},
International Conference of mathematical physics, 1994, Paris.

\smallskip
\noindent
[8] V. Bazhanov, A. Bobenko, N. Reshetikhin, \emph{Quantum Discrete
Sine-\break
Gordon Model at roots of 1: Integrable Quantum System on the
Integrable Classical Background},
Comm. Math. Phys. {\bf 175} (1996), 
${\textnormal{n}}^{\circ}2$,377-400.

\smallskip
\noindent
[9] C. Kassel, \emph{Quantum groups}, Springer-Verlag 1994 p.15.

\smallskip
\noindent
[10] V. Drinfeld, \emph{A new realizations of Yangians and quantized
affine algebras}, Soviet. Math. Dokl., {\bf 36}~:2
(1988), 212-216.

\smallskip
\noindent
[11] C. DeConcini, V.G. Kac, C. Procesi,
\emph{Some quantum analogues of solvable Lie groups},
Oxford University Press.Stud.Math., Tata Inst.
Fundam. Res. 13,41-65 (1995).

\smallskip
\noindent
[12] M. Semenov-Tian-Shansky, \emph{Dressing action 
transformations and\, \break 
Poisson-Lie group actions}, Publ. Math. RIMS {\bf 21} (1985),
1237-1260.

\smallskip
\noindent
[13] E. Frenkel, N. Reshetikhin, \emph{Quantum Affine Algebras and Deformations
of the Virasoro and W-algebras}, Comm. Math. Phys.
{\bf 178} (1996), ${\textnormal{n}}^{\circ}1$, 237-264.
 
\smallskip
\noindent
[14] J. Dixmier, \emph{Alg{\`e}bres enveloppantes}, Gauthier-Villards 1974
p.116.

\bigskip
\noindent
\begin{flushleft}
Cyril Grunspan\\
Ecole Polytechnique\\
Centre de Math{\'e}matiques\\
UMR 7640 du CNRS\\
F-91128 PALAISEAU Cedex\\
\begin{verbatim}
courriel: grunspan@math.polytechnique.fr 
\end{verbatim}
\end{flushleft}
\end{document}